\documentclass[12pt,leqno]{amsart}

\usepackage{amssymb}
\usepackage{amsmath}
\usepackage{amsthm}

\usepackage{enumerate}
\usepackage{array}

\begingroup\makeatletter\ifx\SetFigFont\undefined
\def\x#1#2#3#4#5#6#7\relax{\def\x{#1#2#3#4#5#6}}%
\expandafter\x\fmtname xxxxxx\relax \def\y{splain}%
\ifx\x\y   
\gdef\SetFigFont#1#2#3{%
  \ifnum #1<17\tiny\else \ifnum #1<20\small\else
  \ifnum #1<24\normalsize\else \ifnum #1<29\large\else
  \ifnum #1<34\Large\else \ifnum #1<41\LARGE\else
     \huge\fi\fi\fi\fi\fi\fi
  \csname #3\endcsname}%
\else
\gdef\SetFigFont#1#2#3{\begingroup
  \count@#1\relax \ifnum 25<\count@\count@25\fi
  \def\x{\endgroup\@setsize\SetFigFont{#2pt}}%
  \expandafter\x
    \csname \romannumeral\the\count@ pt\expandafter\endcsname
    \csname @\romannumeral\the\count@ pt\endcsname
  \csname #3\endcsname}%
\fi
\fi\endgroup

\newcommand{\Ciit}[3]{
\begin{picture}(2752,895)(7025,-4136)
\thinlines
\put(7201,-3961){\circle{336}}
\put(8401,-3961){\circle{336}}
\put(9601,-3961){\circle{336}}
\put(7369,-3911){\line( 1, 0){864}}
\put(7369,-4011){\line( 1, 0){864}}
\put(8569,-3911){\line( 1, 0){864}}
\put(8569,-4011){\line( 1, 0){864}}
\put(7141,-3436){\makebox(0,0)[lb]{\smash{\SetFigFont{6}{7.2}{rm}\mbox{$#1$}}}}
\put(8341,-3436){\makebox(0,0)[lb]{\smash{\SetFigFont{6}{7.2}{rm}\mbox{$#2$}}}}
\put(9541,-3436){\makebox(0,0)[lb]{\smash{\SetFigFont{6}{7.2}{rm}\mbox{$#3$}}}}
\end{picture}
}
\begingroup\makeatletter\ifx\SetFigFont\undefined
\def\x#1#2#3#4#5#6#7\relax{\def\x{#1#2#3#4#5#6}}%
\expandafter\x\fmtname xxxxxx\relax \def\y{splain}%
\ifx\x\y   
\gdef\SetFigFont#1#2#3{%
  \ifnum #1<17\tiny\else \ifnum #1<20\small\else
  \ifnum #1<24\normalsize\else \ifnum #1<29\large\else
  \ifnum #1<34\Large\else \ifnum #1<41\LARGE\else
     \huge\fi\fi\fi\fi\fi\fi
  \csname #3\endcsname}%
\else
\gdef\SetFigFont#1#2#3{\begingroup
  \count@#1\relax \ifnum 25<\count@\count@25\fi
  \def\x{\endgroup\@setsize\SetFigFont{#2pt}}%
  \expandafter\x
    \csname \romannumeral\the\count@ pt\expandafter\endcsname
    \csname @\romannumeral\the\count@ pt\endcsname
  \csname #3\endcsname}%
\fi
\fi\endgroup

\newcommand{\Ciita}[3]{
\begin{picture}(2752,895)(7025,-4136)
\thinlines
\put(7201,-3961){\circle{336}}
\put(8401,-3961){\circle{336}}
\put(9601,-3961){\circle{336}}
\put(7369,-3911){\line( 1, 0){764}}
\put(7369,-4011){\line( 1, 0){764}}
\put(7800,-4131){$>$}
\put(8569,-3911){\line( 1, 0){764}}
\put(8569,-4011){\line( 1, 0){764}}
\put(9000,-4131){$>$}
\put(7141,-3436){\makebox(0,0)[lb]{\smash{\SetFigFont{6}{7.2}{rm}\mbox{$#1$}}}}
\put(8341,-3436){\makebox(0,0)[lb]{\smash{\SetFigFont{6}{7.2}{rm}\mbox{$#2$}}}}
\put(9541,-3436){\makebox(0,0)[lb]{\smash{\SetFigFont{6}{7.2}{rm}\mbox{$#3$}}}}
\end{picture}
}
\begingroup\makeatletter\ifx\SetFigFont\undefined
\def\x#1#2#3#4#5#6#7\relax{\def\x{#1#2#3#4#5#6}}%
\expandafter\x\fmtname xxxxxx\relax \def\y{splain}%
\ifx\x\y   
\gdef\SetFigFont#1#2#3{%
  \ifnum #1<17\tiny\else \ifnum #1<20\small\else
  \ifnum #1<24\normalsize\else \ifnum #1<29\large\else
  \ifnum #1<34\Large\else \ifnum #1<41\LARGE\else
     \huge\fi\fi\fi\fi\fi\fi
  \csname #3\endcsname}%
\else
\gdef\SetFigFont#1#2#3{\begingroup
  \count@#1\relax \ifnum 25<\count@\count@25\fi
  \def\x{\endgroup\@setsize\SetFigFont{#2pt}}%
  \expandafter\x
    \csname \romannumeral\the\count@ pt\expandafter\endcsname
    \csname @\romannumeral\the\count@ pt\endcsname
  \csname #3\endcsname}%
\fi
\fi\endgroup

\newcommand{\Ciitar}[3]{
\begin{picture}(2752,895)(7025,-4136)
\thinlines
\put(7201,-3961){\circle{236}}
\put(7201,-3961){\circle{436}}
\put(8401,-3961){\circle{336}}
\put(9601,-3961){\circle{336}}
\put(7519,-3911){\line( 1, 0){714}}
\put(7519,-4011){\line( 1, 0){714}}
\put(7346,-4131){$<$}
\put(8569,-3911){\line( 1, 0){764}}
\put(8569,-4011){\line( 1, 0){764}}
\put(9000,-4131){$>$}
\put(7141,-3436){\makebox(0,0)[lb]{\smash{\SetFigFont{6}{7.2}{rm}\mbox{$#1$}}}}
\put(8341,-3436){\makebox(0,0)[lb]{\smash{\SetFigFont{6}{7.2}{rm}\mbox{$#2$}}}}
\put(9541,-3436){\makebox(0,0)[lb]{\smash{\SetFigFont{6}{7.2}{rm}\mbox{$#3$}}}}
\end{picture}
}

\begingroup\makeatletter\ifx\SetFigFont\undefined
\def\x#1#2#3#4#5#6#7\relax{\def\x{#1#2#3#4#5#6}}%
\expandafter\x\fmtname xxxxxx\relax \def\y{splain}%
\ifx\x\y   
\gdef\SetFigFont#1#2#3{%
  \ifnum #1<17\tiny\else \ifnum #1<20\small\else
  \ifnum #1<24\normalsize\else \ifnum #1<29\large\else
  \ifnum #1<34\Large\else \ifnum #1<41\LARGE\else
     \huge\fi\fi\fi\fi\fi\fi
  \csname #3\endcsname}%
\else
\gdef\SetFigFont#1#2#3{\begingroup
  \count@#1\relax \ifnum 25<\count@\count@25\fi
  \def\x{\endgroup\@setsize\SetFigFont{#2pt}}%
  \expandafter\x
    \csname \romannumeral\the\count@ pt\expandafter\endcsname
    \csname @\romannumeral\the\count@ pt\endcsname
  \csname #3\endcsname}%
\fi
\fi\endgroup

\begingroup\makeatletter\ifx\SetFigFont\undefined
\def\x#1#2#3#4#5#6#7\relax{\def\x{#1#2#3#4#5#6}}%
\expandafter\x\fmtname xxxxxx\relax \def\y{splain}%
\ifx\x\y   
\gdef\SetFigFont#1#2#3{%
  \ifnum #1<17\tiny\else \ifnum #1<20\small\else
  \ifnum #1<24\normalsize\else \ifnum #1<29\large\else
  \ifnum #1<34\Large\else \ifnum #1<41\LARGE\else
     \huge\fi\fi\fi\fi\fi\fi
  \csname #3\endcsname}%
\else
\gdef\SetFigFont#1#2#3{\begingroup
  \count@#1\relax \ifnum 25<\count@\count@25\fi
  \def\x{\endgroup\@setsize\SetFigFont{#2pt}}%
  \expandafter\x
    \csname \romannumeral\the\count@ pt\expandafter\endcsname
    \csname @\romannumeral\the\count@ pt\endcsname
  \csname #3\endcsname}%
\fi
\fi\endgroup

\begingroup\makeatletter\ifx\SetFigFont\undefined
\def\x#1#2#3#4#5#6#7\relax{\def\x{#1#2#3#4#5#6}}%
\expandafter\x\fmtname xxxxxx\relax \def\y{splain}%
\ifx\x\y   
\gdef\SetFigFont#1#2#3{%
  \ifnum #1<17\tiny\else \ifnum #1<20\small\else
  \ifnum #1<24\normalsize\else \ifnum #1<29\large\else
  \ifnum #1<34\Large\else \ifnum #1<41\LARGE\else
     \huge\fi\fi\fi\fi\fi\fi
  \csname #3\endcsname}%
\else
\gdef\SetFigFont#1#2#3{\begingroup
  \count@#1\relax \ifnum 25<\count@\count@25\fi
  \def\x{\endgroup\@setsize\SetFigFont{#2pt}}%
  \expandafter\x
    \csname \romannumeral\the\count@ pt\expandafter\endcsname
    \csname @\romannumeral\the\count@ pt\endcsname
  \csname #3\endcsname}%
\fi
\fi\endgroup

\newcommand{\Ciitariii}[3]{
\begin{picture}(2752,895)(7025,-4136)
\thinlines
\put(7201,-3961){\circle{336}}
\put(8401,-3961){\circle{236}}
\put(8401,-3961){\circle{436}}
\put(9601,-3961){\circle{336}}
\put(7369,-3911){\line( 1, 0){714}}
\put(7369,-4011){\line( 1, 0){714}}
\put(7750,-4131){$>$}
\put(8719,-3911){\line( 1, 0){714}}
\put(8719,-4011){\line( 1, 0){714}}
\put(8550,-4131){$<$}
\put(7141,-3436){\makebox(0,0)[lb]{\smash{\SetFigFont{6}{7.2}{rm}\mbox{$#1$}}}}
\put(8341,-3436){\makebox(0,0)[lb]{\smash{\SetFigFont{6}{7.2}{rm}\mbox{$#2$}}}}
\put(9541,-3436){\makebox(0,0)[lb]{\smash{\SetFigFont{6}{7.2}{rm}\mbox{$#3$}}}}
\end{picture}
}

\begingroup\makeatletter\ifx\SetFigFont\undefined
\def\x#1#2#3#4#5#6#7\relax{\def\x{#1#2#3#4#5#6}}%
\expandafter\x\fmtname xxxxxx\relax \def\y{splain}%
\ifx\x\y   
\gdef\SetFigFont#1#2#3{%
  \ifnum #1<17\tiny\else \ifnum #1<20\small\else
  \ifnum #1<24\normalsize\else \ifnum #1<29\large\else
  \ifnum #1<34\Large\else \ifnum #1<41\LARGE\else
     \huge\fi\fi\fi\fi\fi\fi
  \csname #3\endcsname}%
\else
\gdef\SetFigFont#1#2#3{\begingroup
  \count@#1\relax \ifnum 25<\count@\count@25\fi
  \def\x{\endgroup\@setsize\SetFigFont{#2pt}}%
  \expandafter\x
    \csname \romannumeral\the\count@ pt\expandafter\endcsname
    \csname @\romannumeral\the\count@ pt\endcsname
  \csname #3\endcsname}%
\fi
\fi\endgroup

\begingroup\makeatletter\ifx\SetFigFont\undefined
\def\x#1#2#3#4#5#6#7\relax{\def\x{#1#2#3#4#5#6}}%
\expandafter\x\fmtname xxxxxx\relax \def\y{splain}%
\ifx\x\y   
\gdef\SetFigFont#1#2#3{%
  \ifnum #1<17\tiny\else \ifnum #1<20\small\else
  \ifnum #1<24\normalsize\else \ifnum #1<29\large\else
  \ifnum #1<34\Large\else \ifnum #1<41\LARGE\else
     \huge\fi\fi\fi\fi\fi\fi
  \csname #3\endcsname}%
\else
\gdef\SetFigFont#1#2#3{\begingroup
  \count@#1\relax \ifnum 25<\count@\count@25\fi
  \def\x{\endgroup\@setsize\SetFigFont{#2pt}}%
  \expandafter\x
    \csname \romannumeral\the\count@ pt\expandafter\endcsname
    \csname @\romannumeral\the\count@ pt\endcsname
  \csname #3\endcsname}%
\fi
\fi\endgroup

\newcommand{\Fivta}[5]{
\begin{picture}(5152,895)(4625,-4136)
\thinlines
\put(4801,-3961){\circle{336}}
\put(6001,-3961){\circle{336}}
\put(7201,-3961){\circle{336}}
\put(8401,-3961){\circle{336}}
\put(9601,-3961){\circle{336}}
\put(4969,-3961){\line( 1, 0){864}}
\put(6169,-3961){\line( 1, 0){864}}
\put(7369,-3911){\line( 1, 0){764}}
\put(7369,-4011){\line( 1, 0){764}}
\put(8569,-3961){\line( 1, 0){864}}
\put(7800,-4131){$>$}
\put(4741,-3436){\makebox(0,0)[lb]{\smash{\SetFigFont{6}{7.2}{rm}\mbox{$#1$}}}}
\put(5941,-3436){\makebox(0,0)[lb]{\smash{\SetFigFont{6}{7.2}{rm}\mbox{$#2$}}}}
\put(7141,-3436){\makebox(0,0)[lb]{\smash{\SetFigFont{6}{7.2}{rm}\mbox{$#3$}}}}
\put(8341,-3436){\makebox(0,0)[lb]{\smash{\SetFigFont{6}{7.2}{rm}\mbox{$#4$}}}}
\put(9541,-3436){\makebox(0,0)[lb]{\smash{\SetFigFont{6}{7.2}{rm}\mbox{$#5$}}}}
\end{picture}
}
\begingroup\makeatletter\ifx\SetFigFont\undefined
\def\x#1#2#3#4#5#6#7\relax{\def\x{#1#2#3#4#5#6}}%
\expandafter\x\fmtname xxxxxx\relax \def\y{splain}%
\ifx\x\y   
\gdef\SetFigFont#1#2#3{%
  \ifnum #1<17\tiny\else \ifnum #1<20\small\else
  \ifnum #1<24\normalsize\else \ifnum #1<29\large\else
  \ifnum #1<34\Large\else \ifnum #1<41\LARGE\else
     \huge\fi\fi\fi\fi\fi\fi
  \csname #3\endcsname}%
\else
\gdef\SetFigFont#1#2#3{\begingroup
  \count@#1\relax \ifnum 25<\count@\count@25\fi
  \def\x{\endgroup\@setsize\SetFigFont{#2pt}}%
  \expandafter\x
    \csname \romannumeral\the\count@ pt\expandafter\endcsname
    \csname @\romannumeral\the\count@ pt\endcsname
  \csname #3\endcsname}%
\fi
\fi\endgroup

\newcommand{\Fivtaii}[5]{
\begin{picture}(5152,895)(4625,-4136)
\thinlines
\put(4801,-3961){\circle{336}}
\put(6001,-3961){\circle{336}}
\put(7201,-3961){\circle{336}}
\put(8401,-3961){\circle{336}}
\put(9601,-3961){\circle{336}}
\put(4969,-3961){\line( 1, 0){864}}
\put(6169,-3961){\line( 1, 0){864}}
\put(7469,-3911){\line( 1, 0){764}}
\put(7469,-4011){\line( 1, 0){764}}
\put(7296,-4131){$<$}
\put(8569,-3961){\line( 1, 0){864}}
\put(4741,-3436){\makebox(0,0)[lb]{\smash{\SetFigFont{6}{7.2}{rm}\mbox{$#1$}}}}
\put(5941,-3436){\makebox(0,0)[lb]{\smash{\SetFigFont{6}{7.2}{rm}\mbox{$#2$}}}}
\put(7141,-3436){\makebox(0,0)[lb]{\smash{\SetFigFont{6}{7.2}{rm}\mbox{$#3$}}}}
\put(8341,-3436){\makebox(0,0)[lb]{\smash{\SetFigFont{6}{7.2}{rm}\mbox{$#4$}}}}
\put(9541,-3436){\makebox(0,0)[lb]{\smash{\SetFigFont{6}{7.2}{rm}\mbox{$#5$}}}}
\end{picture}
}
\begingroup\makeatletter\ifx\SetFigFont\undefined
\def\x#1#2#3#4#5#6#7\relax{\def\x{#1#2#3#4#5#6}}%
\expandafter\x\fmtname xxxxxx\relax \def\y{splain}%
\ifx\x\y   
\gdef\SetFigFont#1#2#3{%
  \ifnum #1<17\tiny\else \ifnum #1<20\small\else
  \ifnum #1<24\normalsize\else \ifnum #1<29\large\else
  \ifnum #1<34\Large\else \ifnum #1<41\LARGE\else
     \huge\fi\fi\fi\fi\fi\fi
  \csname #3\endcsname}%
\else
\gdef\SetFigFont#1#2#3{\begingroup
  \count@#1\relax \ifnum 25<\count@\count@25\fi
  \def\x{\endgroup\@setsize\SetFigFont{#2pt}}%
  \expandafter\x
    \csname \romannumeral\the\count@ pt\expandafter\endcsname
    \csname @\romannumeral\the\count@ pt\endcsname
  \csname #3\endcsname}%
\fi
\fi\endgroup

\begingroup\makeatletter\ifx\SetFigFont\undefined
\def\x#1#2#3#4#5#6#7\relax{\def\x{#1#2#3#4#5#6}}%
\expandafter\x\fmtname xxxxxx\relax \def\y{splain}%
\ifx\x\y   
\gdef\SetFigFont#1#2#3{%
  \ifnum #1<17\tiny\else \ifnum #1<20\small\else
  \ifnum #1<24\normalsize\else \ifnum #1<29\large\else
  \ifnum #1<34\Large\else \ifnum #1<41\LARGE\else
     \huge\fi\fi\fi\fi\fi\fi
  \csname #3\endcsname}%
\else
\gdef\SetFigFont#1#2#3{\begingroup
  \count@#1\relax \ifnum 25<\count@\count@25\fi
  \def\x{\endgroup\@setsize\SetFigFont{#2pt}}%
  \expandafter\x
    \csname \romannumeral\the\count@ pt\expandafter\endcsname
    \csname @\romannumeral\the\count@ pt\endcsname
  \csname #3\endcsname}%
\fi
\fi\endgroup

\newcommand{\Giita}[3]{
\begin{picture}(2752,895)(1025,-4136)
\thinlines
\put(1201,-3961){\circle{336}}
\put(2401,-3961){\circle{336}}
\put(3601,-3961){\circle{336}}
\put(1369,-3961){\line( 1, 0){864}}
\put(2769,-3861){\line( 1, 0){704}}
\put(2769,-4061){\line( 1, 0){704}}
\put(2569,-3961){\line( 1, 0){864}}
\put(2496,-4131){$<$}
\put(1141,-3436){\makebox(0,0)[lb]{\smash{\SetFigFont{6}{7.2}{rm}\mbox{$#1$}}}}
\put(2341,-3436){\makebox(0,0)[lb]{\smash{\SetFigFont{6}{7.2}{rm}\mbox{$#2$}}}}
\put(3541,-3436){\makebox(0,0)[lb]{\smash{\SetFigFont{6}{7.2}{rm}\mbox{$#3$}}}}
\end{picture}
}

\begingroup\makeatletter\ifx\SetFigFont\undefined
\def\x#1#2#3#4#5#6#7\relax{\def\x{#1#2#3#4#5#6}}%
\expandafter\x\fmtname xxxxxx\relax \def\y{splain}%
\ifx\x\y   
\gdef\SetFigFont#1#2#3{%
  \ifnum #1<17\tiny\else \ifnum #1<20\small\else
  \ifnum #1<24\normalsize\else \ifnum #1<29\large\else
  \ifnum #1<34\Large\else \ifnum #1<41\LARGE\else
     \huge\fi\fi\fi\fi\fi\fi
  \csname #3\endcsname}%
\else
\gdef\SetFigFont#1#2#3{\begingroup
  \count@#1\relax \ifnum 25<\count@\count@25\fi
  \def\x{\endgroup\@setsize\SetFigFont{#2pt}}%
  \expandafter\x
    \csname \romannumeral\the\count@ pt\expandafter\endcsname
    \csname @\romannumeral\the\count@ pt\endcsname
  \csname #3\endcsname}%
\fi
\fi\endgroup

\newcommand{\Giitaii}[3]{
\begin{picture}(2752,895)(1025,-4136)
\thinlines
\put(1201,-3961){\circle{336}}
\put(2401,-3961){\circle{336}}
\put(3601,-3961){\circle{336}}
\put(1369,-3961){\line( 1, 0){864}}
\put(2519,-3861){\line( 1, 0){704}}
\put(2519,-4061){\line( 1, 0){704}}
\put(2569,-3961){\line( 1, 0){864}}
\put(3000,-4131){$>$}
\put(1141,-3436){\makebox(0,0)[lb]{\smash{\SetFigFont{6}{7.2}{rm}\mbox{$#1$}}}}
\put(2341,-3436){\makebox(0,0)[lb]{\smash{\SetFigFont{6}{7.2}{rm}\mbox{$#2$}}}}
\put(3541,-3436){\makebox(0,0)[lb]{\smash{\SetFigFont{6}{7.2}{rm}\mbox{$#3$}}}}
\end{picture}
}

\begingroup\makeatletter\ifx\SetFigFont\undefined
\def\x#1#2#3#4#5#6#7\relax{\def\x{#1#2#3#4#5#6}}%
\expandafter\x\fmtname xxxxxx\relax \def\y{splain}%
\ifx\x\y   
\gdef\SetFigFont#1#2#3{%
  \ifnum #1<17\tiny\else \ifnum #1<20\small\else
  \ifnum #1<24\normalsize\else \ifnum #1<29\large\else
  \ifnum #1<34\Large\else \ifnum #1<41\LARGE\else
     \huge\fi\fi\fi\fi\fi\fi
  \csname #3\endcsname}%
\else
\gdef\SetFigFont#1#2#3{\begingroup
  \count@#1\relax \ifnum 25<\count@\count@25\fi
  \def\x{\endgroup\@setsize\SetFigFont{#2pt}}%
  \expandafter\x
    \csname \romannumeral\the\count@ pt\expandafter\endcsname
    \csname @\romannumeral\the\count@ pt\endcsname
  \csname #3\endcsname}%
\fi
\fi\endgroup

\newcommand{\Evit}[7]{\raisebox{2.4mm}{\parbox[c]{3.4cm}{$
\begin{picture}(5300,3400)(950,-6600)
\thinlines
\put(1201,-3961){\circle{336}}
\put(2401,-3961){\circle{336}}
\put(3601,-3961){\circle{336}}
\put(4801,-3961){\circle{336}}
\put(6001,-3961){\circle{336}}
\put(3601,-5161){\circle{336}}
\put(3601,-6361){\circle{336}}
\put(1369,-3961){\line( 1, 0){864}}
\put(2569,-3961){\line( 1, 0){864}}
\put(3769,-3961){\line( 1, 0){864}}
\put(4969,-3961){\line( 1, 0){864}}
\put(3601,-4129){\line( 0, -1){864}}
\put(3601,-5329){\line( 0, -1){864}}
\put(1141,-3436){\makebox(0,0)[lb]{\smash{\SetFigFont{6}{7.2}{rm}\mbox{$#1$}}}}
\put(2341,-3436){\makebox(0,0)[lb]{\smash{\SetFigFont{6}{7.2}{rm}\mbox{$#3$}}}}
\put(3541,-3436){\makebox(0,0)[lb]{\smash{\SetFigFont{6}{7.2}{rm}\mbox{$#4$}}}}
\put(4741,-3436){\makebox(0,0)[lb]{\smash{\SetFigFont{6}{7.2}{rm}\mbox{$#5$}}}}
\put(5941,-3436){\makebox(0,0)[lb]{\smash{\SetFigFont{6}{7.2}{rm}\mbox{$#6$}}}}
\put(4024,-5263){\makebox(0,0)[lb]{\smash{\SetFigFont{6}{7.2}{rm}\mbox{$#2$}}}}
\put(4024,-6463){\makebox(0,0)[lb]{\smash{\SetFigFont{6}{7.2}{rm}\mbox{$#7$}}}}
\end{picture}$}}
}

\begingroup\makeatletter\ifx\SetFigFont\undefined
\def\x#1#2#3#4#5#6#7\relax{\def\x{#1#2#3#4#5#6}}%
\expandafter\x\fmtname xxxxxx\relax \def\y{splain}%
\ifx\x\y   
\gdef\SetFigFont#1#2#3{%
  \ifnum #1<17\tiny\else \ifnum #1<20\small\else
  \ifnum #1<24\normalsize\else \ifnum #1<29\large\else
  \ifnum #1<34\Large\else \ifnum #1<41\LARGE\else
     \huge\fi\fi\fi\fi\fi\fi
  \csname #3\endcsname}%
\else
\gdef\SetFigFont#1#2#3{\begingroup
  \count@#1\relax \ifnum 25<\count@\count@25\fi
  \def\x{\endgroup\@setsize\SetFigFont{#2pt}}%
  \expandafter\x
    \csname \romannumeral\the\count@ pt\expandafter\endcsname
    \csname @\romannumeral\the\count@ pt\endcsname
  \csname #3\endcsname}%
\fi
\fi\endgroup

\newcommand{\Eviit}[8]{\raisebox{2.4mm}{\parbox[c]{5cm}{$
\begin{picture}(7700,2200)(950,-5400)
\thinlines
\put(1201,-3961){\circle{336}}
\put(2401,-3961){\circle{336}}
\put(3601,-3961){\circle{336}}
\put(4801,-3961){\circle{336}}
\put(6001,-3961){\circle{336}}
\put(7201,-3961){\circle{336}}
\put(8401,-3961){\circle{336}}
\put(3601,-5161){\circle{336}}
\put(1369,-3961){\line( 1, 0){864}}
\put(2569,-3961){\line( 1, 0){864}}
\put(3769,-3961){\line( 1, 0){864}}
\put(4969,-3961){\line( 1, 0){864}}
\put(6169,-3961){\line( 1, 0){864}}
\put(7369,-3961){\line( 1, 0){864}}
\put(3601,-4129){\line( 0, -1){864}}
\put(1141,-3436){\makebox(0,0)[lb]{\smash{\SetFigFont{6}{7.2}{rm}\mbox{$#1$}}}}
\put(2341,-3436){\makebox(0,0)[lb]{\smash{\SetFigFont{6}{7.2}{rm}\mbox{$#3$}}}}
\put(3541,-3436){\makebox(0,0)[lb]{\smash{\SetFigFont{6}{7.2}{rm}\mbox{$#4$}}}}
\put(4741,-3436){\makebox(0,0)[lb]{\smash{\SetFigFont{6}{7.2}{rm}\mbox{$#5$}}}}
\put(5941,-3436){\makebox(0,0)[lb]{\smash{\SetFigFont{6}{7.2}{rm}\mbox{$#6$}}}}
\put(7141,-3436){\makebox(0,0)[lb]{\smash{\SetFigFont{6}{7.2}{rm}\mbox{$#7$}}}}
\put(8341,-3436){\makebox(0,0)[lb]{\smash{\SetFigFont{6}{7.2}{rm}\mbox{$#8$}}}}
\put(4024,-5263){\makebox(0,0)[lb]{\smash{\SetFigFont{6}{7.2}{rm}\mbox{$#2$}}}}
\end{picture}$}}
}
\begingroup\makeatletter\ifx\SetFigFont\undefined
\def\x#1#2#3#4#5#6#7\relax{\def\x{#1#2#3#4#5#6}}%
\expandafter\x\fmtname xxxxxx\relax \def\y{splain}%
\ifx\x\y   
\gdef\SetFigFont#1#2#3{%
  \ifnum #1<17\tiny\else \ifnum #1<20\small\else
  \ifnum #1<24\normalsize\else \ifnum #1<29\large\else
  \ifnum #1<34\Large\else \ifnum #1<41\LARGE\else
     \huge\fi\fi\fi\fi\fi\fi
  \csname #3\endcsname}%
\else
\gdef\SetFigFont#1#2#3{\begingroup
  \count@#1\relax \ifnum 25<\count@\count@25\fi
  \def\x{\endgroup\@setsize\SetFigFont{#2pt}}%
  \expandafter\x
    \csname \romannumeral\the\count@ pt\expandafter\endcsname
    \csname @\romannumeral\the\count@ pt\endcsname
  \csname #3\endcsname}%
\fi
\fi\endgroup

\newcommand{\Eviiit}[9]{\raisebox{2.4mm}{\parbox[c]{6cm}{$
\begin{picture}(8900,2200)(950,-5400)
\thinlines
\put(1201,-3961){\circle{336}}
\put(2401,-3961){\circle{336}}
\put(3601,-3961){\circle{336}}
\put(4801,-3961){\circle{336}}
\put(6001,-3961){\circle{336}}
\put(7201,-3961){\circle{336}}
\put(8401,-3961){\circle{336}}
\put(9601,-3961){\circle{336}}
\put(3601,-5161){\circle{336}}
\put(1369,-3961){\line( 1, 0){864}}
\put(2569,-3961){\line( 1, 0){864}}
\put(3769,-3961){\line( 1, 0){864}}
\put(4969,-3961){\line( 1, 0){864}}
\put(6169,-3961){\line( 1, 0){864}}
\put(7369,-3961){\line( 1, 0){864}}
\put(8569,-3961){\line( 1, 0){864}}
\put(3601,-4129){\line( 0, -1){864}}
\put(1141,-3436){\makebox(0,0)[lb]{\smash{\SetFigFont{6}{7.2}{rm}\mbox{$#1$}}}}
\put(2341,-3436){\makebox(0,0)[lb]{\smash{\SetFigFont{6}{7.2}{rm}\mbox{$#3$}}}}
\put(3541,-3436){\makebox(0,0)[lb]{\smash{\SetFigFont{6}{7.2}{rm}\mbox{$#4$}}}}
\put(4741,-3436){\makebox(0,0)[lb]{\smash{\SetFigFont{6}{7.2}{rm}\mbox{$#5$}}}}
\put(5941,-3436){\makebox(0,0)[lb]{\smash{\SetFigFont{6}{7.2}{rm}\mbox{$#6$}}}}
\put(7141,-3436){\makebox(0,0)[lb]{\smash{\SetFigFont{6}{7.2}{rm}\mbox{$#7$}}}}
\put(8341,-3436){\makebox(0,0)[lb]{\smash{\SetFigFont{6}{7.2}{rm}\mbox{$#8$}}}}
\put(9541,-3436){\makebox(0,0)[lb]{\smash{\SetFigFont{6}{7.2}{rm}\mbox{$#9$}}}}
\put(4024,-5263){\makebox(0,0)[lb]{\smash{\SetFigFont{6}{7.2}{rm}\mbox{$#2$}}}}
\end{picture}$}}
}

\begingroup\makeatletter\ifx\SetFigFont\undefined
\def\x#1#2#3#4#5#6#7\relax{\def\x{#1#2#3#4#5#6}}%
\expandafter\x\fmtname xxxxxx\relax \def\y{splain}%
\ifx\x\y   
\gdef\SetFigFont#1#2#3{%
  \ifnum #1<17\tiny\else \ifnum #1<20\small\else
  \ifnum #1<24\normalsize\else \ifnum #1<29\large\else
  \ifnum #1<34\Large\else \ifnum #1<41\LARGE\else
     \huge\fi\fi\fi\fi\fi\fi
  \csname #3\endcsname}%
\else
\gdef\SetFigFont#1#2#3{\begingroup
  \count@#1\relax \ifnum 25<\count@\count@25\fi
  \def\x{\endgroup\@setsize\SetFigFont{#2pt}}%
  \expandafter\x
    \csname \romannumeral\the\count@ pt\expandafter\endcsname
    \csname @\romannumeral\the\count@ pt\endcsname
  \csname #3\endcsname}%
\fi
\fi\endgroup

\newcommand{\Ant}[4]{
\begin{picture}(3952,1645)(1025,-3386)
\thinlines
\put(1201,-3961){\circle{336}}
\put(2401,-3961){\circle{336}}
\put(4771,-3961){\circle{336}}
\put(3436,-3961){\circle*{30}}
\put(3586,-3961){\circle*{30}}
\put(3736,-3961){\circle*{30}}
\put(2986,-2461){\circle{336}}
\put(1369,-3961){\line( 1, 0){870}}
\put(2569,-3961){\line( 1, 0){600}}
\put(4003,-3961){\line( 1, 0){600}}
\put(1330,-3831){\line( 6, 5){1510}}
\put(4641,-3831){\line( -6, 5){1510}}
\put(1141,-3436){\makebox(0,0)[lb]{\smash{\SetFigFont{6}{7.2}{rm}\mbox{$#1$}}}}
\put(2341,-3436){\makebox(0,0)[lb]{\smash{\SetFigFont{6}{7.2}{rm}\mbox{$#2$}}}}
\put(4741,-3436){\makebox(0,0)[lb]{\smash{\SetFigFont{6}{7.2}{rm}\mbox{$#3$}}}}
\put(2926,-1936){\makebox(0,0)[lb]{\smash{\SetFigFont{6}{7.2}{rm}\mbox{$#4$}}}}
\end{picture}
}
\begingroup\makeatletter\ifx\SetFigFont\undefined
\def\x#1#2#3#4#5#6#7\relax{\def\x{#1#2#3#4#5#6}}%
\expandafter\x\fmtname xxxxxx\relax \def\y{splain}%
\ifx\x\y   
\gdef\SetFigFont#1#2#3{%
  \ifnum #1<17\tiny\else \ifnum #1<20\small\else
  \ifnum #1<24\normalsize\else \ifnum #1<29\large\else
  \ifnum #1<34\Large\else \ifnum #1<41\LARGE\else
     \huge\fi\fi\fi\fi\fi\fi
  \csname #3\endcsname}%
\else
\gdef\SetFigFont#1#2#3{\begingroup
  \count@#1\relax \ifnum 25<\count@\count@25\fi
  \def\x{\endgroup\@setsize\SetFigFont{#2pt}}%
  \expandafter\x
    \csname \romannumeral\the\count@ pt\expandafter\endcsname
    \csname @\romannumeral\the\count@ pt\endcsname
  \csname #3\endcsname}%
\fi
\fi\endgroup

\newcommand{\Bnt}[5]{
\begin{picture}(5151,2200)(1123,-4136)
\thinlines
\put(1549,-3109){\circle{336}}
\put(1549,-4813){\circle{336}}
\put(2401,-3961){\circle{336}}
\put(4801,-3961){\circle{336}}
\put(6001,-3961){\circle{336}}
\put(3436,-3946){\circle*{30}}
\put(3586,-3946){\circle*{30}}
\put(3736,-3946){\circle*{30}}
\put(1667,-3227){\line( 1, -1){615}}
\put(1667,-4695){\line( 1, 1){615}}
\put(2569,-3961){\line( 1, 0){600}}
\put(4003,-3961){\line( 1, 0){600}}
\put(4969,-3911){\line( 1, 0){864}}
\put(4969,-4011){\line( 1, 0){864}}
\put(926,-3211){\makebox(0,0)[lb]{\smash{\SetFigFont{6}{7.2}{rm}\mbox{$#1$}}}}
\put(926,-4936){\makebox(0,0)[lb]{\smash{\SetFigFont{6}{7.2}{rm}\mbox{$#2$}}}}
\put(2341,-3436){\makebox(0,0)[lb]{\smash{\SetFigFont{6}{7.2}{rm}\mbox{$#3$}}}}
\put(4741,-3436){\makebox(0,0)[lb]{\smash{\SetFigFont{6}{7.2}{rm}\mbox{$#4$}}}}
\put(5941,-3436){\makebox(0,0)[lb]{\smash{\SetFigFont{6}{7.2}{rm}\mbox{$#5$}}}}
\end{picture}
}
\begingroup\makeatletter\ifx\SetFigFont\undefined
\def\x#1#2#3#4#5#6#7\relax{\def\x{#1#2#3#4#5#6}}%
\expandafter\x\fmtname xxxxxx\relax \def\y{splain}%
\ifx\x\y   
\gdef\SetFigFont#1#2#3{%
  \ifnum #1<17\tiny\else \ifnum #1<20\small\else
  \ifnum #1<24\normalsize\else \ifnum #1<29\large\else
  \ifnum #1<34\Large\else \ifnum #1<41\LARGE\else
     \huge\fi\fi\fi\fi\fi\fi
  \csname #3\endcsname}%
\else
\gdef\SetFigFont#1#2#3{\begingroup
  \count@#1\relax \ifnum 25<\count@\count@25\fi
  \def\x{\endgroup\@setsize\SetFigFont{#2pt}}%
  \expandafter\x
    \csname \romannumeral\the\count@ pt\expandafter\endcsname
    \csname @\romannumeral\the\count@ pt\endcsname
  \csname #3\endcsname}%
\fi
\fi\endgroup

\newcommand{\Bnta}[5]{
\begin{picture}(5151,2200)(1123,-4136)
\thinlines
\put(1549,-3109){\circle{336}}
\put(1549,-4813){\circle{336}}
\put(2401,-3961){\circle{336}}
\put(4801,-3961){\circle{336}}
\put(6001,-3961){\circle{336}}
\put(3436,-3946){\circle*{30}}
\put(3586,-3946){\circle*{30}}
\put(3736,-3946){\circle*{30}}
\put(1667,-3227){\line( 1, -1){615}}
\put(1667,-4695){\line( 1, 1){615}}
\put(2569,-3961){\line( 1, 0){600}}
\put(4003,-3961){\line( 1, 0){600}}
\put(4969,-3911){\line( 1, 0){764}}
\put(4969,-4011){\line( 1, 0){764}}
\put(5400,-4131){$>$}
\put(926,-3211){\makebox(0,0)[lb]{\smash{\SetFigFont{6}{7.2}{rm}\mbox{$#1$}}}}
\put(926,-4936){\makebox(0,0)[lb]{\smash{\SetFigFont{6}{7.2}{rm}\mbox{$#2$}}}}
\put(2341,-3436){\makebox(0,0)[lb]{\smash{\SetFigFont{6}{7.2}{rm}\mbox{$#3$}}}}
\put(4741,-3436){\makebox(0,0)[lb]{\smash{\SetFigFont{6}{7.2}{rm}\mbox{$#4$}}}}
\put(5941,-3436){\makebox(0,0)[lb]{\smash{\SetFigFont{6}{7.2}{rm}\mbox{$#5$}}}}
\end{picture}
}
\begingroup\makeatletter\ifx\SetFigFont\undefined
\def\x#1#2#3#4#5#6#7\relax{\def\x{#1#2#3#4#5#6}}%
\expandafter\x\fmtname xxxxxx\relax \def\y{splain}%
\ifx\x\y   
\gdef\SetFigFont#1#2#3{%
  \ifnum #1<17\tiny\else \ifnum #1<20\small\else
  \ifnum #1<24\normalsize\else \ifnum #1<29\large\else
  \ifnum #1<34\Large\else \ifnum #1<41\LARGE\else
     \huge\fi\fi\fi\fi\fi\fi
  \csname #3\endcsname}%
\else
\gdef\SetFigFont#1#2#3{\begingroup
  \count@#1\relax \ifnum 25<\count@\count@25\fi
  \def\x{\endgroup\@setsize\SetFigFont{#2pt}}%
  \expandafter\x
    \csname \romannumeral\the\count@ pt\expandafter\endcsname
    \csname @\romannumeral\the\count@ pt\endcsname
  \csname #3\endcsname}%
\fi
\fi\endgroup

\newcommand{\Bntaii}[5]{
\begin{picture}(5151,2200)(1123,-4136)
\thinlines
\put(1549,-3109){\circle{336}}
\put(1549,-4813){\circle{336}}
\put(2401,-3961){\circle{336}}
\put(4801,-3961){\circle{336}}
\put(6001,-3961){\circle{336}}
\put(3436,-3946){\circle*{30}}
\put(3586,-3946){\circle*{30}}
\put(3736,-3946){\circle*{30}}
\put(1667,-3227){\line( 1, -1){615}}
\put(1667,-4695){\line( 1, 1){615}}
\put(2569,-3961){\line( 1, 0){600}}
\put(4003,-3961){\line( 1, 0){600}}
\put(5069,-3911){\line( 1, 0){764}}
\put(5069,-4011){\line( 1, 0){764}}
\put(4896,-4131){$<$}
\put(926,-3211){\makebox(0,0)[lb]{\smash{\SetFigFont{6}{7.2}{rm}\mbox{$#1$}}}}
\put(926,-4936){\makebox(0,0)[lb]{\smash{\SetFigFont{6}{7.2}{rm}\mbox{$#2$}}}}
\put(2341,-3436){\makebox(0,0)[lb]{\smash{\SetFigFont{6}{7.2}{rm}\mbox{$#3$}}}}
\put(4741,-3436){\makebox(0,0)[lb]{\smash{\SetFigFont{6}{7.2}{rm}\mbox{$#4$}}}}
\put(5941,-3436){\makebox(0,0)[lb]{\smash{\SetFigFont{6}{7.2}{rm}\mbox{$#5$}}}}
\end{picture}
}
\begingroup\makeatletter\ifx\SetFigFont\undefined
\def\x#1#2#3#4#5#6#7\relax{\def\x{#1#2#3#4#5#6}}%
\expandafter\x\fmtname xxxxxx\relax \def\y{splain}%
\ifx\x\y   
\gdef\SetFigFont#1#2#3{%
  \ifnum #1<17\tiny\else \ifnum #1<20\small\else
  \ifnum #1<24\normalsize\else \ifnum #1<29\large\else
  \ifnum #1<34\Large\else \ifnum #1<41\LARGE\else
     \huge\fi\fi\fi\fi\fi\fi
  \csname #3\endcsname}%
\else
\gdef\SetFigFont#1#2#3{\begingroup
  \count@#1\relax \ifnum 25<\count@\count@25\fi
  \def\x{\endgroup\@setsize\SetFigFont{#2pt}}%
  \expandafter\x
    \csname \romannumeral\the\count@ pt\expandafter\endcsname
    \csname @\romannumeral\the\count@ pt\endcsname
  \csname #3\endcsname}%
\fi
\fi\endgroup

\newcommand{\Bntar}[5]{
\begin{picture}(5151,2200)(1123,-4136)
\thinlines
\put(1549,-3109){\circle{336}}
\put(1549,-4813){\circle{336}}
\put(2401,-3961){\circle{336}}
\put(4801,-3961){\circle{336}}
\put(6001,-3961){\circle{236}}
\put(6001,-3961){\circle{436}}
\put(3436,-3946){\circle*{30}}
\put(3586,-3946){\circle*{30}}
\put(3736,-3946){\circle*{30}}
\put(1667,-3227){\line( 1, -1){615}}
\put(1667,-4695){\line( 1, 1){615}}
\put(2569,-3961){\line( 1, 0){600}}
\put(4003,-3961){\line( 1, 0){600}}
\put(4969,-3911){\line( 1, 0){714}}
\put(4969,-4011){\line( 1, 0){714}}
\put(5350,-4131){$>$}
\put(926,-3211){\makebox(0,0)[lb]{\smash{\SetFigFont{6}{7.2}{rm}\mbox{$#1$}}}}
\put(926,-4936){\makebox(0,0)[lb]{\smash{\SetFigFont{6}{7.2}{rm}\mbox{$#2$}}}}
\put(2341,-3436){\makebox(0,0)[lb]{\smash{\SetFigFont{6}{7.2}{rm}\mbox{$#3$}}}}
\put(4741,-3436){\makebox(0,0)[lb]{\smash{\SetFigFont{6}{7.2}{rm}\mbox{$#4$}}}}
\put(5941,-3436){\makebox(0,0)[lb]{\smash{\SetFigFont{6}{7.2}{rm}\mbox{$#5$}}}}
\end{picture}
}

\begingroup\makeatletter\ifx\SetFigFont\undefined
\def\x#1#2#3#4#5#6#7\relax{\def\x{#1#2#3#4#5#6}}%
\expandafter\x\fmtname xxxxxx\relax \def\y{splain}%
\ifx\x\y   
\gdef\SetFigFont#1#2#3{%
  \ifnum #1<17\tiny\else \ifnum #1<20\small\else
  \ifnum #1<24\normalsize\else \ifnum #1<29\large\else
  \ifnum #1<34\Large\else \ifnum #1<41\LARGE\else
     \huge\fi\fi\fi\fi\fi\fi
  \csname #3\endcsname}%
\else
\gdef\SetFigFont#1#2#3{\begingroup
  \count@#1\relax \ifnum 25<\count@\count@25\fi
  \def\x{\endgroup\@setsize\SetFigFont{#2pt}}%
  \expandafter\x
    \csname \romannumeral\the\count@ pt\expandafter\endcsname
    \csname @\romannumeral\the\count@ pt\endcsname
  \csname #3\endcsname}%
\fi
\fi\endgroup

\newcommand{\Cnt}[6]{
\begin{picture}(7552,895)(1025,-4136)
\thinlines
\put(1201,-3961){\circle{336}}
\put(2401,-3961){\circle{336}}
\put(3601,-3961){\circle{336}}
\put(6001,-3961){\circle{336}}
\put(7201,-3961){\circle{336}}
\put(8401,-3961){\circle{336}}
\put(1369,-3911){\line( 1, 0){864}}
\put(1369,-4011){\line( 1, 0){864}}
\put(2569,-3961){\line( 1, 0){864}}
\put(3769,-3961){\line( 1, 0){600}}
\put(5203,-3961){\line( 1, 0){600}}
\put(6169,-3961){\line( 1, 0){864}}
\put(7369,-3911){\line( 1, 0){864}}
\put(7369,-4011){\line( 1, 0){864}}

\put(4636,-3946){\circle*{30}}
\put(4786,-3946){\circle*{30}}
\put(4936,-3946){\circle*{30}}

\put(1141,-3436){\makebox(0,0)[lb]{\smash{\SetFigFont{6}{7.2}{rm}\mbox{$#1$}}}}
\put(2341,-3436){\makebox(0,0)[lb]{\smash{\SetFigFont{6}{7.2}{rm}\mbox{$#2$}}}}
\put(3541,-3436){\makebox(0,0)[lb]{\smash{\SetFigFont{6}{7.2}{rm}\mbox{$#3$}}}}
\put(5941,-3436){\makebox(0,0)[lb]{\smash{\SetFigFont{6}{7.2}{rm}\mbox{$#4$}}}}
\put(7141,-3436){\makebox(0,0)[lb]{\smash{\SetFigFont{6}{7.2}{rm}\mbox{$#5$}}}}
\put(8341,-3436){\makebox(0,0)[lb]{\smash{\SetFigFont{6}{7.2}{rm}\mbox{$#6$}}}}

\end{picture}
}
\begingroup\makeatletter\ifx\SetFigFont\undefined
\def\x#1#2#3#4#5#6#7\relax{\def\x{#1#2#3#4#5#6}}%
\expandafter\x\fmtname xxxxxx\relax \def\y{splain}%
\ifx\x\y   
\gdef\SetFigFont#1#2#3{%
  \ifnum #1<17\tiny\else \ifnum #1<20\small\else
  \ifnum #1<24\normalsize\else \ifnum #1<29\large\else
  \ifnum #1<34\Large\else \ifnum #1<41\LARGE\else
     \huge\fi\fi\fi\fi\fi\fi
  \csname #3\endcsname}%
\else
\gdef\SetFigFont#1#2#3{\begingroup
  \count@#1\relax \ifnum 25<\count@\count@25\fi
  \def\x{\endgroup\@setsize\SetFigFont{#2pt}}%
  \expandafter\x
    \csname \romannumeral\the\count@ pt\expandafter\endcsname
    \csname @\romannumeral\the\count@ pt\endcsname
  \csname #3\endcsname}%
\fi
\fi\endgroup

\newcommand{\Cnta}[6]{
\begin{picture}(7552,895)(1025,-4136)
\thinlines
\put(1201,-3961){\circle{336}}
\put(2401,-3961){\circle{336}}
\put(3601,-3961){\circle{336}}
\put(6001,-3961){\circle{336}}
\put(7201,-3961){\circle{336}}
\put(8401,-3961){\circle{336}}
\put(1369,-3911){\line( 1, 0){764}}
\put(1369,-4011){\line( 1, 0){764}}
\put(1800,-4131){$>$}
\put(2569,-3961){\line( 1, 0){864}}
\put(3769,-3961){\line( 1, 0){600}}
\put(5203,-3961){\line( 1, 0){600}}
\put(6169,-3961){\line( 1, 0){864}}
\put(7369,-3911){\line( 1, 0){764}}
\put(7369,-4011){\line( 1, 0){764}}
\put(7800,-4131){$>$}

\put(4636,-3946){\circle*{30}}
\put(4786,-3946){\circle*{30}}
\put(4936,-3946){\circle*{30}}

\put(1141,-3436){\makebox(0,0)[lb]{\smash{\SetFigFont{6}{7.2}{rm}\mbox{$#1$}}}}
\put(2341,-3436){\makebox(0,0)[lb]{\smash{\SetFigFont{6}{7.2}{rm}\mbox{$#2$}}}}
\put(3541,-3436){\makebox(0,0)[lb]{\smash{\SetFigFont{6}{7.2}{rm}\mbox{$#3$}}}}
\put(5941,-3436){\makebox(0,0)[lb]{\smash{\SetFigFont{6}{7.2}{rm}\mbox{$#4$}}}}
\put(7141,-3436){\makebox(0,0)[lb]{\smash{\SetFigFont{6}{7.2}{rm}\mbox{$#5$}}}}
\put(8341,-3436){\makebox(0,0)[lb]{\smash{\SetFigFont{6}{7.2}{rm}\mbox{$#6$}}}}

\end{picture}
}
\begingroup\makeatletter\ifx\SetFigFont\undefined
\def\x#1#2#3#4#5#6#7\relax{\def\x{#1#2#3#4#5#6}}%
\expandafter\x\fmtname xxxxxx\relax \def\y{splain}%
\ifx\x\y   
\gdef\SetFigFont#1#2#3{%
  \ifnum #1<17\tiny\else \ifnum #1<20\small\else
  \ifnum #1<24\normalsize\else \ifnum #1<29\large\else
  \ifnum #1<34\Large\else \ifnum #1<41\LARGE\else
     \huge\fi\fi\fi\fi\fi\fi
  \csname #3\endcsname}%
\else
\gdef\SetFigFont#1#2#3{\begingroup
  \count@#1\relax \ifnum 25<\count@\count@25\fi
  \def\x{\endgroup\@setsize\SetFigFont{#2pt}}%
  \expandafter\x
    \csname \romannumeral\the\count@ pt\expandafter\endcsname
    \csname @\romannumeral\the\count@ pt\endcsname
  \csname #3\endcsname}%
\fi
\fi\endgroup

\newcommand{\Cntaii}[6]{
\begin{picture}(7552,895)(1025,-4136)
\thinlines
\put(1201,-3961){\circle{336}}
\put(2401,-3961){\circle{336}}
\put(3601,-3961){\circle{336}}
\put(6001,-3961){\circle{336}}
\put(7201,-3961){\circle{336}}
\put(8401,-3961){\circle{336}}
\put(1469,-3911){\line( 1, 0){764}}
\put(1469,-4011){\line( 1, 0){764}}
\put(1296,-4131){$<$}
\put(2569,-3961){\line( 1, 0){864}}
\put(3769,-3961){\line( 1, 0){600}}
\put(5203,-3961){\line( 1, 0){600}}
\put(6169,-3961){\line( 1, 0){864}}
\put(7369,-3911){\line( 1, 0){764}}
\put(7369,-4011){\line( 1, 0){764}}
\put(7800,-4131){$>$}

\put(4636,-3946){\circle*{30}}
\put(4786,-3946){\circle*{30}}
\put(4936,-3946){\circle*{30}}

\put(1141,-3436){\makebox(0,0)[lb]{\smash{\SetFigFont{6}{7.2}{rm}\mbox{$#1$}}}}
\put(2341,-3436){\makebox(0,0)[lb]{\smash{\SetFigFont{6}{7.2}{rm}\mbox{$#2$}}}}
\put(3541,-3436){\makebox(0,0)[lb]{\smash{\SetFigFont{6}{7.2}{rm}\mbox{$#3$}}}}
\put(5941,-3436){\makebox(0,0)[lb]{\smash{\SetFigFont{6}{7.2}{rm}\mbox{$#4$}}}}
\put(7141,-3436){\makebox(0,0)[lb]{\smash{\SetFigFont{6}{7.2}{rm}\mbox{$#5$}}}}
\put(8341,-3436){\makebox(0,0)[lb]{\smash{\SetFigFont{6}{7.2}{rm}\mbox{$#6$}}}}

\end{picture}
}
\begingroup\makeatletter\ifx\SetFigFont\undefined
\def\x#1#2#3#4#5#6#7\relax{\def\x{#1#2#3#4#5#6}}%
\expandafter\x\fmtname xxxxxx\relax \def\y{splain}%
\ifx\x\y   
\gdef\SetFigFont#1#2#3{%
  \ifnum #1<17\tiny\else \ifnum #1<20\small\else
  \ifnum #1<24\normalsize\else \ifnum #1<29\large\else
  \ifnum #1<34\Large\else \ifnum #1<41\LARGE\else
     \huge\fi\fi\fi\fi\fi\fi
  \csname #3\endcsname}%
\else
\gdef\SetFigFont#1#2#3{\begingroup
  \count@#1\relax \ifnum 25<\count@\count@25\fi
  \def\x{\endgroup\@setsize\SetFigFont{#2pt}}%
  \expandafter\x
    \csname \romannumeral\the\count@ pt\expandafter\endcsname
    \csname @\romannumeral\the\count@ pt\endcsname
  \csname #3\endcsname}%
\fi
\fi\endgroup

\newcommand{\Cntaiii}[6]{
\begin{picture}(7552,895)(1025,-4136)
\thinlines
\put(1201,-3961){\circle{336}}
\put(2401,-3961){\circle{336}}
\put(3601,-3961){\circle{336}}
\put(6001,-3961){\circle{336}}
\put(7201,-3961){\circle{336}}
\put(8401,-3961){\circle{336}}
\put(1369,-3911){\line( 1, 0){764}}
\put(1369,-4011){\line( 1, 0){764}}
\put(1800,-4131){$>$}
\put(2569,-3961){\line( 1, 0){864}}
\put(3769,-3961){\line( 1, 0){600}}
\put(5203,-3961){\line( 1, 0){600}}
\put(6169,-3961){\line( 1, 0){864}}
\put(7469,-3911){\line( 1, 0){764}}
\put(7469,-4011){\line( 1, 0){764}}
\put(7296,-4131){$<$}

\put(4636,-3946){\circle*{30}}
\put(4786,-3946){\circle*{30}}
\put(4936,-3946){\circle*{30}}

\put(1141,-3436){\makebox(0,0)[lb]{\smash{\SetFigFont{6}{7.2}{rm}\mbox{$#1$}}}}
\put(2341,-3436){\makebox(0,0)[lb]{\smash{\SetFigFont{6}{7.2}{rm}\mbox{$#2$}}}}
\put(3541,-3436){\makebox(0,0)[lb]{\smash{\SetFigFont{6}{7.2}{rm}\mbox{$#3$}}}}
\put(5941,-3436){\makebox(0,0)[lb]{\smash{\SetFigFont{6}{7.2}{rm}\mbox{$#4$}}}}
\put(7141,-3436){\makebox(0,0)[lb]{\smash{\SetFigFont{6}{7.2}{rm}\mbox{$#5$}}}}
\put(8341,-3436){\makebox(0,0)[lb]{\smash{\SetFigFont{6}{7.2}{rm}\mbox{$#6$}}}}

\end{picture}
}
\begingroup\makeatletter\ifx\SetFigFont\undefined
\def\x#1#2#3#4#5#6#7\relax{\def\x{#1#2#3#4#5#6}}%
\expandafter\x\fmtname xxxxxx\relax \def\y{splain}%
\ifx\x\y   
\gdef\SetFigFont#1#2#3{%
  \ifnum #1<17\tiny\else \ifnum #1<20\small\else
  \ifnum #1<24\normalsize\else \ifnum #1<29\large\else
  \ifnum #1<34\Large\else \ifnum #1<41\LARGE\else
     \huge\fi\fi\fi\fi\fi\fi
  \csname #3\endcsname}%
\else
\gdef\SetFigFont#1#2#3{\begingroup
  \count@#1\relax \ifnum 25<\count@\count@25\fi
  \def\x{\endgroup\@setsize\SetFigFont{#2pt}}%
  \expandafter\x
    \csname \romannumeral\the\count@ pt\expandafter\endcsname
    \csname @\romannumeral\the\count@ pt\endcsname
  \csname #3\endcsname}%
\fi
\fi\endgroup

\newcommand{\Cntar}[6]{
\begin{picture}(7552,895)(1025,-4136)
\thinlines
\put(1201,-3961){\circle{236}}
\put(1201,-3961){\circle{436}}
\put(2401,-3961){\circle{336}}
\put(3601,-3961){\circle{336}}
\put(6001,-3961){\circle{336}}
\put(7201,-3961){\circle{336}}
\put(8401,-3961){\circle{336}}
\put(1519,-3911){\line( 1, 0){714}}
\put(1519,-4011){\line( 1, 0){714}}
\put(1346,-4131){$<$}
\put(2569,-3961){\line( 1, 0){864}}
\put(3769,-3961){\line( 1, 0){600}}
\put(5203,-3961){\line( 1, 0){600}}
\put(6169,-3961){\line( 1, 0){864}}
\put(7369,-3911){\line( 1, 0){764}}
\put(7369,-4011){\line( 1, 0){764}}
\put(7800,-4131){$>$}

\put(4636,-3946){\circle*{30}}
\put(4786,-3946){\circle*{30}}
\put(4936,-3946){\circle*{30}}

\put(1141,-3436){\makebox(0,0)[lb]{\smash{\SetFigFont{6}{7.2}{rm}\mbox{$#1$}}}}
\put(2341,-3436){\makebox(0,0)[lb]{\smash{\SetFigFont{6}{7.2}{rm}\mbox{$#2$}}}}
\put(3541,-3436){\makebox(0,0)[lb]{\smash{\SetFigFont{6}{7.2}{rm}\mbox{$#3$}}}}
\put(5941,-3436){\makebox(0,0)[lb]{\smash{\SetFigFont{6}{7.2}{rm}\mbox{$#4$}}}}
\put(7141,-3436){\makebox(0,0)[lb]{\smash{\SetFigFont{6}{7.2}{rm}\mbox{$#5$}}}}
\put(8341,-3436){\makebox(0,0)[lb]{\smash{\SetFigFont{6}{7.2}{rm}\mbox{$#6$}}}}

\end{picture}
}

\begingroup\makeatletter\ifx\SetFigFont\undefined
\def\x#1#2#3#4#5#6#7\relax{\def\x{#1#2#3#4#5#6}}%
\expandafter\x\fmtname xxxxxx\relax \def\y{splain}%
\ifx\x\y   
\gdef\SetFigFont#1#2#3{%
  \ifnum #1<17\tiny\else \ifnum #1<20\small\else
  \ifnum #1<24\normalsize\else \ifnum #1<29\large\else
  \ifnum #1<34\Large\else \ifnum #1<41\LARGE\else
     \huge\fi\fi\fi\fi\fi\fi
  \csname #3\endcsname}%
\else
\gdef\SetFigFont#1#2#3{\begingroup
  \count@#1\relax \ifnum 25<\count@\count@25\fi
  \def\x{\endgroup\@setsize\SetFigFont{#2pt}}%
  \expandafter\x
    \csname \romannumeral\the\count@ pt\expandafter\endcsname
    \csname @\romannumeral\the\count@ pt\endcsname
  \csname #3\endcsname}%
\fi
\fi\endgroup

\newcommand{\Cntarii}[6]{
\begin{picture}(7552,895)(1025,-4136)
\thinlines
\put(1201,-3961){\circle{236}}
\put(1201,-3961){\circle{436}}
\put(2401,-3961){\circle{336}}
\put(3601,-3961){\circle{336}}
\put(6001,-3961){\circle{336}}
\put(7201,-3961){\circle{336}}
\put(8401,-3961){\circle{336}}
\put(1519,-3911){\line( 1, 0){714}}
\put(1519,-4011){\line( 1, 0){714}}
\put(1346,-4131){$<$}
\put(2569,-3961){\line( 1, 0){864}}
\put(3769,-3961){\line( 1, 0){600}}
\put(5203,-3961){\line( 1, 0){600}}
\put(6169,-3961){\line( 1, 0){864}}
\put(7469,-3911){\line( 1, 0){764}}
\put(7469,-4011){\line( 1, 0){764}}
\put(7296,-4131){$<$}

\put(4636,-3946){\circle*{30}}
\put(4786,-3946){\circle*{30}}
\put(4936,-3946){\circle*{30}}

\put(1141,-3436){\makebox(0,0)[lb]{\smash{\SetFigFont{6}{7.2}{rm}\mbox{$#1$}}}}
\put(2341,-3436){\makebox(0,0)[lb]{\smash{\SetFigFont{6}{7.2}{rm}\mbox{$#2$}}}}
\put(3541,-3436){\makebox(0,0)[lb]{\smash{\SetFigFont{6}{7.2}{rm}\mbox{$#3$}}}}
\put(5941,-3436){\makebox(0,0)[lb]{\smash{\SetFigFont{6}{7.2}{rm}\mbox{$#4$}}}}
\put(7141,-3436){\makebox(0,0)[lb]{\smash{\SetFigFont{6}{7.2}{rm}\mbox{$#5$}}}}
\put(8341,-3436){\makebox(0,0)[lb]{\smash{\SetFigFont{6}{7.2}{rm}\mbox{$#6$}}}}

\end{picture}
}

\begingroup\makeatletter\ifx\SetFigFont\undefined
\def\x#1#2#3#4#5#6#7\relax{\def\x{#1#2#3#4#5#6}}%
\expandafter\x\fmtname xxxxxx\relax \def\y{splain}%
\ifx\x\y   
\gdef\SetFigFont#1#2#3{%
  \ifnum #1<17\tiny\else \ifnum #1<20\small\else
  \ifnum #1<24\normalsize\else \ifnum #1<29\large\else
  \ifnum #1<34\Large\else \ifnum #1<41\LARGE\else
     \huge\fi\fi\fi\fi\fi\fi
  \csname #3\endcsname}%
\else
\gdef\SetFigFont#1#2#3{\begingroup
  \count@#1\relax \ifnum 25<\count@\count@25\fi
  \def\x{\endgroup\@setsize\SetFigFont{#2pt}}%
  \expandafter\x
    \csname \romannumeral\the\count@ pt\expandafter\endcsname
    \csname @\romannumeral\the\count@ pt\endcsname
  \csname #3\endcsname}%
\fi
\fi\endgroup

\newcommand{\Cntariii}[6]{
\begin{picture}(7552,895)(1025,-4136)
\thinlines
\put(1201,-3961){\circle{236}}
\put(1201,-3961){\circle{436}}
\put(2401,-3961){\circle{336}}
\put(3601,-3961){\circle{336}}
\put(6001,-3961){\circle{336}}
\put(7201,-3961){\circle{336}}
\put(8401,-3961){\circle{236}}
\put(8401,-3961){\circle{436}}
\put(1519,-3911){\line( 1, 0){714}}
\put(1519,-4011){\line( 1, 0){714}}
\put(1346,-4131){$<$}
\put(2569,-3961){\line( 1, 0){864}}
\put(3769,-3961){\line( 1, 0){600}}
\put(5203,-3961){\line( 1, 0){600}}
\put(6169,-3961){\line( 1, 0){864}}
\put(7369,-3911){\line( 1, 0){714}}
\put(7369,-4011){\line( 1, 0){714}}
\put(7750,-4131){$>$}

\put(4636,-3946){\circle*{30}}
\put(4786,-3946){\circle*{30}}
\put(4936,-3946){\circle*{30}}

\put(1141,-3436){\makebox(0,0)[lb]{\smash{\SetFigFont{6}{7.2}{rm}\mbox{$#1$}}}}
\put(2341,-3436){\makebox(0,0)[lb]{\smash{\SetFigFont{6}{7.2}{rm}\mbox{$#2$}}}}
\put(3541,-3436){\makebox(0,0)[lb]{\smash{\SetFigFont{6}{7.2}{rm}\mbox{$#3$}}}}
\put(5941,-3436){\makebox(0,0)[lb]{\smash{\SetFigFont{6}{7.2}{rm}\mbox{$#4$}}}}
\put(7141,-3436){\makebox(0,0)[lb]{\smash{\SetFigFont{6}{7.2}{rm}\mbox{$#5$}}}}
\put(8341,-3436){\makebox(0,0)[lb]{\smash{\SetFigFont{6}{7.2}{rm}\mbox{$#6$}}}}

\end{picture}
}

\begingroup\makeatletter\ifx\SetFigFont\undefined
\def\x#1#2#3#4#5#6#7\relax{\def\x{#1#2#3#4#5#6}}%
\expandafter\x\fmtname xxxxxx\relax \def\y{splain}%
\ifx\x\y   
\gdef\SetFigFont#1#2#3{%
  \ifnum #1<17\tiny\else \ifnum #1<20\small\else
  \ifnum #1<24\normalsize\else \ifnum #1<29\large\else
  \ifnum #1<34\Large\else \ifnum #1<41\LARGE\else
     \huge\fi\fi\fi\fi\fi\fi
  \csname #3\endcsname}%
\else
\gdef\SetFigFont#1#2#3{\begingroup
  \count@#1\relax \ifnum 25<\count@\count@25\fi
  \def\x{\endgroup\@setsize\SetFigFont{#2pt}}%
  \expandafter\x
    \csname \romannumeral\the\count@ pt\expandafter\endcsname
    \csname @\romannumeral\the\count@ pt\endcsname
  \csname #3\endcsname}%
\fi
\fi\endgroup

\newcommand{\Dnt}[6]{\parbox[c]{3.4cm}{$
\begin{picture}(4803,2200)(4525,-5050)
\thinlines
\put(5149,-3109){\circle{336}}
\put(5149,-4813){\circle{336}}
\put(6001,-3961){\circle{336}}
\put(8401,-3961){\circle{336}}
\put(9253,-3109){\circle{336}}
\put(9253,-4813){\circle{336}}
\put(7036,-3946){\circle*{30}}
\put(7186,-3946){\circle*{30}}
\put(7336,-3946){\circle*{30}}
\put(5267,-3227){\line( 1, -1){615}}
\put(5267,-4695){\line( 1, 1){615}}
\put(6169,-3961){\line( 1, 0){600}}
\put(7603,-3961){\line( 1, 0){600}}
\put(8519,-3843){\line( 1, 1){615}}
\put(8519,-4079){\line( 1, -1){615}}
\put(4476,-3211){\makebox(0,0)[lb]{\smash{\SetFigFont{6}{7.2}{rm}\mbox{$#1$}}}}
\put(4476,-4936){\makebox(0,0)[lb]{\smash{\SetFigFont{6}{7.2}{rm}\mbox{$#2$}}}}
\put(5941,-3436){\makebox(0,0)[lb]{\smash{\SetFigFont{6}{7.2}{rm}\mbox{$#3$}}}}
\put(8341,-3436){\makebox(0,0)[lb]{\smash{\SetFigFont{6}{7.2}{rm}\mbox{$#4$}}}}
\put(9676,-3211){\makebox(0,0)[lb]{\smash{\SetFigFont{6}{7.2}{rm}\mbox{$#5$}}}}
\put(9676,-4936){\makebox(0,0)[lb]{\smash{\SetFigFont{6}{7.2}{rm}\mbox{$#6$}}}}
\end{picture}$}
}

\newcommand{\Lg}{\mbox{$\mathfrak g$}}

\newcommand{\Lk}{\mbox{$\mathfrak k$}}
\newcommand{\Lp}{\mbox{$\mathfrak p$}}
\newcommand{\La}{\mbox{$\mathfrak a$}}
\newcommand{\Lm}{\mbox{$\mathfrak m$}}

\newcommand{\D}{\mbox{$\mathcal D$}}

\newcommand{\Pf}{{\em Proof}. }
\newcommand{\EPf}{\hfill$\square$}
\newcommand{\Z}{\mbox{$\mathbb Z$}}
\newcommand{\R}{\mbox{$\mathbb R$}}

\newcommand{\inn}[2]{\mbox{$\mathcal{h} #1,#2 \mathcal{i}$}}

\newcommand{\ontop}[2]{\genfrac{}{}{0pt}{}{#1}{#2}}

\newtheorem{thm}[equation]{Theorem}

\newtheorem{cor}[equation]{Corollary}
\newtheorem{prop}[equation]{Proposition}
\newtheorem{lem}[equation]{Lemma}

\theoremstyle{remark}
\newtheorem{rem}[equation]{Remark}

\theoremstyle{remark}
\newtheorem{defn}[equation]{Definition}

\begin{document}
 
\title{Homogeneous structures and rigidity of\\ isoparametric submanifolds
in Hilbert space} 
\date{\today}

\author{Claudio Gorodski}
\address{Instituto de Matem\'atica e Estat\'\i stica, Universidade de 
S\~ao Paulo, Brazil}
\email{gorodski@ime.usp.br}

\author{Ernst Heintze}
\address{Institut f\"ur Mathematik, Universit\"at Augsburg, Germany}  
\email{ernst.heintze@math.uni-augsburg.de}

\begin{abstract}
We study isoparametric 
submanifolds of rank at least two 
in a separable Hilbert space,
which are known to be homogeneous 
by the main result in~\cite{HL},
and associate
to such a submanifold $M$ and a point $x$ in $M$ 
a canonical homogeneous structure $\Gamma_x$ (a certain bilinear map
defined on a subspace of $T_xM\times T_xM$). 
We prove that $\Gamma_x$ together with the second fundamental form $\alpha_x$ encodes all the information about~$M$, 
and deduce from this the rigidity result that $M$ 
is completely determined by $\alpha_x$ and 
$(\nabla\alpha)_x$, thereby making such submanifolds accessible to 
classification.  
As an essential step, we show that the one-parameter groups of isometries
constructed in~\cite{HL} to prove their homogeneity
induce smooth 
and hence everywhere defined Killing fields, implying the continuity of 
$\Gamma$ (this result
also seems to close a gap in~\cite{Ch}). Here an important tool is 
the introduction of affine root systems of isoparametric submanifolds.
\end{abstract}

\thanks{The first author was partially supported by the
CNPq grant 302472/2009-6 and the FAPESP project 2007/03192-7.} 

\dedicatory{Dedicated to Richard Palais on the occasion
     of his $80$th birthday}
\maketitle


\section{Introduction}

In ~\cite{HPTT} Richard Palais and his co-authors discussed among other things 
relations between isoparametric submanifolds in Hilbert space and affine 
Kac-Moody algebras. The present paper may be seen as a continuation of that 
line of research. It contributes to the conjecture that all isoparametric 
submanifolds of rank at least two in an infinite dimensional Hilbert space 
arise as principal orbits of isotropy representations of symmetric spaces 
of affine Kac-Moody type (which are obtained from involutions of the second 
kind of affine Kac-Moody groups).

We begin with a simple example in finite dimensions 
that motivates our main construction in infinite dimensions.
Let $M=G/K$ be a homogeneous space embedded in a Euclidean space $V$ 
as an orbit of a compact connected group $G$ of isometries of $V$, where
$K$ is the isotropy subgroup at $x\in M$. Then the Lie algebra of $G$ 
admits a reductive decomposition as $\Lg=\Lk+\Lm$, where $\Lk$ is the 
Lie algebra of $K$. Each element of $\Lm$ is a Killing field
on $V$ and this defines an 
isomorphism $\Lm \rightarrow T_xM$ (by evaluating the Killing field at $x$) 
whose inverse
we denote by $X\mapsto\check X$.
Now the interesting point is that the bilinear mapping
$\bar\Gamma:=\bar\Gamma_x:T_xM\times T_xV\to T_xV$ 
defined by
$\bar\Gamma_XY := \frac{d}{dt}|_{t=0}
(\exp t\check X)_*(Y)$
determines $M$ completely (as well as the reductive complement).
In fact, the Killing fields $\check X \in \Lm$ are completely determined
by their value and derivative at $x$, 
which are $X:=\check X(x)$ and $\bar\Gamma_X$.  
Thus $\bar\Gamma$ determines $\Lm$ and, since the subgroup corresponding to the Lie algebra generated by $\Lm$ acts transitively on the submanifold, 
$M$ is determined as well.
  
In the special case in which 
$M$ has flat normal bundle and $G$ induces parallel translation
along it (the only case in which we will be interested),
we have $\bar\Gamma_X\xi=-A_\xi X$ and thus by skew-symmetry of $\bar\Gamma _X$, $\bar\Gamma_XY=\Gamma_XY+\alpha_x(X,Y)$ 
for all $X$, $Y\in T_xM$ and $\xi\in\nu_xM$,
where $\alpha_x$ is the second fundamental form of $M$ at $x$,
$A_\xi$ is the shape operator in the direction of $\xi$, and 
\[ \Gamma=\Gamma_x:T_xM\times T_xM\to T_xM \]
is the tangential component of $\bar\Gamma$, i.e.
\[  \Gamma_XY := \left(\frac{d}{dt}\Big|_{t=0}
(\exp t\check X)_*(Y)\right)^\top. \]
Thus the pair $(\alpha_x,\Gamma_x)$ contains the same information as  
$\bar\Gamma_x$ and  
we call $\Gamma$ a \emph{homogeneous 
structure} for $M$. 

The main goal of this paper is to carry over these ideas to
isoparametric submanifolds of rank at least $2$ in Hilbert space 
with the ultimate goal of proving the 
above mentioned conjecture in a forthcoming paper
by restricting the possibilities for
$\alpha$ and $\Gamma$ so much that only the examples
coming from affine Kac-Moody algebras remain. 
The restrictions on $\alpha$ are essentially already known 
(as $\alpha$ is determined by the affine Weyl group and the multiplicities) 
and many restrictions for $\Gamma$ will be derived in this paper.

Next we explain the main results of this paper in more detail
(we refer to section~\ref{prelim} for the relevant terminology and notation). 
Let $M$ be a connected, complete, full, irreducible isoparametric 
submanifold of rank at least $2$ in an separable infinite dimensional 
Hilbert space~$V$. By the main result in~\cite{HL}, $M$ is extrinsically 
homogeneous, but it is unknown whether the group of isometries of $V$ 
preserving $M$ is a Banach-Lie subgroup of the group of isometries of $V$,
not to speak of a reductive complement to the isotropy subalgebra. 
On the other hand, it is known from~\cite{HL} that 
there exist canonically defined one-parameter groups $\{F_X^t\}_t$
of isometries of $V$ leaving $M$ invariant, for each 
$X\in E_{\mathbf i}(x)$ with $\mathbf i \not=\mathbf 0$ and each $x\in M$. 
The restriction of
$F_X^t$ to any curvature sphere through $x$ (including 
$S_{\mathbf 0}(x)=x+E_{\mathbf 0}(x)$) is differentiable and $X$ is the
initial direction of the curve $t\mapsto F_X^t(x)$.
Thus we can define
\[  \Gamma_XY :=
(\Gamma_x)_XY :=
 \left(\frac{d}{dt}\Big|_{t=0}
 (F_X^t)_*(Y)\right)^\top \]
for all $X\in E_\mathbf i(x),\,Y\in E_\mathbf j (x)$ and 
$\mathbf i$, $\mathbf j \in \mathbf I$ with $\mathbf i \not= \mathbf 0$, 
in analogy with the finite dimensional situation. We call $\Gamma$ the \emph{(canonical) homogeneous structure} for $M$ at $x$.
The drawback of this definition is that, 
for each $X\in E_\mathbf i$, 
$\Gamma_X:T_xM\rightarrow T_xM$, or equivalently 
$\frac{d}{dt}|_{t=0}(F_X^t)_*:V\rightarrow V$,
is in principle only densely defined  
and might not be continuous and hence not extendable to the whole 
vector space. 
Geometrically speaking, it might happen that 
$(F_X^t)_*$ rotates in certain 
two-dimensional subspaces faster and faster, and analytically that the 
infinitesimal generator of $(F_X^t)_*$ is unbounded (cf.~Remark~\ref{stone}). 
That this does not occur is one of the main results of this paper.

\medskip

\noindent\textbf{Theorem A.} 
\textit{Each $\Gamma_X$ is continuous and thus extends to  $T_xM$. 
Equivalently, the one-parameter groups $\{F_X^t\}_t$ are smooth curves
in the Banach-Lie group of isometries of $V$. Moreover,
$\Gamma$ is continuous as a bilinear map.}

\medskip

An immediate consequence is that the tangent vectors
to the orbits of $F_X^t$ yield Killing vector fields defined 
on the entire Hilbert space. As a side remark, 
we point out that this result fills apparently a gap in~\cite{Ch}
where the existence of globally defined Killing fields was taken for granted. 

From Theorem~A we conclude, as in finite dimensions, that $M$ is completely determined by $\alpha_x$ and $\Gamma_x$ for any $x\in M$.
Moreover it is easily observed that $\Gamma_x$ and $(\nabla\alpha)_x$ are 
closely related. Actually we can show that they contain equivalent
information if $\alpha_x$ is given (Theorem~\ref{thm:nabla-alpha-gamma}) and therefore obtain the following rigidity result.

\medskip

\noindent\textbf{Theorem B.} 
\textit{For any $x\in M$, $\alpha_x$ and $(\nabla\alpha)_x$  
determine $M$ completely.}

\medskip

Our proof of Theorem~B also applies to finite dimensional 
homogeneous isoparametric submanifolds, and the result seems
to be even new there. In more geometric terms the theorem states that $M$ 
is completely determined by the curvature spheres at a single point $x$ 
(the information $\alpha_x$ contains) and how for each
finite dimensional curvature sphere $S_\mathbf i (x)$ any other 
curvature sphere $S_\mathbf j$ evolves along 
$S_\mathbf i (x)$ infinitesimally (the 
information $(\nabla \alpha)_x$ contains).

The proof of Theorem~A requires several steps which cover almost
the entire  paper and can be outlined as follows. We 
fix $x \in M$, $\mathbf i\in \mathbf I$ with $\mathbf i\neq\mathbf 0$
and $X_\mathbf i \in E_ \mathbf i =E_\mathbf i (x)$. By construction 
and the smoothness properties of $F_{X_\mathbf i}^t,\, \Gamma_{X_\mathbf i}$ 
is continuous on $E_\mathbf 0$. Thus it is enough to find a constant $C$ 
that depends only on $X_\mathbf i$ such that
\[ ||\Gamma_{X_{\mathbf i}}Y|| \leqq C ||Y|| \]
for all $Y$ belonging to the algebraic span of the $E_{\mathbf j}(x)$,
$\mathbf j \not=\mathbf 0$.
This problem can be split into two parts, namely to find a constant $C_1$ 
that works for all $Y$ such that $X_{\mathbf i}$ and $Y$ are both 
tangent to a finite dimensional slice through~$x$, and a 
constant $C_2$ that works for all $Y$ tangent to the 
infinite dimensional rank one slice containing $E_{\mathbf i}$. 
The first part can be easily solved as in that case 
one can estimate $\Gamma$ in the slice: since the slice is a homogeneous finite 
dimensional isoparametric submanifold and such submanifolds 
are classified, there 
is no problem to find $C_1$ which works uniformly for all slices.

The second part is more difficult. A special case, which we eventually
prove to be sufficient, is to estimate 
$||\Gamma_{X_\mathbf i}Y_{\mathbf j}||$ for all $Y_{\mathbf j}\in E_{\mathbf j}$ 
and all $\mathbf j\neq\mathbf 0$ such that $v_{\mathbf i}$, $v_{\mathbf j}$ 
are parallel. 
Essentially from the Gauss equation we obtain the 
formula
\[ \langle\Gamma_{X_{\mathbf i}}Y_{\mathbf j},\Gamma_{Y_{\mathbf j}}X_{\mathbf i}\rangle=\frac12
\langle v_{\mathbf i},v_{\mathbf j}\rangle\,||X_{\mathbf i}||^2\,||Y_{\mathbf j}||^2 \]
and the Codazzi equation allows us to interchange the two arguments of $\Gamma$ if one restricts $\Gamma_{X_{\mathbf i}}Y_{\mathbf j}$ to its components $(\Gamma_{X_{\mathbf i}}Y_{\mathbf j})_{E_\mathbf k}$ in 
$E_\mathbf k$. In fact, for all 
$\mathbf i \not= \mathbf j $
\[ (\Gamma_{Y_{\mathbf j}}X_{\mathbf i})_{E_{\mathbf k}}= c_{ijk}\,(\Gamma_{X_{\mathbf i}}Y_{\mathbf j})_{E_{\mathbf k}}, \]
where $c_{ijk}=\lambda$ if $v_\mathbf j -v_\mathbf k =\lambda (v_\mathbf i -v_\mathbf k)$ for some $\lambda$ and $c_{ijk}=0$ otherwise
(in particular, $(\Gamma_{X_{\mathbf i}}Y_{\mathbf j})_{E_{\mathbf k}}=0$ 
if $v_{\mathbf i}$, $v_{\mathbf j}$, $v_{\mathbf k}$ are not 
colinear).
The two formulae combined together yield an explicit value for
the sum over $\mathbf k$ of 
$||(\Gamma_{X_{\mathbf i}}Y_{\mathbf j})_{E_{\mathbf k}}||^2$ with certain coefficients
which, however, in general cannot be immediately used to estimate 
the length of $\Gamma_{X_\mathbf i}Y_{\mathbf j}$ as the coefficients may have different 
signs or tend to zero.
However, in the particular case that $\Gamma_{X_{\mathbf i}}Y_{\mathbf j}$ is 
contained in $E_{\mathbf 0}$ for all $\mathbf j\neq\mathbf 0$ with
$v_\mathbf j$ parallel to $v_\mathbf i$, all terms with $\mathbf k\neq\mathbf 0$
vanish yielding
\[ ||\Gamma_{X_{\mathbf i}}Y_{\mathbf j}||^2=\frac12
|| v_{\mathbf i}||^2\,||X_{\mathbf i}||^2\,||Y_{\mathbf j}||^2 \]
if $\mathbf i\not=\mathbf j$ and thus the continuity of $\Gamma_{X_{\mathbf 
i}}$.
Therefore Theorem~A follows in many cases from the next result, which 
also gives interesting information on $\Gamma$ itself.
\medskip

\noindent\textbf{Theorem C.} 
\textit{Assume the affine Weyl group of $M$ is not of type 
$\tilde B_n$ or $\tilde C_n$. Then
\[ \Gamma_{E_{\mathbf i}}E_{\mathbf j}\subset E_{\mathbf 0} \]
if $v_{\mathbf i}$ and $v_{\mathbf j}$ are 
parallel.}
 
\medskip
 
The proof of Theorem~C in turn requires several steps. Crucial ingredients
are a density theorem for the image of $\Gamma$ (Theorem~\ref{p}), a
formula for $\Gamma_{X_\mathbf i}\Gamma_{Y_\mathbf j}
Z_\mathbf k-\Gamma_{Y_\mathbf j}\Gamma_{X_\mathbf i}Z_\mathbf k$ 
which is obtained from a careful analysis of the Gauss 
equation (Corollary~\ref{stand-cpt}) and 
a simple lemma from plane geometry (Lemma~\ref{euclid}).

In the case of a general affine Weyl group, 
the statement of Theorem~C does not hold as it is and
the results about the image of $\Gamma$ 
are more technical to describe. 
To do that, we have to refine the information on 
the affine 
Weyl group by associating an affine root system to the isoparametric 
submanifold $M$ (cf.~section~\ref{section:ars}). Like for finite root systems, this is 
described by a Dynkin diagram which is obtained from the Coxeter graph of 
the affine Weyl group by attaching arrows to the double and triple links and 
additional concentric circles around those vertices that correspond to 
roots for which also twice the root is a root. These concentric circles can 
only occur in the $\tilde B_n$ or $\tilde C_n$ cases and correspond to 
reducible eigenspaces $E_\mathbf i$, i.~e.~eigenspaces which split as
$E_\mathbf i =E_\mathbf i ^\prime\oplus E_\mathbf i^{\prime\prime}$ 
under the isotropy representation of the isometry group of $M$, with 
$\dim E_\mathbf i^{\prime}$ even and
$\dim E_\mathbf i^{\prime\prime}=1$ or $3$ (cf.~section ~\ref{prelim}).

At this point it is convenient to change the notation and identify
$\mathbf I\setminus\{\mathbf 0\}$ with $\mathcal A\times\mathbb Z$, 
where $\mathcal A$ parametrizes the infinite dimensional 
rank~$1$ slices (through~$x$) and $\mathbb Z$ parametrizes, 
for each $\alpha\in \mathcal A$, the finite dimensional 
curvature distributions $E_{\alpha,i}$ in this rank~$1$ slice
in such a way that the corresponding focal hyperplanes $H_{\alpha,i}$ 
(which are parallel to each other) occur in consecutive order.
We then have the following two general results which 
include Theorem~C as special case 
and are proved by similar but slightly more refined arguments. They suffice
to finish the proof of Theorem~A in all cases.

\medskip

\noindent\textbf{Theorem D.} 
\textit{Let $(\alpha,i)\in \mathcal A\times \Z$ with $E_{\alpha,i}$ 
irreducible and $j\in \Z$. Then
\[ \Gamma_{E_{\alpha, i}}E_{\alpha, j}\subset E_{\mathbf 0} \]
unless the Dynkin diagram of $M$ is 
\[ \Cnta{}{}{}{}{}{} \quad\mbox{or}\quad \Cntar{}{}{}{}{}{}, \]
$i-j$ is even, and the hyperplane $H_{\alpha,i}$ is conjugate under the affine 
Weyl group to the focal hyperplane corresponding the right extremal 
vertex while $H_{\alpha,\frac{i+j}2}$ is conjugate to the left one.
In all cases we have
\begin{enumerate}[(i)]
 \item $\Gamma_{E_{\alpha, i}}E_{\alpha, j}\subset E_{\mathbf 0}\oplus 
 E_{\alpha, \frac{i+j}{2}} $   if $i-j$ is even.
 \item $\Gamma_{E_{\alpha, i}}E_{\alpha, j}\subset E_{\mathbf 0}$   if $i-
 j$ is divisible by $4$.
   \item $\Gamma_{E_{\alpha, i}}E_{\alpha, j}\subset E_{\mathbf 0}\oplus 
   E_{\alpha, 2i-j}\oplus E_{\alpha, 2j-i}$   
   if $i-j$ is odd.
\end{enumerate}}

\medskip

\noindent\textbf{Theorem E.} 
\textit{Let $(\alpha,i)\in\mathcal A\times \Z $ with 
$E_{\alpha,i}=E_{\alpha, i} ^\prime\oplus E_{\alpha, i}^{\prime\prime}$ 
reducible and $j\in \Z$.
Then also  $E_{\alpha,k}$ is reducible whenever $i-k$ is even and
\begin{enumerate}[(i)]
 \item $\Gamma_{E_{\alpha, i}^{\prime\prime}}E_{\alpha, j}\subset 
 E_{\mathbf 0}\oplus 
E_{\alpha, 2i-j}^\prime$   if $i-j$ is even.
\item $\Gamma_{E_{\alpha, i}}E_{\alpha, j}^{\prime\prime}\subset E_{\mathbf 
0}\oplus 
 E_{\alpha,2j-i }^\prime$   if $i-j$ is even.
\item $\Gamma_{E_{\alpha, i}^{\prime\prime}}E_{\alpha,
  j}^{\prime\prime}\subset E_{\mathbf 0}$   if $i-j$ is even.
\item $\Gamma_{E_{\alpha, i}^\prime }E_{\alpha, j}^\prime \subset 
E_{\mathbf 0}\oplus E_{\alpha, \frac{i+j}{2}}^{\prime\prime}$   if $i-j$ is 
divisible by four.
\item $\Gamma_{E_{\alpha, i}}E_{\alpha, j}\subset E_{\mathbf 0}\oplus 
E_{\alpha, 2i-j}\oplus E_{\alpha, 2j-i}$   
   if $i-j$ is odd.
\end{enumerate}}

\medskip

As an application of Theorems~D and~E, we show that 
$E_{\mathbf 0}$  is always infinite dimensional.

It is a pleasure to thank Jost Eschenburg for several helpful discussions.
Some of the questions treated here 
were also treated in the PhD thesis of K. Weinl~\cite{weinl}.

\section{Preliminaries}\label{prelim}

We recall some basic facts about isoparametric 
submanifolds of Hilbert space (cf.~\cite{Te4,PT2,HL}) 
and introduce terminology
and notation that will be used throughout the paper. 

A submanifold $M$ of an infinite dimensional separable Hilbert space
$V$ is called \emph{proper Fredholm} if the normal exponential map
$\nu M\to V$ restricted to any finite normal disk bundle is a 
proper Fredholm map. A proper Fredholm submanifold $M$ is called
\emph{isoparametric} if its normal bundle is globally flat, and the 
shape operators along any parallel normal vector field are conjugate. 
Here globally flat means that every normal vector can be uniquely 
extended to a parallel normal vector field along the 
whole of $M$. The Fredholm condition implies that the codimension of
$M$ is finite and its shape operators are 
compact (self-adjoint) operators. 
Since the normal bundle $\nu M$ is flat, the Ricci equation
yields a splitting of the tangent bundle as $TM=\oplus_{{\mathbf i}\in \mathbf I}E_{\mathbf i}$ (closure of the algebraic direct sum) 
into the simultaneous eigendistributions $E_{\mathbf i}$ of the shape operators,
where $\mathbf I$ is a countable index set containing $\mathbf 0$; each $E_{\mathbf i}$ is 
called a \emph{curvature distribution}.
For each normal vector
$\xi\in\nu M$, the corresponding shape operator satisfies 
$A_\xi|_{E_{\mathbf i}}=\langle\xi,v_{\mathbf i}\rangle\mathrm{id}_{E_{\mathbf i}}$, where $v_{\mathbf i}$ 
is a globally defined parallel normal vector field on $M$; each 
$v_{\mathbf i}$ is called a \emph{curvature normal}. Unless explicitly
stated, we will denote the zero curvature normal by $v_{\mathbf 0}$ (if it 
occurs) and the corresponding curvature distribution by $E_{\mathbf 0}$. 
For convenience, we also set $\mathbf I^*=\mathbf I\setminus\{\mathbf 0\}$.
Note that the substantial codimension of $M$ equals the number 
of linearly independent curvature normals; this number is called
the \emph{rank} of $M$. We will always assume that $M$ is full in $V$,
that is, not contained in a proper affine subspace. It then follows
that the curvature normals of $M$ span the normal space. 

As a consequence of the Codazzi equations, each curvature 
distribution $E_{\mathbf i}$ is integrable with totally geodesic 
leaves.
The leaf of $E_{\mathbf i}$ passing through $x$ is denoted by $S_{\mathbf i}(x)$.
Due to the compactness of the shape operators, the dimension
of $E_{\mathbf i}$ is finite if $\mathbf i\neq\mathbf 0$, and then $m_{\mathbf i}=\dim E_{\mathbf i}$
is called a \emph{multiplicity}, but $E_{\mathbf 0}$ can have infinite 
dimension. Since the leaves of the $E_{\mathbf i}$ are umbilic in $V$
(by definition), it follows that 
$S_{\mathbf i}(x)$ is a round sphere for $\mathbf i\neq\mathbf 0$ 
(centered at $c_{\mathbf i}(x):=x+(v_{\mathbf i}(x)/||v_{\mathbf i}||^2)$,
with radius $1/||v_{\mathbf i}||$) and $S_{\mathbf 0}(x)=x+E_{\mathbf 0}(x)$ is a closed
affine subspace. 

We will always assume that $M$ is complete. A complete 
isoparametric submanifold $M$ of $V$ determines a singular 
foliation of $V$ by parallel submanifolds, where the regular leaves 
are also isoparametric of the same codimension as $M$ and 
the singular leaves are the focal manifolds of $M$. 
Indeed a parallel submanifold of $M$, denoted by $M_\xi$,
is determined by a parallel
normal vector field $\xi$ along $M$ so that the map
$\pi_\xi:x\mapsto x+\xi(x)$ from $M$ to $M_\xi$ is a submersion.
This map has differential $\mathrm{id}-A_\xi$ so its kernel 
is $\oplus\{E_{\mathbf i}\;|\;\mbox{${\mathbf i}\in \mathbf I$, $\langle\xi,v_{\mathbf i}\rangle=1$}\}$.
We thus see that the focal set of $M$ decomposes
into focal manifolds, and $M_\xi$ is a focal manifold precisely 
if $\ker d\pi_\xi$
is not zero in which case $\pi_\xi$ is called a \emph{focal map}.
For $x\in M$, the affine
normal space $x+\nu_xM$ meets the focal set of $M$ along the union 
over ${\mathbf i}\in \mathbf I^*$ of the affine hyperplanes
$H_{\mathbf i}(x):=x+\{\xi\in\nu_xM\;|\;
\langle\xi,v_{\mathbf i}\rangle=1\}$, called
\emph{focal hyperplanes} (with respect to~$x$).
An important result of Terng~\cite{Te4}
states that the group $W$ which is generated by reflections
in the $H_{\mathbf i}(x)$, ${\mathbf i}\in \mathbf I^*$, is an 
affine Weyl group acting on $x+\nu_xM$.
We will always assume that $M$ is irreducible, i.e.~it cannot
be split as an extrinsic product of lower dimensional
isoparametric submanifolds. It then follows that 
$W$ acts irreducibly on $x+\nu_xM$ (cf.~\cite{HL2}) 
and thus is isomorphic to one of $\tilde A_n$, $\tilde B_n$, 
$\tilde C_n$, $\tilde D_n$, $\tilde E_6$,
$\tilde E_7$,  $\tilde E_8$, $\tilde F_4$, $\tilde G_2$. 
It follows that the set of focal hyperplanes decomposes
into finitely many families of parallel equidistant
hyperplanes in $x+\nu_xM$. For each family, the corresponding
curvature normals are thus of the form $v_k=(d_0+kd)^{-1}v$,
$k\in\Z$, where $v$ is a unit vector, $d$ is the distance between
two consecutive hyperplanes and $d_0$ is the distance from $x$
to the first hyperplane of this family in the direction of $v$.
It also follows that there are only finitely many different
multiplicities, since those are preserved by the action of $W$
on the set of hyperplanes. 
Later in the paper it will be convenient to identify $\mathbf I^*$
with $\mathcal A\times \Z$, where $\mathcal A$ is a finite index set parametrizing 
the families of parallel curvature normals and, for each 
$\alpha\in\mathcal A$, $\Z$ parametrizes the curvature normals in that 
family so that the corresponding focal hyperplanes 
$\{H_{\alpha,i}(x)\}_{i\in\mathbb Z}$ 
are in consecutive order. Whereas typical indices in $\mathbf I$ 
will be denoted by $\mathbf i$, $\mathbf j$, $\mathbf k$...,
typical indices in $\mathcal A\times \Z$ will be denoted 
by $(\alpha,i)$, $(\beta,j)$, $(\gamma,k)$,...

Since the second fundamental form $\alpha$ of $M$ is related
to the shape operators by 
$\langle \alpha(X,Y),\xi\rangle=\langle A_\xi X,Y\rangle$ for 
$X$, $Y\in T_xM$, $\xi\in\nu_xM$, it follows that 
\begin{equation}\label{alpha}
\alpha(X_{\mathbf i},Y_{\mathbf j})=\langle X_{\mathbf i},Y_{\mathbf j}\rangle\,v_{\mathbf i}
\end{equation}
for $X_{\mathbf i}\in E_{\mathbf i}(x)$, $Y_{\mathbf j}\in E_{\mathbf j}(x)$. 
Hence the focal hyperplanes in $x+\nu_xM$ 
together with the multiplicities $m_{\mathbf i}=\dim E_{\mathbf i}$, ${\mathbf i}\in \mathbf I^*$, essentially
determine the second fundamental form, up to
passing to a parallel isoparametric submanifold. 
In turn the focal hyperplanes
are already determined by the affine Weyl group, up to 
scaling of the ambient metric.
Thus the affine Weyl group together
with the multiplicities essentially determine the second
fundamental form. Such data is usually encoded in a 
Coxeter graph with multiplicities (cf.~section~\ref{section:ars} for the 
Coxeter graph).

Another fundamental invariant of $M$ is the covariant 
derivative of the second fundamental form. By taking derivatives,
it follows from~(\ref{alpha}) that 
\begin{equation}\label{nabla-alpha}
 \nabla_{X_{\mathbf i}}\alpha(Y_{\mathbf j},Z_{\mathbf k})=
\langle \nabla_{X_{\mathbf i}}Y_{\mathbf j},Z_{\mathbf k}\rangle(v_{\mathbf j}-v_{\mathbf k}) 
\end{equation}
for $\mathbf i$, $\mathbf j$, $\mathbf k\in \mathbf I$ and $X_{\mathbf i}\in E_{\mathbf i}(x)$, $Y_{\mathbf j}\in E_{\mathbf j}(x)$, 
$Z_{\mathbf k}\in E_{\mathbf k}(x)$. One uses the Codazzi equation, which is the symmetry
of $\nabla\alpha$ in all three arguments, to derive strong
restrictions on $M$. For instance if $P$ is an affine subspace
of $x+\nu_xM$ then $\mathcal D_P(x):=\oplus_{v_{\mathbf i}(x)\in P}E_{\mathbf i}(x)$
(closure of algebraic direct sum)
is an integrable distribution with totally geodesic leaves. 
Moreover the leaf through $x\in M$, to be denoted by $L_P(x)$, 
is a complete full isoparametric submanifold of the affine
subspace 
$W_P(x):=x+\mathcal D_P(x)+\mathrm{span}\{v_{\mathbf i}(x)\;|\;v_{\mathbf i}(x)\in P\}$ 
of rank $\dim P$ if $0\in P$ and $\dim X+1$ if $0\not\in P$. 
By taking in particular  $P=\{v_{\mathbf i}\}$ for some ${\mathbf i}\in \mathbf I$ we see
that $E_{\mathbf i}$ is integrable and $S_{\mathbf i}(x)$ is totally geodesic in $M$
as mentioned above. In general $L_P(x)$ is called
a \emph{slice} through $x$. It is finite dimensional precisely if
$0\not\in P$. In this case it can be focalized, that is,
there exists a parallel normal vector field $\xi_P$ with the 
property that  
$\langle\xi_P,v_{\mathbf i}\rangle=1$ if and only if $v_{\mathbf i}(x)\in P$, so that the kernel
of the differential $(\pi_{\xi_P})_*$ of the focal map
is $\mathcal D_P$ and the connected components
of the fibers of $\pi_{\xi_P}$ are the leaves $L_P$. 
Note also in this case that $L_P(x)$ is contained in the round hypersphere 
in $W_P(x)$ of center $c_P(x):=x+\xi_P(x)$ and radius $||\xi_P||$. 
On the other hand,
if we start with a point $y$ in $x+\nu_xM$ lying in some
focal hyperplanes, then there is a unique
finite dimensional 
slice $L_P(x)$ focalizing at $y$ coming from $P=\mathrm{affine\
span}\{v_{\mathbf i}(x)\;|\;y\in H_{\mathbf i}(x)\}$. In this case
$\mathcal D_P=\oplus_{y\in H_{\mathbf i}(x)}E_{\mathbf i}$, namely, $P$ 
does not contain any other curvature normals besides 
those $v_{\mathbf i}(x)$ with $y\in H_{\mathbf i}(x)$. 

In addition to the above assumptions (connectedness,
completeness, fullness, irreducibility) we henceforth
assume that $M$ has rank at least~$2$. Then the main result
of~\cite{HL} asserts that $M$ is homogeneous. 
In order to explain this result we first recall that 
for each curvature sphere $S_{\mathbf i}(x)$, ${\mathbf i}\in \mathbf I^*$, there 
exists a distinguished compact connected Lie group $\Phi_{\mathbf i}^*=\Phi_{\mathbf i}^*(x)$ 
which acts isometrically on the affine span $W_{\mathbf i}(x)$ of 
$S_{\mathbf i}(x)$ and has $S_{\mathbf i}(x)$ as an orbit. In fact $\Phi_{\mathbf i}^*$ is 
the identity component of the normal holonomy group 
of the focal manifold obtained from $M$ by focalizing the distribution
$E_{\mathbf i}$. The action of $\Phi_{\mathbf i}^*$ on $W_{\mathbf i}(x)$ is equivalent to the 
isotropy representation of a symmetric space of rank $1$
different from the Cayley projective plane,
where we view $W_{\mathbf i}(x)$ as a vector space with the origin
at the center $c_{\mathbf i}(x)$ of $S_{\mathbf i}(x)$. Hence
$\Phi_{\mathbf i}^*$ is isomorphic to one of 
$SO(m_{\mathbf i}+1)$, $U(\frac{m_{\mathbf i}+1}2)$ ($m_{\mathbf i}$ odd) or 
$(Sp(\frac{m_{\mathbf i}+1}4)\times Sp(1))/\Z_2$ ($m_1\equiv3\mod4$).
Note that the isotropy group of~$\Phi_{\mathbf i}^*$ at~$x$
acts irreducibly on $E_{\mathbf i}(x)$ in the 
first case only and that otherwise $E_{\mathbf i}(x)$ 
splits into two irreducible subspaces as 
$E_{\mathbf i}(x)=E_{\mathbf i}'(x)\oplus E_{\mathbf i}''(x)$
with $\dim E_{\mathbf i}'(x)$ even and $\dim E_{\mathbf i}''(x)=1$ or~$3$.
It is useful to note that the group 
$\Phi_{\mathbf i}^*$ is already determined by any irreducible finite dimensional 
slice of rank $2$ containing $S_{\mathbf i}(x)$, and it follows
that $\Phi_{\mathbf i}^*\cong  SO(m_{\mathbf i}+1)$ unless $W$ is isomorphic
to $\tilde B_n$ or $\tilde C_n$. 
We also mention that the construction of $\Phi_{\mathbf i}^*$
can be generalized by replacing $S_{\mathbf i}(x)$ by any slice $L_P(x)$
to yield a compact connected Lie group $\Phi_P^*(x)$ acting isometrically
on $W_P(x)$ and having $L_P(x)$ as an orbit. If $0\not\in P$,
the group  $\Phi_P^*(x)$ is the normal holonomy group of a 
focal manifold and acts on $W_P(x)$ as the isotropy
representation of a symmetric space of rank equal to the 
rank of $L_P(x)$ as an isoparametric submanifold 
(viewing $W_P(x)$ as a vector space with origin at ~$c_P(x)$). 
The later statement is also known as the Homogeneous Slice Theorem~\cite{HL}.

The action of the group $\Phi_{\mathbf i}^*(x)$ on $S_{\mathbf i}(x)$
induces in a natural way an invariant connection on $S_{\mathbf i}(x)$ 
which coincides with the Levi-Civit\`a connection 
if $\Phi_{\mathbf i}^*\cong SO(m_{\mathbf i}+1)$. It is defined via a reductive
decomposition on the Lie algebra level 
$\mathbf L(\Phi_{\mathbf i}^*)=\mathbf L((\Phi_{\mathbf i}^*)_x)+\Lm_{\mathbf i}$ where
$(\Phi_{\mathbf i}^*)_x$ denotes the isotropy group of $\Phi_{\mathbf i}^*$ at $x$ and
$\Lm_{\mathbf i}$ is the orthogonal complement of 
$\mathbf L((\Phi_{\mathbf i}^*)_x)$
with respect the inner product 
$\langle X,Y\rangle=-B(X,Y)-\mathrm{trace}(\rho_*X)(\rho_*Y)$ for $X$, $Y\in\mathbf L(\Phi_{\mathbf i}^*)$,
where $B$ denotes the Killing form of $\mathbf L(\Phi_{\mathbf i}^*)$
and $\rho$ denotes its representation on $W_{\mathbf i}(x)$. 
The invariant connection of course has the property 
that $\gamma_X(t):=\exp tX(x)$ is a geodesic for any $X\in\Lm_{\mathbf i}$ 
and $(\exp tX)_*$ induces parallel translation along $\gamma_X$. 

Finally
the proof of the homogeneity of $M$ is based on the explicit
construction of certain one-parameter groups $F_X^t$ of isometries
of $V$ leaving $M$ invariant. In fact for each $x\in M$, ${\mathbf i}\in \mathbf I^*$ 
and $X\in E_{\mathbf i}(x)$, there exists such a one-parameter group 
having the following properties, by which they are also 
determined. 
\begin{enumerate}
\item $F_X^t(x):=\gamma_X(t)$ is the geodesic in $S_{\mathbf i}(x)$ 
with initial speed $X$ with respect to the above invariant connection
and $(F_ X^t)_*Y$ is a parallel vector field along $\gamma_X$ for any 
$Y\in E_{\mathbf i}(x)$, i.e.~$(F_X^t)_*|_{W_{\mathbf i}(x)}=(\exp t\check 
X)_*$ where 
$\check X$ is the unique vector in $\Lm_{\mathbf i}$ such that 
$\gamma_X(t)=\exp t\check X(x)$ has initial speed $X$. 
\item For each ${\mathbf j}\in \mathbf I\setminus\{\mathbf i\}$ and $y\in 
S_{\mathbf j}(x)$,  
$F_X^t(y)$ is the unique smooth curve which is
everywhere orthogonal to $E_{\mathbf j}$ 
and satisfies $F_X^t(y)\in S_{\mathbf j}(\gamma_X(t))$ for all $t$. 
Moreover
$(t,y)\mapsto F_X^t(y)$ is smooth on $\R\times S_{\mathbf j}(x)$. 
\item $(F_X^t)_*\xi$ is parallel in $\nu M$ along 
$\gamma_X$ for all $\xi\in\nu_xM$. 
\end{enumerate}
By means of these one-parameter groups one proves
that the group of isometries of $V$ preserving $M$ is transitive
on $Q_x$, for $x\in M$, where $Q_x$ is defined as the set 
of points of $M$ that can be reached from $x$ by following
piecewise differentiable curves along finite dimensional
curvature spheres. One finally deduces the homogeneity of $M$
from the results that this group is also transitive on the closure of $Q_x$ 
in~$M$~\cite[Prop.~4.4]{HL} and that $Q_x$ is dense 
in $M$~\cite[Thm.~B]{HL}. 

\emph{In this paper, $V$ will always denote a 
finite or infinite dimensional
separable Hilbert space, and $M$ 
a connected, complete, full and irreducible isoparametric submanifold of $V$
with rank at least $2$. In case $M$ is finite dimensional and 
of rank $2$, we assume in addition that $M$ is homogeneous.}
 Thus 
$M$ is always homogeneous, and more precisely, an orbit of its 
stabilizer in the isometry group of~$V$. 
Moreover, for $x\in M$, ${\mathbf i}\in \mathbf I^*$ and 
$X\in E_{\mathbf i}(x)$ we have 
the one-parameter groups $F_ X^t$ leaving $M$ invariant 
described above (in the finite dimensional case, 
this follows from~\cite[Lemma~4.2]{HL}). In the finite dimensional
case we also recall that $E_{\mathbf 0}=0$ by irreducibility.
Although in this paper 
we are ultimately interested in the infinite dimensional case, 
the reason to also treat finite dimensional isoparametric submanifolds is that they appear as (necessarily homogeneous) slices of the infinite dimensional 
ones. 

\section{The homogeneous structure $\Gamma$}
\setcounter{equation}{0}

\subsection{Definition of $\Gamma$ and basic properties}
We now introduce the fundamental object of this paper.

\begin{defn}\label{gamma} 
For each $x\in M$, ${\mathbf i}\in \mathbf I^*$, ${\mathbf j}\in \mathbf I$ and $X_{\mathbf i}\in E_{\mathbf i}(x)$, $Y_{\mathbf j}\in E_{\mathbf j}(x)$
let 
\[ (\Gamma_x)_{X_{\mathbf i}}Y_{\mathbf j}:=\left(\frac{d}{dt}\Big |_{t=0}(F_{X_{\mathbf i}}^t)_*(Y_{\mathbf j})\right)^\top \]
where $(\cdot)^\top$ denotes the tangential part. We call
$\Gamma$ the \emph{homogeneous structure} of $M$. 
\end{defn}

If $x$ is fixed we write $\Gamma$ instead of $\Gamma_x$
and $E_{\mathbf i}$ instead of $E_{\mathbf i}(x)$.

\begin{rem}\label{stone}
By a theorem of Stone (see e.~g.~\cite[p.~327]{conway}), 
for a given one-parameter group $\varphi_t$ of unitary automorphisms of a 
complex Hilbert space, there exists a self-adjoint operator $A$ such that 
$\varphi_t=\exp itA$ if and only if $\varphi_t$ is strongly continuous, 
that is, in case $t\mapsto\varphi_tv$ is continuous 
for all $v$. The domain of definition of $A$ is then the 
set of vectors $v$ such that $\lim_{t\to0}\frac{\varphi_t v-v}{t}$
exists and $iAv$ is this limit. The operator 
$A$ is called the infinitesimal generator of $\varphi_t$. 
In our setting, Stone's theorem can be applied to the 
complexification of the one-parameter group $(F_{X_{\mathbf i}}^t)_*$.
In fact, we know that $(F_{X_{\mathbf i}}^t)_*$ is strongly continuous on a dense 
subset of $V$ (namely, on the algebraic span of $\nu_xM$ and the $E_{\mathbf i}(x)$ for 
${\mathbf i}\in \mathbf I$) and thus strongly continuous on $V$ (since the one-parameter
group $(F_{X_{\mathbf i}}^t)_*$ consists of isometries). 
Thus $\Gamma_{X_{\mathbf i}}$ 
is essentially the infinitesimal generator of $(F_{X_{\mathbf i}}^t)_*$. 
In particular
the continuity of $\Gamma_{X_{\mathbf i}}$ is equivalent to the domain of 
this infinitesimal generator being the entire Hilbert space. 
\end{rem}

\begin{lem}\label{elem-gamma}
\begin{enumerate}
\item[(i)] $\frac{d}{dt}\big|_{t=0}(F_{X_{\mathbf i}}^t)_*Y_{\mathbf j} = \Gamma_{X_{\mathbf i}}Y_{\mathbf j}
+\alpha(X_{\mathbf i},Y_{\mathbf j})$
\item[(ii)] $\frac{d}{dt}\big|_{t=0}(F_{X_{\mathbf i}}^t)_*\xi=-A_{\xi}X_{\mathbf i}$
\end{enumerate}
for ${\mathbf i}\in \mathbf I^*$, ${\mathbf j}\in \mathbf I$ and $X_{\mathbf i}\in E_{\mathbf i}$, $Y_{\mathbf j}\in E_{\mathbf j}$, $\xi\in\nu_xM$.
\end{lem}

\Pf Part (ii) follows from the definition of $F^t_{X_{\mathbf i}}$, and 
part~(i) follows from~(ii) by taking inner product with~$Y_{\mathbf j}$.  \EPf

\begin{lem}\label{gamma-skew-and-invariant}
\begin{enumerate}
\item[(i)] $\Gamma_X$ is skew-symmetric: $\langle\Gamma_{X_{\mathbf i}}Y_{\mathbf j},Z_{\mathbf k}\rangle+
\langle Y_{\mathbf j},\Gamma_{X_{\mathbf i}}Z_{\mathbf k}\rangle=0$ for all ${\mathbf i}\in \mathbf I^*$, $\mathbf j$, $\mathbf k\in\mathbf I$,
and $X_{\mathbf i}\in E_{\mathbf i}$, $Y_{\mathbf j}\in E_{\mathbf j}$, $Z_{\mathbf k}\in E_{\mathbf k}$.
\item[(ii)] $\Gamma$ is invariant under isometries: 
$\Gamma_{g_*X_{\mathbf i}}g_*Y_{\mathbf j}|_{gx}=g_*\Gamma_{X_{\mathbf i}}Y_{\mathbf j}|_x$ 
for any extrinsic isometry $g$ of $M$,
$X_{\mathbf i}\in E_{\mathbf i}$, $Y_{\mathbf j}\in E_{\mathbf j}$ and ${\mathbf i}\in \mathbf I^*$, ${\mathbf j}\in\mathbf I$. 
\end{enumerate}
\end{lem}

\Pf Part~(i) follows from $\langle (F^t_{X_{\mathbf i}})_*Y_{\mathbf j},(F^t_{X_{\mathbf i}})_*Z_{\mathbf k}\rangle
=\langle Y_{\mathbf j},Z_{\mathbf k}\rangle $ by taking derivative at $t=0$. 

Part (ii) follows from $F_{g_*X_{\mathbf i}}^t=gF^t_{X_{\mathbf i}}g^{-1}$ which in turn 
is a direct consequence of the definition of $F_{X_{\mathbf i}}^t$. \EPf 

\begin{lem}\label{gamma-restriction}
\begin{enumerate}
\item[(i)]
Let $L$ be a slice 
of $M$ through $x$, choose ${\mathbf i}\in\mathbf I$ with
$E_{\mathbf i}(x)\subset T_xL$ and~$X\in E_{\mathbf i}(x)$.  
Then $F_X^t(L)=L$.
\item[(ii)] If, in addition, $L$ is irreducible
and has rank at least~$2$,
$W$ denotes its affine span in~$V$, and  
${}^L\!F^t_X$ the one-parameter group of isometries of $W$ 
associated to $X$, then 
\[ F_X^t|_W={}^L\!F_X^t \]
for all $t$, and the homogeneous structure of $L$ is the restriction of that of $M$.
\end{enumerate}
\end{lem}

\Pf (i) We have $L=L_P(x)$ for some affine subspace $P\subset\nu_xM$
and $F_X^t(L_P(x))=L_P(F_X(x))=L_P(x)$ since $F_X^t(x)\in L_P(x)$. 

(ii)The result is a consequence of~\cite[Lemma~1.4]{HL},
and the fact that the group $\Phi_{\mathbf i}^*$ which acts 
transitively on $S_{\mathbf i}(x)$ is already determined by $L$. The last assertion
follows from~\cite[Remark, p.~162]{HL} in case
$L$ is finite dimensional, and otherwise by applying the 
same remark twice to the 
inclusions $S_{\mathbf i}(x)\subset L' \subset M$ and 
$S_{\mathbf i}(x)\subset L'\subset L$, where $L'$ is a finite
dimensional rank $2$ slice. \EPf

\subsection{The homogeneous structure of parallel isoparametric submanifolds}

Let $\xi$ be a parallel normal vector field along $M$ and denote
by $\pi:M\to M_\xi$ the endpoint map $x\mapsto x+\xi(x)$. 
Assume that $M_\xi$ is also isoparametric
and denote its homogeneous structure by ${}^\xi\Gamma$. 

\begin{lem}\label{parallel}
For all ${\mathbf i}\in \mathbf I^*$, ${\mathbf j}\in\mathbf I$, $X_{\mathbf i}\in E_{\mathbf i}$ and $Y_{\mathbf j}\in E_{\mathbf j}$:
\begin{enumerate}
\item[(i)] ${}^\xi\Gamma_{\pi_*X_{\mathbf i}}Y_{\mathbf j}\vert_{\pi(x)}=\Gamma_{X_{\mathbf i}}Y_{\mathbf j}\vert_x$;
\item[(ii)] $\pi_*{X_{\mathbf i}}=\pm\frac{||v_{\mathbf i}||}{||v_{\mathbf i}^\xi||}X_{\mathbf i}$, where 
$v_{\mathbf i}^\xi$ is the curvature normal of $M_\xi$ with respect to $\pi_*E_{\mathbf i}$
(this subspace is a curvature distribution of $M_\xi$ that
coincides with $E_{\mathbf i}$ up to parallel translation in $V$,
but may correspond to a different index in $\mathbf I^*$ if $M_\xi=M$). 
\end{enumerate}
\end{lem}

\Pf (i) The diffeomorphism $\pi$ maps curvature spheres
of $M$ through $x$ to curvature spheres of $M_\xi$ through $\pi(x)$.
In fact $\pi_*=\mathrm{id}-A_\xi$ maps curvature distributions
to curvature distributions and actually preserves $E_{\mathbf k}(x)$ 
as a subspace of $V$ for all $\mathbf k\in\mathbf I$.  
Since $F_{X_{\mathbf i}}^t$ is an isometry of $V$ preserving $M$
and inducing parallel transport along $\nu M$, it immediately
follows that $F_{X_{\mathbf i}}^t$ commutes with $\pi$ and in particular
preserves the parallel submanifold $M_\xi$. Note that the initial
speed of $t\mapsto F_{X_{\mathbf i}}^t(\pi(x))$ is $\pi_* X_{\mathbf i}$.  
We claim that $F_{X_{\mathbf i}}^t$ is the canonical one-parameter
group of isometries of $V$ preserving $M_\xi$ 
associated to $\pi_*X_{\mathbf i}\in T_{\pi(x)}M_\xi$.
The result then follows by  
differentiating at $t=0$ and taking tangential parts. 
In turn, the claim is proved by checking~(a), (b) and~(c)
in section~\ref{prelim}, and these conditions follow from the fact
that $M$ and $M_\xi$ have parallel curvature distributions. 

(ii) This follows from the fact that $\pi$ maps curvature spheres of $M$
through $x$ of radii $1/||v_{\mathbf i}||$ 
to curvature spheres of $M_\xi$ through $\pi(x)$
of radii $1/||v^\xi_{\mathbf i}||$. \EPf

\subsection{Finite dimensional case}\label{fin-dim}

The case of finite dimensional
slices is particularly interesting, since they are congruent 
to principal orbits of isotropy representations
of symmetric spaces, for which we have an effective 
way of computing the homogeneous structure. 

Let $M$ be a principal orbit of the isotropy
representation of an irreducible 
symmetric space $G/K$ of rank at least~$2$, say of noncompact 
type, with $K$ connected.
Let $\Lg=\Lk+\Lp$ be the decomposition of the Lie algebra $\Lg$ 
of $G$ into the eigenspaces of the involution.
Equip $\Lp$ with an $\mathrm{Ad}_{\mathfrak k}$-invariant 
inner product. 
Choose a maximal Abelian subspace $\La$ of $\Lp$,
Then $M$ is the adjoint orbit $\mathrm{Ad}_{\mathfrak k}(x)$ for
some regular point $x\in\La$ and $\nu_xM=\La$.
Consider the usual root space decompositions 
$\Lk=\Lk_0+\sum_{\lambda\in\Lambda}\Lk_\lambda$ and 
$\Lp=\La+\sum_{\lambda\in\Lambda}\Lp_\lambda$ where 
$\Lambda$ denotes the root system with respect to $\La$ and
$\Lk_\lambda=\Lk_{-\lambda}$, $\Lp_\lambda=\Lp_{-\lambda}$.
Then $\Lp_+:=\sum_{\lambda\in\Lambda}\Lp_\lambda=T_xM$, 
the index set $\mathbf I$ can be identified with the set 
$\Lambda^+_{\mathrm {red}}$
of positive roots
$\lambda$ such that $\frac12\lambda$ is not a root, and for 
such a root we have $E_\lambda=\Lp_\lambda\oplus\Lp_{2\lambda}$, where 
we use the convention that $\Lp_{2\lambda}=0$ if $2\lambda\not\in\Lambda$. 
Thus 
$M$ carries the canonical homogeneous structure given by the reductive
complement $\Lk^+=\sum_{\lambda\in\Lambda}\Lk_\lambda$, 
and this coincides with the homogeneous 
structure in our definition. In fact, in the proof of~\cite[Lemma~4.2]{HL}
it was shown that for any $X\in\Lp_\lambda\oplus\Lp_{2\lambda}$, 
$F_X^t$ equals $\exp t\check X$, where $\check X\in\Lk_\lambda+\Lk_{2\lambda}$ 
is the unique element in $\Lk^+$ satisfying $[\check X,x]=X$.  
Hence  
\begin{equation}\label{gamma-bracket}
\Gamma_XY=[\check X,Y]_{\mathfrak p^+}
\end{equation}
(component in $\Lp_+$) for $X$, $Y\in T_xM=\Lp^+$.

As an application of the above discussion, we prove a
result that will be used in the proof of 
Theorem~\ref{thm:nabla-alpha-gamma}.

\begin{prop}
If $v_{\mathbf i}$, $v_{\mathbf j}$, $v_{\mathbf k}$ are pairwise different
and they span an affine subspace $P$ in $\nu_xM$ that does not contain $0$, then
\begin{equation}\label{gamma-homo}
 \Gamma_{X_{\mathbf i}}\Gamma_{Y_{\mathbf j}}Z_{\mathbf k}-\Gamma_{Y_{\mathbf j}}\Gamma_{X_{\mathbf i}}Z_{\mathbf k}
=\Gamma_{\Gamma_{X_{\mathbf i}}Y_{\mathbf j}-\Gamma_{Y_{\mathbf j}}X_{\mathbf i}}Z_{\mathbf k} 
\end{equation}
where $X_{\mathbf i}\in E_{\mathbf i}(x)$, $Y_{\mathbf j}\in E_{\mathbf j}(x)$, $Z_{\mathbf k}\in E_{\mathbf k}(x)$. 
\end{prop}

\Pf The assumption that $v_{\mathbf i}$, $v_{\mathbf j}$, $v_{\mathbf k}$ span $P$
allows us to restrict to the corresponding finite dimensional
slice, which can be assumed to be irreducible 
and is thus congruent to a principal orbit
of the isotropy representation of a symmetric space, so
formula~(\ref{gamma-bracket}) can be used. The Jacobi identity
gives 
\begin{equation}\label{gamma-jacobi}
 [\check X_{\mathbf i},[\check Y_{\mathbf j},Z_{\mathbf k}]]+[\check Y_{\mathbf j},[Z_{\mathbf k},\check X_{\mathbf i}]]
+ [Z_{\mathbf k},[\check X_{\mathbf i},\check Y_{\mathbf j}]] = 0. 
\end{equation}
Note that $\Gamma_{Y_{\mathbf j}}Z_{\mathbf k}=[\check Y_{\mathbf j},Z_{\mathbf k}]_{\mathfrak p_+}=[\check Y_{\mathbf j},Z_{\mathbf k}]$
since $v_{\mathbf j}$, $v_{\mathbf k}$ are different. 
Taking $\Lp_+$-components in~(\ref{gamma-jacobi}), 
we now see that the first two terms in the formula thus obtained
correspond to the left hand side of~(\ref{gamma-homo}), and it remains to
see that $-[Z_{\mathbf k},[\check X_{\mathbf i},\check Y_{\mathbf j}]]_{\mathfrak p_+}$ corresponds to
the right hand side. In fact, this also follows from Jacobi, namely
\begin{eqnarray*}
[[\check X_{\mathbf i},\check Y_{\mathbf j}],x]&=&[X_{\mathbf i},\check Y_{\mathbf j}]+[\check X_{\mathbf i},Y_{\mathbf j}] \\
                           &=&\Gamma_{X_{\mathbf i}}Y_{\mathbf j}-\Gamma_{Y_{\mathbf j}}X_{\mathbf i} 
\end{eqnarray*}
using that $v_{\mathbf i}\neq v_{\mathbf j}$, so
$[[\check X_{\mathbf i},\check Y_{\mathbf j}],Z_{\mathbf k}]_{\mathfrak p_+}=\Gamma_{\Gamma_{X_{\mathbf i}}Y_{\mathbf j}-\Gamma_{Y_{\mathbf j}}X_{\mathbf i}}Z_{\mathbf k}$. \EPf

\subsection{First results on the image of $\Gamma$}
We start with some basic results about the components 
of $\Gamma_{E_{\mathbf i}}E_{\mathbf j}$ in the eigenspaces $E_{\mathbf k}$. 
Later in section~\ref{image} they will be considerably refined. 

\begin{prop}\label{gamma-ei-ei}
Let ${\mathbf i}\in \mathbf I^*$. 
\begin{enumerate}
\item[(i)] If $E_{\mathbf i}$ is irreducible then $\Gamma_{E_{\mathbf i}}E_{\mathbf i}=0$.
\item[(ii)] If $E_{\mathbf i}$ is not irreducible then 
$\Gamma_{E_{\mathbf i}''}E_{\mathbf i}''=0$,
$\Gamma_{E_{\mathbf i}'}E_{\mathbf i}''\subset E_{\mathbf i}'$, $\Gamma_{E_{\mathbf i}''}E_{\mathbf i}'\subset E_{\mathbf i}'$
and $\Gamma_{E_{\mathbf i}'}E_{\mathbf i}'\subset E_{\mathbf i}''$. 
\end{enumerate}
\end{prop}

\Pf (i) Note that $\Gamma_{E_{\mathbf i}}E_{\mathbf i}\subset E_{\mathbf i}$ simply by the definition
of $F^t_{X_{\mathbf i}}$ for $X_{\mathbf i}\in E_{\mathbf i}$. The irreducibility of $E_{\mathbf i}$ means
$\Phi_{\mathbf i}^*\cong SO(m_{\mathbf i}+1)$, where $m_{\mathbf i}=\dim S_{\mathbf i}(x)$. Thus the 
connection induced by $\Phi_{\mathbf i}^*$ on $S_{\mathbf i}(x)$ coincides with the 
standard Levi-Civit\`a connection. Since $(F_X^t)_*Y$ is parallel
along $t\mapsto F_X^t(x)$ in $S_{\mathbf i}(x)$ for all $X$, $Y\in E_{\mathbf i}(x)$,
$\Gamma_XY=\left(\frac{d}{dt}\big| _{t=0}(F_X^t)_*Y)\right)^\top=0$. 
Alternatively, we could also prove this result by employing a 
line of reasoning like in part~(ii). 

(ii) We use Lemma~\ref{gamma-restriction}.
Let $L$ be an irreducible finite dimensional rank two slice through $x$
which contains $E_{\mathbf i}$. Then $L$ is congruent to a principal
orbit of the isotropy representation of a symmetric space 
of type $BC_2$, and the statements follow from 
the bracket relations
in the corresponding Lie algebra and
formula~(\ref{gamma-bracket}). \EPf

\begin{prop}\label{gamma-perp}
We have $\Gamma_{E_{\mathbf i}}E_{\mathbf j}\perp E_{\mathbf j}$ for ${\mathbf i}\in \mathbf I^*$, ${\mathbf j}\in\mathbf I$ and 
$\mathbf i\neq \mathbf j$. 
\end{prop}

\Pf Let $X_{\mathbf i}\in E_{\mathbf i}$, $Y_{\mathbf j}\in E_{\mathbf j}$ and let $c(s_1,s_2)$  
be a differentiable parametrized surface in $S_{\mathbf j}(x)$ with $c(0,0)=x$  
and $\frac{\partial c}{\partial s_1}(0,0)=Y_{\mathbf j}$. Put 
$\varphi(s_1,s_2,t):=F_{X_{\mathbf i}}^t(c(s_1,s_2))$. Then 
$\frac{\partial\varphi}{\partial s_1}(0,0,t)=(F_{X_{\mathbf i}}^t)_*(Y_{\mathbf j})$ and 
thus $\frac{D}{\partial t}\frac{\partial\varphi}{\partial s_1}(0,0,0)=\Gamma_{X_{\mathbf i}}Y_{\mathbf j}$.
Hence
\[ \left\langle \Gamma_{X_{\mathbf i}}Y_{\mathbf j},\frac{\partial c}{\partial s_2}(0,0)\right\rangle
=\left\langle \frac{D}{\partial s_1}\frac{\partial\varphi}{\partial t}(0,0,0),
  \frac{\partial c}{\partial s_2}(0,0)\right\rangle =
-\left\langle \frac{\partial\varphi}{\partial t}(0,0,0),
 \frac{D}{\partial s_1}\frac{\partial c}{\partial s_2}(0,0)\right\rangle =0, \]
where we have used that 
$\frac{\partial\varphi}{\partial t}$ is orthogonal to $E_{\mathbf j}$ and $S_{\mathbf j}(x)$ 
is totally geodesic in $M$. As $\frac{\partial c}{\partial s_2}(0,0)$
can be chosen arbitrarily in $E_{\mathbf j}(x)$, the result follows. \EPf

\medskip

The next result shows that $\Gamma$ determines $\nabla\alpha$.

	\begin{prop}\label{gamma-nabla-alpha}
\begin{enumerate}
\item[(i)]  $\nabla_{X_{\mathbf i}}\alpha(Y_{\mathbf j},Z_{\mathbf k}) = 
\langle\Gamma_{X_{\mathbf i}}Y_{\mathbf j},Z_{\mathbf k}\rangle(v_{\mathbf j}-v_{\mathbf k})$;
\item[(ii)] $\Gamma_{X_{\mathbf i}}Y_{\mathbf j}=\nabla_{X_{\mathbf i}}\tilde Y_{\mathbf j}\ \mbox{mod $E_{\mathbf j}$}$;
\end{enumerate}
for all ${\mathbf i}\in \mathbf I^*$, $\mathbf j$, $\mathbf k\in\mathbf I$, $X_{\mathbf i}\in E_{\mathbf i}$, $Y_{\mathbf j}\in E_{\mathbf j}$, $Z_{\mathbf k}\in E_{\mathbf k}$.
Here $\tilde Y_{\mathbf j}$ is any smooth local extension 
of $Y_{\mathbf j}$ to a section of $E_{\mathbf j}$. 
\end{prop}

\Pf (i) follows from $\alpha((F^t_{X_{\mathbf i}})_*Y_{\mathbf j},(F^t_{X_{\mathbf i}})_*Z_{\mathbf k})
=(F^t_{X_{\mathbf i}})_*\alpha(Y_{\mathbf j},Z_{\mathbf k})$ by taking derivative with respect to the normal
connection at $t=0$ and 
using e.g.~$\alpha(X,Z_{\mathbf k})=\langle X,Z_{\mathbf k}\rangle v_{\mathbf k}$ for all
$X\in T_xM$ and the parallelism in the normal bundle of the right hand side. 

Part (ii) follows by comparing (i) with the formula (\ref{nabla-alpha}) for 
$\nabla\alpha$.
\EPf 

\begin{cor}
$\Gamma$ is also $\mathbb R$-linear in the lower 
argument, and thus can be extended as a bilinear map to 
$\sum_{\mathbf i\in\mathbf I^*}E_{\mathbf i} \times 
\sum_{\mathbf j\in\mathbf I} E_{\mathbf j}$, where $\sum$ means the algebraic sum.
\end{cor}

\Pf If $\mathbf i=\mathbf j$ this follows from the definition.
Otherwise we use Propositions~\ref{gamma-perp}
and~\ref{gamma-nabla-alpha}(ii). \EPf

\subsection{Consequences of the Codazzi equation}

Important consequences of Proposition~\ref{gamma-nabla-alpha}
are some formulae for permuting the arguments of $\Gamma$, which are 
obtained via the Codazzi equation. We introduce the 
following notation for $\mathbf i$, $\mathbf j$, 
$\mathbf k\in\mathbf I$ with $\mathbf i\neq\mathbf k$:
\[ \frac{v_{\mathbf j}-v_{\mathbf k}}{v_{\mathbf i}-v_{\mathbf k}}=\left\{\begin{array}{ll}
\lambda & \mbox{if $v_{\mathbf j}-v_{\mathbf k}=\lambda(v_{\mathbf i}-v_{\mathbf k})$}, \\
0 & \mbox{if $v_{\mathbf j}-v_{\mathbf k}$ is not a multiple of $v_{\mathbf i}-v_{\mathbf k}$.}
\end{array} \right. \]
We also denote by $(X)_{E_{\mathbf k}}$ the component of 
$X\in T_xM$ in $E_{\mathbf k}$.

\begin{prop}\label{gamma-permute}
\[ (\Gamma_{Y_{\mathbf j}}{X_{\mathbf i}})_{E_{\mathbf k}}=\frac{v_{\mathbf j}-v_{\mathbf k}}{v_{\mathbf i}-v_{\mathbf k}}
\,(\Gamma_{X_{\mathbf i}}{Y_{\mathbf j}})_{E_{\mathbf k}} \]
for all $\mathbf i$, ${\mathbf j}\in \mathbf I^*$, $\mathbf k\in\mathbf I$ 
with $\mathbf k\neq\mathbf i$ and 
$X_{\mathbf i}\in E_{\mathbf i}$, $Y_{\mathbf j}\in E_{\mathbf j}$. 
\end{prop}

\Pf Due to Proposition~\ref{gamma-nabla-alpha}(i),
\begin{eqnarray*}
\langle\Gamma_{X_{\mathbf i}}Y_{\mathbf j},Z_{\mathbf k}\rangle(v_{\mathbf j}-v_{\mathbf k}) & = & 
\nabla_{X_{\mathbf i}}\alpha(Y_{\mathbf j},Z_{\mathbf k})\\
&=&\nabla_{Y_{\mathbf j}}\alpha(X_{\mathbf i},Z_{\mathbf k})\\
&=&\langle\Gamma_{Y_{\mathbf j}}X_{\mathbf i},Z_{\mathbf k}\rangle(v_{\mathbf i}-v_{\mathbf k})
\end{eqnarray*}
which yields the desired result. \EPf  

\begin{lem}\label{gamma-colinear}
\begin{enumerate}
\item[(i)] If $(\Gamma_{X_{\mathbf i}}Y_{\mathbf j})_{E_{\mathbf k}}\neq0$ for ${\mathbf i}\in \mathbf I^*$, $\mathbf j$, $\mathbf k\in\mathbf I$,
then $v_{\mathbf i}$, $v_{\mathbf j}$, $v_{\mathbf k}$ are colinear.
\item[(ii)] For $\mathbf i$, $\mathbf j$, $\mathbf k\in \mathbf I^*$,
the condition $(\Gamma_{E_{\mathbf i}}E_{\mathbf j})_{E_{\mathbf k}}\neq0$ is symmetric
in $\mathbf i$, $\mathbf j$, $\mathbf k$. 
\end{enumerate}
\end{lem}

\Pf Part~(i) is clear from Codazzi. For part~(ii), note that 
the symmetry is obvious if $\mathbf i=\mathbf j=\mathbf k$ and a consequence 
of Codazzi if the indices are mutually different. 
In the remaining cases $(\Gamma_{E_{\mathbf i}}E_{\mathbf j})_{E_{\mathbf k}}=0$ 
by Codazzi and Proposition~\ref{gamma-perp}. \EPf

\begin{rem}
If $v_{\mathbf i}$, $v_{\mathbf j}$, $v_{\mathbf k}$ are colinear then the corresponding
focal hyperplanes $H_{\mathbf i}$, $H_{\mathbf j}$, $H_{\mathbf k}$
share a point in common. The converse holds if $M$ has rank~$2$.
In fact for $x\in M$,
$y=x+\xi\in H_{\mathbf i}(x)\cap H_{\mathbf j}(x)\cap H_{\mathbf k}(x)$
if and only if $\langle\xi,v_{\mathbf i}(x)\rangle=\langle\xi,v_{\mathbf j}(x)\rangle=
\langle\xi,v_{\mathbf k}(x)\rangle=1$ where $\xi\in\nu_xM$. 
\end{rem}

\section{Density of the image of $\Gamma$ and equivalence
between $\Gamma$ and $\nabla\alpha$}
\setcounter{equation}{0}

In the first theorem, we consider the algebraic span of the 
subsets $\Gamma_{E_{\mathbf i}(x)}E_{\mathbf j}$ for 
$\mathbf i$, $\mathbf j\in \mathbf I^*$, 
$\mathbf i\neq \mathbf j$.

\begin{thm}\label{p}
$\sum_{\ontop{\mathbf i,{\mathbf j}\in \mathbf I^*}{\mathbf i\neq \mathbf j}}\Gamma_{E_{\mathbf i}(x)}E_{\mathbf j}(x)$ is dense in $T_xM$ and moreover
contains $\sum_{{\mathbf i}\in \mathbf I^*}E_{\mathbf i}(x)$.   
\end{thm}

\Pf Set $\mathcal D(x)=\sum_{\ontop{\mathbf i,{\mathbf j}\in \mathbf I^*}{ \mathbf i
\neq \mathbf j}}\Gamma_{E_{\mathbf i}(x)}E_{\mathbf j}(x)$
for all $x\in M$. 
Then $\D$ is a possibly nonsmooth distribution of $M$ which 
however is invariant under all isometries of $V$ leaving $M$
invariant (Lemma~\ref{gamma-skew-and-invariant}(ii)). 

We first show that $E_{\mathbf k}\subset\mathcal D(x)$ for each $\mathbf 
k\in \mathbf I^*$.
Here we may assume that $M$ is finite dimensional, by passing to 
an irreducible finite dimensional slice of rank at least $2$ 
containing $E_{\mathbf k}$. Fix $x\in M$ and suppose $Y\in T_xM$ 
is orthogonal to $\mathcal D(x)$. 
We are going to show that for each $\mathbf l\in\mathbf I$,
the $E_{\mathbf l}$-component $Y_{\mathbf l}$ lies in $T_yM$ for all $y\in M$
and thus has to vanish, as otherwise the line 
in the direction of $Y_{\mathbf l}$ could be split off, 
contradicting the irreducibility of~$M$. 

So we fix $\mathbf l$ and, by Theorem~D in~\cite{HOT}, 
select indices $\mathbf i_1,\ldots,\mathbf i_r\in\mathbf I$ 
different from $\mathbf l$ such that any point $y\in M$ can be reached 
from $x$ by following piecewise smooth curves whose smooth arcs 
are contained in curvature spheres $S_{\mathbf i}$ with $\mathbf i\in
\{\mathbf i_1,\ldots,\mathbf i_r\}$. 
For such an~$\mathbf i$, consider $X\in E_{\mathbf i}(x)$. Then $Y(t):=(F_{X_{\mathbf i}}^t)_*Y$
is a smooth extension of $Y$ along the curve $\gamma(t):=F_{X_{\mathbf i}}^t(x)$ 
in $S_{\mathbf i}(x)$ (here we use the finite dimensionality of $M$),
which is everywhere orthogonal to $\mathcal D$. 
Next we split $Y(t)$ into its $E_{\mathbf i}$- and $E_{\mathbf i}^\perp$-components 
and note that the latter, which is given by $(F_X^t)_*(Y-Y_{\mathbf i})$,
is constant. In fact
\[ \frac{d}{dt}\Big|_{t=s}(F_X^t)_*(Y-Y_{\mathbf i})
=(F_X^s)_*(\Gamma_X(Y-Y_{\mathbf i})+\alpha(X,Y-Y_{\mathbf i})) = 0 \]
using $\alpha(E_{\mathbf i},E_{\mathbf i}^\perp)=0$, $\langle\Gamma_{E_{\mathbf i}}E_{\mathbf i}^\perp,E_{\mathbf i}\rangle
=\langle\Gamma_{E_{\mathbf i}}E_{\mathbf i},E_{\mathbf i}^\perp\rangle=0$, and 
$\langle\Gamma_{E_{\mathbf i}}Y,E_{\mathbf i}^\perp\rangle=\langle Y,\Gamma_{E_{\mathbf i}}E_{\mathbf i}^\perp\rangle
=0$ by the choice of $Y$. In particular the $E_{\mathbf l}$-component 
of $Y(t)$ is the constant vector $Y_{\mathbf l}$. Changing $X\in E_{\mathbf i}(x)$
we obtain in this way that, for any $y\in S_{\mathbf i}(x)$, 
the $E_{\mathbf l}$-component $Y_{\mathbf l}$ of $Y$ lies in $T_yM$ for 
all $y\in S_{\mathbf i}(x)$ and that moreover $Y_{\mathbf l}$ is the $E_{\mathbf l}$-component
of a vector $\tilde Y\in T_yM$ which is orthogonal to $\mathcal D(y)$. 
Repeating this argument with $\tilde Y$ in place of $Y$ and 
any $\mathbf i'\in\{\mathbf i_1,\ldots,\mathbf i_r\}$ in place of $\mathbf i$ eventually proves 
by induction the existence for any $y\in M$ of a vector in $T_yM$
whose $E_{\mathbf l}$-component is $Y_{\mathbf l}$. This finishes the proof
of the second statement. 

The proof of the first statement is similar but now simpler. 
$M$ is infinite dimensional and we suppose $Y\in T_xM$ is
orthogonal to $\mathcal D(x)$. Then $Y\in E_{\mathbf 0}(x)$ by the 
first part. A calculation similar to the above shows
that $Y(t):=(F_X^t)_*Y$ is constant for all $X\in E_{\mathbf i}(x)$, ${\mathbf i}\in \mathbf I^*$. 
As the set of points that can be reached from $x$ by following
piecewise smooth curves whose smooth arcs lie in finite
dimensional curvature spheres is dense in $M$, we deduce that the constant
vector $Y$ is everywhere tangent to $M$ and thus has to vanish so 
as not to contradict the irreducibility of $M$. \EPf

\begin{cor}\label{cor:p}
\begin{enumerate}
\item[(i)] For each $\mathbf k\in \mathbf I^*$, $E_{\mathbf k}=\sum
\{(\Gamma_{E_{\mathbf i}}E_{\mathbf j})_{E_{\mathbf k}}\;|\;v_{\mathbf i},v_{\mathbf j}\not\in\mathbb R v_{\mathbf k}\}$.
\item[(ii)]$\sum
\{(\Gamma_{E_{\mathbf i}}E_{\mathbf j})_{E_{\mathbf 0}}\;|\;\mbox{$\mathbf i$, ${\mathbf j}\in \mathbf I^*$ with $v_{\mathbf i}$, $v_{\mathbf j}$
lin.~dep.}\}$
is dense in $E_{\mathbf 0}$. 
\end{enumerate}
\end{cor}

\Pf (i) Since $\mathbf k\neq\mathbf 0$ we may assume $M$ is finite dimensional
by passing to an appropriate slice. By Theorem~\ref{p},
$E_{\mathbf k}=\sum_{\mathbf i\neq\mathbf j}(\Gamma_{E_{\mathbf i}}E_{\mathbf j})_{E_{\mathbf k}}$. Now for $\mathbf i\neq\mathbf j$ and
$\mathbf k\in\{\mathbf i,\mathbf j\}$ we have $\langle\Gamma_{E_{\mathbf i}}E_{\mathbf j},E_{\mathbf k}\rangle=0$. Hence 
$E_{\mathbf k}=\sum_{\ontop{\mathbf i\neq\mathbf j}{\mathbf i,\mathbf j\neq\mathbf k}}(\Gamma_{E_{\mathbf i}}E_{\mathbf j})_{E_{\mathbf k}}$ proving~(i).

(ii) The statement follows directly from Theorem~\ref{p}
as $(\Gamma_{E_{\mathbf i}}E_{\mathbf j})_{E_{\mathbf 0}}=0$ if $v_{\mathbf i}$, $v_{\mathbf j}$ are 
linearly independent by Codazzi. \EPf

\begin{thm}\label{thm:nabla-alpha-gamma}
If $\alpha_x$ is given, then $\Gamma_x$ and $(\nabla\alpha)_x$
determine one the other.  
\end{thm}

\Pf Suppose $\alpha_x$ and $\Gamma_x$ are given. 
Proposition~\ref{gamma-nabla-alpha}(i) can be rewritten as saying
\[ \nabla_{X_{\mathbf i}}\alpha(Y_{\mathbf j},Z_{\mathbf k})=-\alpha(\Gamma_{X_{\mathbf i}}Y_{\mathbf j},Z_{\mathbf k})
-\alpha(Y_{\mathbf j},\Gamma_{X_{\mathbf i}}Z_{\mathbf k}). \]
Since $\nabla\alpha$ is a tensor (and thus continuous
in all three arguments), $\nabla_X\alpha(Y,Z)$ is thus determined  
for $X$, $Y$, $Z\in T_xM$ with 
$X\perp E_{\mathbf 0}$, and hence by Codazzi for all $X$, $Y$, $Z\in T_x M$
as $\nabla_X\alpha(Y,Z)=0$ if two arguments lie in $E_{\mathbf 0}$. 

Conversely, we will show that $\Gamma_x$ is determined
by $\alpha_x$ and $\nabla\alpha_x$. So let $\tilde M$ be a 
second isoparametric submanifold in the same Hilbert space $V$.
We assume that $x\in M\cap\tilde M$, $T_xM=T_x\tilde M$, 
$\alpha_x=\tilde\alpha_x$
and $\nabla\alpha_x=\nabla\tilde\alpha_x$, and want to show that 
$\Gamma_x=\tilde\Gamma_x$. The assumption that the second 
fundamental forms coincide implies that $M$ and $\tilde M$ have the
same curvature distributions and the same curvature normals at $x$,
so more precisely we need to show that $\Gamma_x$ and 
$\tilde\Gamma_x$ coincide as mappings $E_{\mathbf i}(x)\times E_{\mathbf j}(x)\to T_xM$ 
for all ${\mathbf i}\in \mathbf I^*$, ${\mathbf j}\in\mathbf I$.  

We first note that 
$\Gamma_x=\tilde\Gamma_x$ on $E_{\mathbf i}\times E_{\mathbf j}$ 
if $\mathbf i\neq\mathbf j$. Indeed
$(\Gamma_{X_{\mathbf i}}Y_{\mathbf j})_{E_{\mathbf j}}=(\tilde\Gamma_{X_{\mathbf i}}Y_{\mathbf j})_{E_{\mathbf j}}=0$ by 
Proposition~\ref{gamma-perp}, 
and $(\Gamma_{X_{\mathbf i}}Y_{\mathbf j})_{E_{\mathbf k}}$ (resp.~$(\tilde\Gamma_{X_{\mathbf i}}Y_{\mathbf j})_{E_{\mathbf k}}$) 
for $\mathbf k\neq\mathbf j$ is determined 
by $\nabla\alpha_x=\nabla\tilde\alpha_x$ in view of
Proposition~\ref{gamma-nabla-alpha}(i),
where $X_{\mathbf i}\in  E_{\mathbf i}=E_{\mathbf i}(x)$, $Y_{\mathbf j}\in E_{\mathbf j}=E_{\mathbf j}(x)$.  

Moreover for $\mathbf i\neq\mathbf 0$ we have
$\Gamma_{E_{\mathbf i}}E_{\mathbf i}=\tilde\Gamma_{E_{\mathbf i}}E_{\mathbf i}=0$ 
if $E_{\mathbf i}$ is irreducible due to Proposition~\ref{gamma-ei-ei}.
So we only need to understand 
$\Gamma_x$, $\tilde\Gamma_x:E_{\mathbf i}(x)\times E_{\mathbf i}(x)\to E_{\mathbf i}(x)$ in case
$E_{\mathbf i}$ is reducible. Fix such an ${\mathbf i}\in \mathbf I^*$. 
Decompose $E_{\mathbf i}=E_{\mathbf i}'\oplus E_{\mathbf i}''$. 
Again by Proposition~\ref{gamma-ei-ei}
and using Lemma~\ref{gamma-skew-and-invariant}(i),
there are two cases that need to be discussed:
\begin{equation}\label{gamma-alpha}
 \Gamma,\ \tilde\Gamma: E_{\mathbf i}'\times E_{\mathbf i}''\to E_{\mathbf i}'
\quad\mbox{and}\quad
 \Gamma,\ \tilde\Gamma: E_{\mathbf i}''\times E_{\mathbf i}' \to E_{\mathbf i}'. 
\end{equation}
The discussion is similar in both cases so we restrict it
to the first one.

There exists an irreducible
finite dimensional rank $2$ slice $L$ of $M$ (resp.~$\tilde L$
of $\tilde M$) through $x$ 
containing $E_{\mathbf i}$, which of course is of nonreduced type
and hence of type $BC_2$.
Recall that $L$, $\tilde L$ are homogeneous.
Since $\alpha_x$, $(\nabla\alpha)_x$
and $\Gamma_x$ restrict to the corresponding objects for 
$L$, and similarly for $\tilde L$, 
it is enough to assume that $M=L$, $\tilde M= \tilde L$  
are finite dimensional, homogeneous of type $BC_2$, 
which we henceforth do.

Let $X_{\mathbf i}\in E_{\mathbf i}'$, $W_{\mathbf i}\in E_{\mathbf i}''$. Owing to Theorem~\ref{p},
we can write $W_{\mathbf i}$ as a finite sum
\begin{equation}\label{wi}
W_{\mathbf i}=\sum_{\ell}\Gamma_{Y_{\mathbf j_{\ell}}}Z_{\mathbf k_{\ell}}
\end{equation}
where for each $\ell$, $v_{\mathbf i}$, $v_{\mathbf j_{\ell}}$,
$v_{\mathbf k_{\ell}}$ are pairwise different curvature normals of $L$, $\tilde L$ (Note that 
the pair $(v_{\mathbf j_{\ell}},v_{\mathbf k_{\ell}})$ can repeat for 
different $\ell$.) 
Using formula~(\ref{gamma-homo}), we get
\begin{eqnarray*}
\Gamma_{X_{\mathbf i}}W_{\mathbf i} & = & \sum_{\ell}\Gamma_{X_{\mathbf i}}\Gamma_{Y_{\mathbf j_{\ell}}}Z_{\mathbf k_{\ell}} \\
& = & \sum_{\ell}\left(\Gamma_{Y_{\mathbf j_{\ell}}}
\Gamma_{X_{\mathbf i}}Z_{\mathbf k_{\ell}}
+ \Gamma_{\Gamma_{X_{\mathbf i}}Y_{\mathbf j_{\ell}}-\Gamma_{Y_{\mathbf j_{\ell}}}X_{\mathbf i}}Z_{\mathbf k_{\ell}}\right),
\end{eqnarray*}
and a similar formula for $\tilde\Gamma$. 
In order to finish the proof, we need only to check the claim
that on the right 
hand side of this formula $\Gamma$ is being computed always 
on a pair of vectors lying in different curvature distributions,
as we already know that $\Gamma$ and $\tilde\Gamma$ coincide
for such pairs of vectors.  

Recall that the positive root system of type $BC_2$ can be described as 
$\Lambda^+=\{\theta_1,\theta_2,2\theta_1,2\theta_2,\theta_1+\theta_2,
\theta_1-\theta_2\}$ so that we can write $\mathbf I=\{\theta_1,\theta_2,
\theta_1+\theta_2,\theta_1-\theta_2\}$ and assume that the 
index~$\mathbf i$ corresponds to~$\theta_1$. Now in the formula~(\ref{wi}) 
we may assume that $(\mathbf j_{\ell},\mathbf k_{\ell})=(\theta_1+\theta_2,\theta_1-\theta_2)$
or $(\theta_1-\theta_2,\theta_1+\theta_2)$ for all $\ell$, and the claim 
follows from remarking that $\Gamma_{X_{\mathbf i}}Z_{\mathbf k_{\ell}}$, $\Gamma_{X_{\mathbf i}}Y_{\mathbf j_{\ell}}$,
$\Gamma_{Y_{\mathbf j_{\ell}}}X_{\mathbf i}$ all lie in $E'_{\theta_2}$. \EPf

\section{Implications of the Gauss equation for $\Gamma$}
\setcounter{equation}{0}

The Gauss equation yields another sort of 
formulae for permuting arguments of $\Gamma$. 
Let $\mathbf i$, $\mathbf j$, $\mathbf k\in \mathbf I^*$, $\mathbf l\in\mathbf I$
and let $X_{\mathbf i}\in E_{\mathbf i}$, $Y_{\mathbf j}\in E_{\mathbf j}$, $Z_{\mathbf k}\in E_{\mathbf k}$, $W_{\mathbf l}\in E_{\mathbf l}$.
On one hand the Gauss equation yields
\begin{eqnarray*}
 \langle R(X_{\mathbf i},Y_{\mathbf j})Z_{\mathbf k},W_{\mathbf l}\rangle & = & 
\langle \alpha(X_{\mathbf i},W_{\mathbf l}),\alpha(Y_{\mathbf j},Z_{\mathbf k})\rangle
-\langle\alpha(X_{\mathbf i},Z_{\mathbf k}),\alpha(Y_{\mathbf j},W_{\mathbf l})\rangle \\
& = & 
-\langle X_{\mathbf i}\wedge Y_{\mathbf j},Z_{\mathbf k}\wedge W_{\mathbf l}\rangle\langle v_{\mathbf i},v_{\mathbf j}\rangle.
\end{eqnarray*}
On the other hand, extending these vectors to smooth local sections 
of the corresponding eigenbundles and using the same letters
for the extensions, we get

\begin{eqnarray} \nonumber
\langle R(X_{\mathbf i},Y_{\mathbf j})Z_{\mathbf k},W_{\mathbf l}\rangle 
& = & \langle \nabla_{X_{\mathbf i}}\nabla_{Y_{\mathbf j}}Z_{\mathbf k},W_{\mathbf l}\rangle
-\langle \nabla_{Y_{\mathbf j}}\nabla_{X_{\mathbf i}}Z_{\mathbf k},W_{\mathbf l}\rangle
-\langle\nabla_{[X_{\mathbf i},Y_{\mathbf j}]}Z_{\mathbf k},W_{\mathbf l}\rangle \\ \nonumber
& = & X_{\mathbf i}\langle\nabla_{Y_{\mathbf j}}Z_{\mathbf k},W_{\mathbf l}\rangle
-\langle\nabla_{Y_{\mathbf j}}Z_{\mathbf k},\nabla_{X_{\mathbf i}}W_{\mathbf l}\rangle \\ \nonumber
&&\qquad-Y_{\mathbf j}\langle\nabla_{X_{\mathbf i}}Z_{\mathbf k},W_{\mathbf l}\rangle
+\langle\nabla_{X_{\mathbf i}}Z_{\mathbf k},\nabla_{Y_{\mathbf j}}W_{\mathbf l}\rangle
-\langle\nabla_{[X_{\mathbf i},Y_{\mathbf j}]}Z_{\mathbf k},W_{\mathbf l}\rangle \\ \nonumber
& = & \langle\Gamma_{X_{\mathbf i}}Z_{\mathbf k},\Gamma_{Y_{\mathbf j}}W_{\mathbf l}\rangle 
-\langle\Gamma_{Y_{\mathbf j}}Z_{\mathbf k},\Gamma_{X_{\mathbf i}}W_{\mathbf l}\rangle 
-\langle\nabla_{[X_{\mathbf i},Y_{\mathbf j}]}Z_{\mathbf k},W_{\mathbf l}\rangle \\  \label{6}
&&\qquad +X_{\mathbf i}\langle\nabla_{Y_{\mathbf j}}Z_{\mathbf k},W_{\mathbf l}\rangle
-\langle\nabla'_{Y_{\mathbf j}}Z_{\mathbf k},\nabla_{X_{\mathbf i}}W_{\mathbf l}\rangle 
-\langle\Gamma_{Y_{\mathbf j}}Z_{\mathbf k},\nabla'_{X_{\mathbf i}}W_{\mathbf l}\rangle \\ \nonumber
&&\qquad  -Y_{\mathbf j}\langle\nabla_{X_{\mathbf i}}Z_{\mathbf k},W_{\mathbf l}\rangle
+\langle\nabla'_{X_{\mathbf i}}Z_{\mathbf k},\nabla_{Y_{\mathbf j}}W_{\mathbf l}\rangle 
+\langle\Gamma_{X_{\mathbf i}}Z_{\mathbf k},\nabla'_{Y_{\mathbf j}}W_{\mathbf l}\rangle,
\end{eqnarray}
where e.g.~$\nabla'_{Y_{\mathbf j}}Z_{\mathbf k}=\nabla_{Y_{\mathbf j}}Z_{\mathbf k}-\Gamma_{Y_{\mathbf j}}Z_{\mathbf k}$. 
Note that $\nabla'$ satisfies the properties of a metric connection
(Leibniz rule and compatibility with the metric) 
when both arguments are sections of curvature distributions,
and that $\nabla'_{Y_{\mathbf j}}Z_{\mathbf k}\in E_{\mathbf k}$ and 
$\nabla'_{Y_{\mathbf j}}Z_{\mathbf k}=(\nabla_{Y_{\mathbf j}}Z_{\mathbf k})_{E_{\mathbf k}}$ if $\mathbf j\neq\mathbf k$. Moreover, $\Gamma$ is parallel
with respect to $\nabla'$ in the following sense.

\begin{lem}\label{gamma-parallel}
$\nabla'_{X_{\mathbf i}}\left(\Gamma_{Y_{\mathbf j}}Z_{\mathbf k}\right)_{E_{\mathbf l}}=\left(\Gamma_{\nabla'_{X_\mathbf i}Y_{\mathbf j}}Z_{\mathbf k}\right)_{E_{\mathbf l}}+
\left(\Gamma_{Y_{\mathbf j}}\nabla'_{X_{\mathbf i}}Z_{\mathbf k}\right)_{E_{\mathbf l}}$. 
\end{lem}

\Pf To check this formula at $p\in M$, it suffices to differentiate at~$t=0$ 
the relevant vector fields along $\gamma(t)=F_{X_{\mathbf i}}^t(p)$. 
From the definition of $\Gamma$, we see that 
$(F^t_{X_{\mathbf i}})_*Y_{\mathbf j}$ is $\nabla'$-parallel along $\gamma$, as
well as $\left(\Gamma_{(F^t_{X_{\mathbf i}})_*Y_{\mathbf j}}
(F^t_{X_{\mathbf i}})_*Z_{\mathbf k}\right)_{E_{\mathbf l}}$, since this is 
equal to 
$\left((F_{X_{\mathbf i}^t})_*\Gamma_{\mathbf Y_j}Z_{\mathbf k}\right)_{E_{\mathbf l}}=(F_{X_{\mathbf i}}^t)_*\left(\Gamma_{\mathbf Y_j}Z_{\mathbf k}\right)_{E_{\mathbf l}}$ by Lemma~\ref{gamma-skew-and-invariant}(ii). 
This already shows that the result holds for $\nabla'$-parallel vector
fields $Y_{\mathbf j}$, $Z_{\mathbf k}$. Hence the result follows
in general,
since $E_{\mathbf j}$ and $E_{\mathbf k}$ are finite dimensional, and so
an arbitrary smooth vector field 
along $\gamma$ with values in one of those curvature distributions 
can be written as a finite linear combination of parallel vector fields
with smooth functions as coefficients. \EPf 

\medskip

To take care of the 
term with the bracket we introduce the following notation.
For any $A\in T_xM$ and $Z_{\mathbf k}$, $W_{\mathbf l}$ as above (with in particular
$\mathbf k\neq\mathbf l$), we put
\begin{equation}\label{gamma-first}
 \langle\Gamma_AZ_{\mathbf k},W_{\mathbf l}\rangle:=
\sum_{\mathbf m\in \mathbf I}\langle\Gamma_{Z_{\mathbf k}}W_{\mathbf l},
A_{\mathbf m}\rangle
\frac{v_{\mathbf l}-v_{\mathbf m}}{v_{\mathbf k}-v_{\mathbf l}}, 
\end{equation}
where $A_{\mathbf m}$ is the component of $A$ in $E_{\mathbf m}$. Note that
the sum converges as the coefficients 
$\frac{v_{\mathbf l}-v_{\mathbf m}}{v_{\mathbf k}-v_{\mathbf l}}$ are bounded,
and that this definition
indeed extends the domain of $\Gamma$ relative to the first
argument in view of
Proposition~\ref{gamma-permute}. The proof of that 
proposition also shows that 
\begin{equation}\label{gamma-nabla}
\langle\Gamma_AZ_{\mathbf k},W_{\mathbf l}\rangle=
\langle\nabla_AZ_{\mathbf k},W_{\mathbf l}\rangle.
\end{equation}

\begin{lem}\label{bracket}
\[ \langle\nabla_{[X_{\mathbf i},Y_{\mathbf j}]}Z_{\mathbf k},W_{\mathbf l}\rangle = \langle\Gamma_{\Gamma_{X_{\mathbf i}}Y_{\mathbf j}-\Gamma_{Y_{\mathbf j}}X_{\mathbf i}}Z_{\mathbf k},W_{\mathbf l}\rangle
+\langle\Gamma_{\nabla'_{X_{\mathbf i}}Y_{\mathbf j}}Z_{\mathbf k},W_{\mathbf l}\rangle-\langle\Gamma_{\nabla'_{X_{\mathbf j}}Y_{\mathbf i}}Z_{\mathbf k},W_{\mathbf l}\rangle. \]
\end{lem}

\Pf We have
$\langle\nabla_{[X_{\mathbf i},Y_{\mathbf j}]}Z_{\mathbf k},W_{\mathbf l}\rangle = \langle\nabla_{\Gamma_{X_{\mathbf i}}Y_{\mathbf j}-\Gamma_{Y_{\mathbf j}}X_{\mathbf i}}Z_{\mathbf k},W_{\mathbf l}\rangle+\langle\nabla_{\nabla'_{X_{\mathbf i}}Y_{\mathbf j}-\nabla'_{Y_{\mathbf j}}X_{\mathbf i}}Z_{\mathbf k},W_{\mathbf l}\rangle$, and the result follows from equation~(\ref{gamma-nabla}). \EPf

\medskip

Next we extend the domain of $\Gamma$ relative to its second argument by 
putting
\begin{equation}\label{gamma-second}
\langle \Gamma_{X_{\mathbf i}}Y,A_{\mathbf m}\rangle 
=-\langle Y,\Gamma_{X_{\mathbf i}}A_{\mathbf m}\rangle 
\end{equation}
for all $Y\in T_xM$, $\mathbf m\in\mathbf I$ and  
$A_{\mathbf m}\in E_{\mathbf m}$
to formulate the main result of this section.

\begin{thm}\label{gauss}
For any $\mathbf i$, $\mathbf j$, $\mathbf k\in \mathbf I^*$
and $\mathbf l\in \mathbf I$ with $\mathbf k\neq\mathbf l$, and
$X_{\mathbf i}\in E_{\mathbf i}$, $Y_{\mathbf j}\in E_{\mathbf j}$,
$Z_{\mathbf k}\in E_{\mathbf k}$, $W_{\mathbf l}\in E_{\mathbf l}$,
\[ \left\langle \left( [\Gamma_{X_{\mathbf i}},\Gamma_{Y_{\mathbf j}}]
-\Gamma_{\Gamma_{X_{\mathbf i}}Y_{\mathbf j}-\Gamma_{Y_{\mathbf j}}X_{\mathbf i}}\right) Z_{\mathbf k},W_{\mathbf l}\right\rangle 
=-\langle X_{\mathbf i}\wedge Y_{\mathbf j},Z_{\mathbf k}\wedge W_{\mathbf l}\rangle\langle v_{\mathbf i},v_{\mathbf j}\rangle. \]
\end{thm}

\Pf We combine equation~(\ref{6}) with Lemma~\ref{bracket} to write
\begin{eqnarray*}
\lefteqn{-\langle X_{\mathbf i}\wedge Y_{\mathbf j},Z_{\mathbf k}\wedge W_{\mathbf l}\rangle\langle v_{\mathbf i},v_{\mathbf j}\rangle  =  
 \left\langle \left( [\Gamma_{X_{\mathbf i}},\Gamma_{Y_{\mathbf j}}]
-\Gamma_{\Gamma_{X_{\mathbf i}}Y_{\mathbf j}-\Gamma_{Y_{\mathbf j}}X_{\mathbf i}}\right) Z_{\mathbf k},W_{\mathbf l}\right\rangle}  \\
&&\qquad+X_{\mathbf i}\langle\nabla_{Y_{\mathbf j}}Z_{\mathbf k},W_{\mathbf l}\rangle
+\langle\nabla'_{X_{\mathbf i}}Z_{\mathbf k},\nabla_{Y_{\mathbf j}}W_{\mathbf l}\rangle 
-\langle\Gamma_{\nabla'_{X_{\mathbf i}Y_{\mathbf j}}}Z_{\mathbf k},W_{\mathbf l}\rangle
-\langle\Gamma_{Y_{\mathbf j}}Z_{\mathbf k},\nabla'_{X_{\mathbf i}}W_{\mathbf l}\rangle \\ 
&&\qquad-Y_{\mathbf j}\langle\nabla_{X_{\mathbf i}}Z_{\mathbf k},W_{\mathbf l}\rangle
-\langle\nabla'_{Y_{\mathbf j}}Z_{\mathbf k},\nabla_{X_{\mathbf i}}W_{\mathbf l}\rangle 
+\langle\Gamma_{\nabla'_{Y_{\mathbf j}X_{\mathbf i}}}Z_{\mathbf k},W_{\mathbf l}\rangle
+\langle\Gamma_{X_{\mathbf i}}Z_{\mathbf k},\nabla'_{Y_{\mathbf j}}W_{\mathbf l}\rangle.
\end{eqnarray*}
The four last terms on the right hand side vanish (and similarly the 
four terms preceeding those) by Lemma~\ref{gamma-parallel} as
\[ \langle\nabla_{X_{\mathbf i}}Z_{\mathbf k},W_{\mathbf l}\rangle
=\langle\Gamma_{X_{\mathbf i}}Z_{\mathbf k},W_{\mathbf l}\rangle \]
(due to $\mathbf k\neq \mathbf l$) and
\[ \langle\nabla'_{Y_{\mathbf j}}Z_{\mathbf k},\nabla_{X_{\mathbf i}}W_{\mathbf l}\rangle = \langle\nabla'_{Y_{\mathbf j}}Z_{\mathbf k},\nabla'_{X_{\mathbf i}}W_{\mathbf l}\rangle 
+ \langle\nabla'_{Y_{\mathbf j}}Z_{\mathbf k},\Gamma_{X_{\mathbf i}}W_{\mathbf l}\rangle = 
-\langle \Gamma_{X_{\mathbf i}}\nabla'_{Y_{\mathbf j}}Z_{\mathbf k},
W_{\mathbf l}\rangle \]
(due to $\nabla'_{Y_{\mathbf j}}Z_{\mathbf k}\in E_{\mathbf k}$,
$\nabla'_{X_{\mathbf i}}W_{\mathbf l}\in E_{\mathbf l}$ 
and~$\mathbf k\neq\mathbf l$). \EPf

\medskip

We will use only special cases of this theorem. 
By using~ definitions~(\ref{gamma-first}), (\ref{gamma-second}) 
and Proposition~\ref{gamma-permute},
\[ (\Gamma_{X_{\mathbf i}}Y_{\mathbf j}-\Gamma_{Y_{\mathbf j}}X_{\mathbf i})_{E_{\mathbf m}}
= (\Gamma_{X_{\mathbf i}}Y_{\mathbf j})_{E_{\mathbf m}}\frac{v_{\mathbf i}-v_{\mathbf j}}
{v_{\mathbf i}-v_{\mathbf m}}, \]
for all $\mathbf m\neq\mathbf i$, and 
under the assumption that 
$(\Gamma_{X_{\mathbf i}}Y_{\mathbf j})_{E_{\mathbf i}}\perp
\Gamma_{Z_{\mathbf k}}W_{\mathbf l}$
(which holds for instance if $\mathbf i\neq \mathbf j$
or $v_{\mathbf i}$, $v_{\mathbf k}$, $v_{\mathbf l}$ are not colinear),
we can reformulate 
Theorem~\ref{gauss} more explicitly as
\begin{eqnarray}\label{eqn:gauss}\nonumber
\langle \Gamma_{X_{\mathbf i}}Z_{\mathbf k},\Gamma_{Y_{\mathbf j}}W_{\mathbf l}\rangle
-\langle\Gamma_{Y_{\mathbf j}}Z_{\mathbf k},\Gamma_{X_{\mathbf i}}W_{\mathbf l}\rangle
&=&\sum_{\mathbf m\in\mathbf I\setminus\{\mathbf i\}}
\langle \Gamma_{Z_{\mathbf k}}W_{\mathbf l},(\Gamma_{X_{\mathbf i}}Y_{\mathbf j})_{E_{\mathbf m}}
\rangle
\frac{v_{\mathbf l}-v_{\mathbf m}}{v_{\mathbf k}-v_{\mathbf l}}
\frac{v_{\mathbf i}-v_{\mathbf j}}{v_{\mathbf i}-v_{\mathbf m}} \\
&&\qquad  -\langle X_{\mathbf i}\wedge Y_{\mathbf j},Z_{\mathbf k}\wedge W_{\mathbf l}\rangle\langle v_{\mathbf i},v_{\mathbf j}\rangle.
\end{eqnarray}
The following two corollaries of Theorem~\ref{gauss}
are obtained from equation~(\ref{eqn:gauss}): the first one, by taking 
$\mathbf i=\mathbf l\neq\mathbf k=\mathbf j$; 
the second one, in case $\mathbf i\neq\mathbf j$,
by using that $(\Gamma_{X_{\mathbf i}}Y_{\mathbf j})_{E_{\mathbf m}}\neq0$ 
only if $v_{\mathbf i}$, $v_{\mathbf j}$, $v_{\mathbf m}$ are colinear.

\begin{cor}\label{gamma-cartan}
Let $\mathbf i$, ${\mathbf j}\in \mathbf I^*$, $\mathbf i\neq \mathbf j$. Then 
\[ \langle \Gamma_{X_{\mathbf i}}Y_{\mathbf j},\Gamma_{Z_{\mathbf j}}W_{\mathbf i}\rangle
+\langle \Gamma_{X_{\mathbf i}}Z_{\mathbf j},\Gamma_{Y_{\mathbf j}}W_{\mathbf i}\rangle
=\langle X_{\mathbf i},W_{\mathbf i}\rangle\langle Y_{\mathbf j},Z_{\mathbf j}\rangle\langle v_{\mathbf i},v_{\mathbf j}\rangle \]
for all $X_{\mathbf i}$, $W_{\mathbf i}\in E_{\mathbf i}$ and $Y_{\mathbf j}$, $Z_{\mathbf j}\in E_{\mathbf j}$. 
In particular
\[ \langle\Gamma_{X_{\mathbf i}}Y_{\mathbf j},\Gamma_{Y_{\mathbf j}}X_{\mathbf i}\rangle 
=\frac12\langle v_{\mathbf i},v_{\mathbf j}\rangle ||X_{\mathbf i}||^2||Y_{\mathbf j}||^2. \]
\end{cor}

\begin{cor}\label{stand-cpt}
Let $\mathbf i$, $\mathbf j$, $\mathbf k\in \mathbf I^*$, $\mathbf l\in\mathbf I$, and assume 
that $v_{\mathbf i}$, $v_{\mathbf k}$, $v_{\mathbf l}$ 
are not colinear. Then there 
exists a constant $c\in\R$ such that 
\[ \langle \Gamma_{X_{\mathbf i}}Z_{\mathbf k},\Gamma_{Y_{\mathbf j}}W_{\mathbf l}\rangle
=\langle \Gamma_{Y_{\mathbf j}}Z_{\mathbf k},\Gamma_{X_{\mathbf i}}W_{\mathbf l}\rangle
+c\cdot\langle \Gamma_{X_{\mathbf i}}Y_{\mathbf j},\Gamma_{Z_{\mathbf k}}W_{\mathbf l}\rangle \]
for all $X_{\mathbf i}\in E_{\mathbf i}$, $Y_{\mathbf j}\in E_{\mathbf j}$, $Z_{\mathbf k}\in E_{\mathbf k}$, $W_{\mathbf l}\in E_{\mathbf l}$. 
Moreover if $\mathbf i=\mathbf j$ or the lines spanned by $v_{\mathbf i}$, $v_{\mathbf j}$ and $v_{\mathbf k}$, $v_{\mathbf l}$
do not meet at a curvature normal, then $c=0$. On the other hand,
if those lines meet at $v_{\mathbf m}$ we have 
$c=\dfrac{v_{\mathbf l}-v_{\mathbf m}}{v_{\mathbf k}-v_{\mathbf l}}\cdot
\dfrac{v_{\mathbf i}-v_{\mathbf j}}{v_{\mathbf i}-v_{\mathbf m}}$.
\end{cor}

\begin{cor}\label{gamma-eiej-eiek-perp}
Let $\mathbf i$, $\mathbf j$, $\mathbf k\in \mathbf I^*$ with
$\mathbf k\neq\mathbf i$, $\mathbf j$ and 
$\frac{v_{\mathbf j}}{v_{\mathbf k}}\neq\frac{v_{\mathbf i}-v_{\mathbf j}}%
{v_{\mathbf k}-v_{\mathbf i}}$. Assume further that
$(\Gamma_{E_{\mathbf i}}E_{\mathbf j})_{E_{\mathbf m}}\perp \Gamma_{E_{\mathbf i}}E_{\mathbf k}$
for all $\mathbf m\in \mathbf I^*$ and $(\Gamma_{E_{\mathbf i}}E_{\mathbf j})_{E_{\mathbf k}}=0$
(both conditions hold if 
e.g.~$\Gamma_{E_{\mathbf i}}E_{\mathbf j}\subset E_{\mathbf 0}$).
Then 
\[ \Gamma_{E_{\mathbf i}}E_{\mathbf j}\perp\Gamma_{E_{\mathbf i}}E_{\mathbf k}. \] 
\end{cor}

\Pf We may assume $\mathbf i\neq\mathbf j$ since 
$\Gamma_{E_{\mathbf i}}E_{\mathbf i}=(\Gamma_{E_{\mathbf i}}E_{\mathbf i})_{E_{\mathbf i}}
\perp  \Gamma_{E_{\mathbf i}}E_{\mathbf k}$ by assumption. If $\mathbf 0$ does not
belong to one of the lines through $v_{\mathbf i}$, $v_{\mathbf j}$
and $v_{\mathbf i}$, $v_{\mathbf k}$, then 
$\Gamma_{E_{\mathbf i}}E_{\mathbf j}\perp E_{\mathbf 0}$
or $\Gamma_{E_{\mathbf i}}E_{\mathbf k}\perp E_{\mathbf 0}$ and the 
assertion follows again from the assumptions. Thus we are left with 
the case $v_{\mathbf i}$, $v_{\mathbf j}$, $v_{\mathbf k}$ are multiples
of each other.  

Let $X_{\mathbf i}$, $W_{\mathbf i}\in E_{\mathbf i}$, $Y_{\mathbf j}\in E_{\mathbf j}$
and $Z_{\mathbf k}\in E_{\mathbf k}$. Then equation~(\ref{eqn:gauss}) with
$\mathbf l=\mathbf i$ yields
\[ \langle \Gamma_{X_{\mathbf i}}Z_{\mathbf k}, \Gamma_{Y_{\mathbf j}}W_{\mathbf i}\rangle
=\langle \Gamma_{Z_{\mathbf k}}W_{\mathbf i}, \Gamma_{X_{\mathbf i}}Y_{\mathbf j}\rangle
\frac{v_{\mathbf i}-v_{\mathbf j}}{v_{\mathbf k}-v_{\mathbf i}}, \] 
and hence
\begin{equation}\label{horiz-perm}
 \langle \Gamma_{X_{\mathbf i}}Z_{\mathbf k}, \Gamma_{W_{\mathbf i}}Y_{\mathbf j}\rangle
\frac{v_{\mathbf j}}{v_{\mathbf i}}
=\langle \Gamma_{X_{\mathbf i}}Y_{\mathbf j}, \Gamma_{W_{\mathbf i}}Z_{\mathbf k}\rangle
\frac{v_{\mathbf k}}{v_{\mathbf i}}\frac{v_{\mathbf i}-v_{\mathbf j}}{v_{\mathbf k}-v_{\mathbf i}}. 
\end{equation} 
Here we 
have used that $(\Gamma_{E_{\mathbf i}}E_{\mathbf j})_{E_{\mathbf k}}=0$
implies $(\Gamma_{E_{\mathbf j}}E_{\mathbf k})_{E_{\mathbf i}}=0$, and
$(\Gamma_{E_{\mathbf i}}E_{\mathbf j})_{E_{\mathbf m}}\perp \Gamma_{E_{\mathbf i}}E_{\mathbf k}$ implies 
$(\Gamma_{E_{\mathbf i}}E_{\mathbf j})_{E_{\mathbf m}}\perp \Gamma_{E_{\mathbf k}}E_{\mathbf i}$
and
$(\Gamma_{E_{\mathbf j}}E_{\mathbf i})_{E_{\mathbf m}}\perp \Gamma_{E_{\mathbf i}}E_{\mathbf k}$.
Choosing $W_{\mathbf i}=X_{\mathbf i}$ in~(\ref{horiz-perm}) gives
$\langle \Gamma_{X_{\mathbf i}}Y_{\mathbf j}, \Gamma_{X_{\mathbf i}}Z_{\mathbf k}\rangle=0$,
and then by polarization, 
$\langle \Gamma_{X_{\mathbf i}}Y_{\mathbf j}, \Gamma_{X'_{\mathbf i}}Z_{\mathbf k}\rangle
+\langle \Gamma_{X'_{\mathbf i}}Y_{\mathbf j}, \Gamma_{X_{\mathbf i}}Z_{\mathbf k}\rangle=0$
for all $X'_{\mathbf i}\in E_{\mathbf i}$. Again by~(\ref{horiz-perm}), we have
$ \langle \Gamma_{X'_{\mathbf i}}Y_{\mathbf j}, \Gamma_{X_{\mathbf i}}Z_{\mathbf k}\rangle=
\langle \Gamma_{X_{\mathbf i}}Y_{\mathbf j}, \Gamma_{X'_{\mathbf i}}Z_{\mathbf k}\rangle
\frac{v_{\mathbf k}}{v_{\mathbf j}}\frac{v_{\mathbf i}-v_{\mathbf j}}{v_{\mathbf k}-v_{\mathbf i}}$, and owing to $\mathbf k\neq\mathbf j$, we get 
$\frac{v_{\mathbf k}}{v_{\mathbf j}}\frac{v_{\mathbf i}-v_{\mathbf j}}{v_{\mathbf k}-v_{\mathbf i}}\neq-1$, from which 
$\langle \Gamma_{X_{\mathbf i}}Y_{\mathbf j}, \Gamma_{X'_{\mathbf i}}Z_{\mathbf k}\rangle
=0$, as desired. \EPf

\medskip

From Proposition~\ref{gamma-permute} and Corollary~\ref{gamma-cartan} we get

\begin{cor}\label{norm1}
\[ \sum_{\mathbf k\in\mathbf I\setminus\{\mathbf i\}}||(\Gamma_{X_{\mathbf i}}Y_{\mathbf j})_{E_{\mathbf k}}||^2
\frac{v_{\mathbf k}-v_{\mathbf j}}{v_{\mathbf k}-v_{\mathbf i}} =\frac12\langle v_{\mathbf i},v_{\mathbf j}\rangle||X_{\mathbf i}||^2||Y_{\mathbf j}||^2 \]
for all $\mathbf i$, ${\mathbf j}\in \mathbf I^*$, $\mathbf i\neq\mathbf j$ and $X_{\mathbf i}\in E_{\mathbf i}$, $Y_{\mathbf j}\in E_{\mathbf j}$.  
\end{cor}

\section{Reduction to rank one slices}\label{sec:red}
\setcounter{equation}{0}

In this section we prove Proposition~\ref{gamma-continuous-rk1},
which reduces the proof of the continuity of $\Gamma_{X_{\mathbf i}}$ to
that on the infinite dimensional rank one slice 
containing $E_{\mathbf i}$. We also explain our strategy
to reduce the proof of
continuity of $\Gamma_{X_{\mathbf i}}$ along that 
slice to finding a constant $C$ such that 
$||\Gamma_{X_{\mathbf i}}Y_{\mathbf j}||\leq 
C\,||v_{\mathbf i}||\,||X_{\mathbf i}||\,||Y_{\mathbf j}||$ for all 
${\mathbf j}\in\mathbf I\setminus\{{\mathbf i}\}$ 
with $v_{\mathbf j}\in\R v_{\mathbf i}$ and $Y_{\mathbf j}\in E_{\mathbf j}$
(Lemma~\ref{finite-non-orthog}).

We first consider a special case which is simpler. 

\begin{prop}
Assume the affine 
Weyl group of $M$ to be of type $\tilde A - \tilde D - \tilde E$.
Let $\mathbf i$, ${\mathbf j}\in \mathbf I^*$ where we assume that
\[ v_{\mathbf i},\ v_{\mathbf j}\quad\mbox{are linearly independent,} \] 
and choose $X_{\mathbf i}\in E_{\mathbf i}$, $Y_{\mathbf j}\in E_{\mathbf j}$. If $v_{\mathbf i}\perp v_{\mathbf j}$ then
$\Gamma_{X_{\mathbf i}}Y_{\mathbf j}=0$; otherwise
\[ ||\Gamma_{X_{\mathbf i}}Y_{\mathbf j}||=\frac12||v_{\mathbf i}||\,||X_{\mathbf i}||\,||Y_{\mathbf j}||. \]
\end{prop}

\Pf Let $P$ be the affine span of $\{v_{\mathbf i}(x),v_{\mathbf j}(x)\}$ for 
some $x\in M$. Then $L_P(x)$ is a slice of
rank two and, by the assumption on $M$,
of type $A_1\times A_1$ or $A_2$.
In the first case $v_{\mathbf i}$ and $v_{\mathbf j}$ are orthogonal, $P$ does not contain
any other curvature normal and $\Gamma_{X_{\mathbf i}}Y_{\mathbf j}=0$. In the second case
there is precisely one more curvature normal $v_{\mathbf k}\in P$ for $\mathbf k\in \mathbf I^*\setminus
\{\mathbf i,\mathbf j\}$. For the sake of simplicity, by applying 
a translation of $V$ to $M$ we may assume
that the origin of $V$ lies at the intersection of 
$H_{\mathbf i}$, $H_{\mathbf j}$, $H_{\mathbf k}$.
Let $u_{\mathbf i}$, $u_{\mathbf j}$, $u_{\mathbf k}$ unit vectors orthogonal to 
$H_{\mathbf i}$, $H_{\mathbf j}$, $H_{\mathbf k}$, respectively. 
We may assume $u_{\mathbf k}=u_{\mathbf i}+u_{\mathbf j}$ by eventually multiplying some 
of these vectors by $-1$. As $c_{\mathbf i}(x)=x+\frac{v_{\mathbf i}}{||v_{\mathbf i}||^2}\in H_{\mathbf i}$,
we have $v_\mathbf r=-\frac{u_\mathbf r}{\langle x,u_\mathbf r\rangle}$
for $\mathbf r=\mathbf i$, $\mathbf j$, $\mathbf k$.

Thus Corollary~\ref{norm1} yields
\[ ||\Gamma_{X_{\mathbf i}}Y_{\mathbf j}||^2=||(\Gamma_{X_{\mathbf i}}Y_{\mathbf j})_{E_{\mathbf k}}||^2=
\frac12\langle v_{\mathbf i},v_{\mathbf j}\rangle \frac{v_{\mathbf k}-v_{\mathbf i}}{v_{\mathbf k}-v_{\mathbf j}}||X_{\mathbf i}||^2||Y_{\mathbf j}||^2
=\frac14||v_{\mathbf i}||^2||X_{\mathbf i}||^2||Y_{\mathbf j}||^2 \]
by a straightforward computation. \EPf

\medskip

If the affine Weyl group of $M$ is not of type
$\tilde A - \tilde D - \tilde E$, 
Corollary~\ref{norm1} does not suffice to estimate
$||\Gamma_{X_{\mathbf i}}Y_{\mathbf j}||$ for linearly independent $v_{\mathbf i}$, $v_{\mathbf j}$ as the 
sum on the left hand side in its statement contains more than one 
term and in general with different signs. 
We circumvent this difficulty by using
the classification of homogeneous compact rank two 
isoparametric submanifolds. 
 
\begin{prop}\label{gamma-cont-fin}
There exists a positive constant $C$ such that 
\[ ||\Gamma_{X_{\mathbf i}}Y||\leq C\, ||v_{\mathbf i}||\,||X_{\mathbf i}||\,||Y|| \]
whenever $X_{\mathbf i}$, $Y$ are
tangent to a finite dimensional slice through~$x$, 
$\mathbf i\in \mathbf I^*$ and $X_{\mathbf i} \in E_{\mathbf i}$.
\end{prop}

\Pf Assume first that such a slice $L$ is irreducible. 
Any curvature sphere of $M$ is contained in a finite dimensional
irreducible slice of rank at least two, so we may assume that the rank
of $L$ is at least two and apply Lemma~\ref{gamma-restriction}
to compute $\Gamma$ along $L$.
Since the rank of $L$ is bounded by the rank of $M$ and
multiplicities of $L$ are also multiplicities of $M$,
there are only finitely many possibilities 
for $L$ up to parallel translation (notice that scaling of the ambient metric
can be viewed as parallel translation in the radial direction).
Moreover, for given $L$, there exists the desired constant 
simply by finite dimensionality, so 
we have only to check that the same constant
works for the isoparametric submanifolds parallel to $L$.
However, this is a direct consequence of Lemma~\ref{parallel}. 

If~$L$ is a product of irreducible slices, then the result follows
from the case above together with the remark that 
$\Gamma_XY=0$ whenever $X$, $Y$ are tangent to different factors.
To prove it, we may assume that $X\in E_{\mathbf i}$,
$Y\in E_{\mathbf j}$ where $\mathbf i$, $\mathbf j\in\mathbf I^*$,
$\mathbf i\neq \mathbf j$.  Since the line through $v_{\mathbf i}$,
$v_{\mathbf j}$ contains no other curvature normals, 
we finish by noting that $\Gamma_{E_{\mathbf i}}E_{\mathbf j}\perp
E_{\mathbf i}$, $E_{\mathbf j}$. \EPf

\begin{prop}\label{gamma-continuous-rk1}
Let ${\mathbf i}\in \mathbf I^*$ and $X_{\mathbf i}\in E_{\mathbf i}=E_{\mathbf i}(x)$. Then 
$\Gamma_{X_{\mathbf i}}$ is continuous (that is, can be extended 
continuously to $T_xM$) if and only if $\Gamma_{X_{\mathbf i}}$ 
is continuous on 
$\sum_{\mathbf j}\{E_{\mathbf j}\;|\;\mbox{${\mathbf j}\in \mathbf I^*$, 
$v_{\mathbf j}\in\R v_{\mathbf i}$}\}$. 
\end{prop}
 
\Pf Write $T_xM=V_0\oplus V_1\oplus V_2$ where $V_0=E_{\mathbf 0}(x)$ 
and $V_1$, $V_2$ are the closures of
\[ \sum_{\mathbf j}\{E_{\mathbf j}\;|\;\mbox{${\mathbf j}\in \mathbf I^*$, 
$v_{\mathbf j}\not\in\R v_{\mathbf i}$}\}\quad\mbox{and}\quad
\sum_{\mathbf j}\{E_{\mathbf j}\;|\;\mbox{${\mathbf j}\in \mathbf I^*$, 
$v_{\mathbf j}\in\R v_{\mathbf i}$}\}, \]
respectively. By the closed graph theorem
$\Gamma_{X_{\mathbf i}}$ is continuous on $V_0$ (as $E_{\mathbf 0}(x)$ is a 
closed subspace lying in the domain of definition of $\Gamma_{X_{\mathbf i}}$
and $\Gamma_{X_{\mathbf i}}$ is skew-symmetric). We orthogonally
decompose $V_1$ into $\Gamma_{X_{\mathbf i}}$-invariant
subspaces  
\[ \sum\{E_{\mathbf j}\;|\;\mbox{${\mathbf j}\in \mathbf I^*
\setminus\{\mathbf i\}$, 
$v_{\mathbf j}\in\ell$}\}, \]
where $\ell$ runs over the lines in $\nu_xM\setminus\{0\}$
passing through $v_{\mathbf i}$. Of course for
only countably many lines these subspaces are nonzero.
Each such $\ell$ determines a finite dimensional rank two slice
through $x$ to which $E_i$ is tangent so, 
by Proposition~\ref{gamma-cont-fin},
$\Gamma_{X_{\mathbf i}}$ is continuous on~$V_1$ with 
$||\Gamma_{X_{\mathbf i}}||$ bounded by a constant independent of $\ell$.
The result follows.\hfill\mbox{} \EPf

\medskip

Identify $\mathbf I^*$ with $\mathcal A\times\Z$ as explained 
in section~\ref{prelim} and set $\mathbf i=(\alpha,i)$. 
Since $L_P(x)$,
$P=\R v_{\alpha,i}(x)$, is totally geodesic, we have 
$\Gamma_{X_{\alpha,i}}Y_{\alpha,j}\in E_{\mathbf 0}\oplus
\bigoplus_{k\in\mathbb Z}E_{\alpha,k}$. 
In these terms, Proposition~\ref{gamma-continuous-rk1}
reduces the proof of continuity of $\Gamma_{X_{\alpha,i}}$
to finding a constant $C$ such that 
$||\Gamma_{X_{\alpha,i}}Y||\leq C||Y||$ for all
$Y\in\sum_{j\in\mathbb Z}E_{\alpha,j}$; 
since $E_{\alpha,i}$ is finite dimensional, we 
may even take $Y\in\sum_{j\in\mathbb Z\setminus\{i\}}E_{\alpha,j}$.
The following lemma gives a sufficient condition 
to further simplify the proof in that the estimate 
only needs to be checked for all $Y\in E_{\alpha,j}$ (and all 
$j\in\mathbb Z\setminus\{i\}$).
Grosso modo it is required that for all $j\neq i$, the subspace 
$\Gamma_{E_{\alpha,i}}E_{\alpha,j}$ be orthogonal
to $\Gamma_{E_{\alpha,i}}E_{\alpha,k}$ for all 
but finitely many $k\in\mathbb Z$.

\begin{lem}\label{finite-non-orthog}
Let $W$ be a Hilbert space with an orthogonal 
decomposition $W=\oplus_{i\in\mathbb Z} W_i$, and let
$f:\sum_{i\in\mathbb Z} W_i\to W$ be a linear map. Assume 
there exists a constant $C>0$ such that $||fw_i||\leq C||w_i||$
for all $i\in\mathbb Z$ and $w_i\in W_i$, and that there 
exist injective maps $m_1,\ldots,m_r:\mathbb Z\to\mathbb Z$ such that 
$f(W_i)\perp f(W_j)$ unless $j\in\{m_1(i),\ldots,m_r(i)\}$. 
Then $||f||\leq \sqrt r C$ and thus $f$ can be continuously extended to $W$.
\end{lem}

\Pf Let $w=\sum_{i\in\mathbb Z}w_i$ where $w_i\in W_i$
and $w_i$ is nonzero for only finitely many indices.  
Then
\begin{eqnarray*}
||f(w)||^2 & = & \sum_{i,j\in\mathbb Z}\langle f(w_i),f(w_j)\rangle \\
&=& \sum_{i\in\mathbb Z}\sum_{k=1}^r\langle f(w_i),f(w_{m_k(i)})\rangle \\
&\leq&\sum_{i\in\mathbb Z}\sum_{k=1}^r||f(w_i)||||f(w_{m_k(i)})|| \\
&\leq&C^2\sum_{k=1}^r\sum_{i\in\mathbb Z}||w_i||||w_{m_k(i)}|| \\
&\leq&C^2\sum_{k=1}^r\left(\sum_{i\in\mathbb Z}||w_i||^2\sum_{i\in\mathbb Z}
||w_{m_k(i)}||^2\right)^{1/2} \\
&\leq&C^2\sum_{k=1}^r(||w||^2||w||^2)^{1/2}\\
&=&rC^2||w||^2,
\end{eqnarray*}
as we wished. \EPf

\medskip

Let $v\in P$ be a unit vector and define $\lambda_{\alpha,k}\in\R$ 
by $v_{\alpha,k}=\frac1{\lambda_{\alpha,k}}v$. Since 
$x+\lambda_{\alpha,k} v = x+\frac{v_{\alpha,k}}{||v_{\alpha,k}||^2}\in
H_{\alpha,k}$, we see that $\lambda_{\alpha,k}$ is the directed 
distance from $x$ to $H_{\alpha,k}$. Moreover we have
\[ \frac{v_{\alpha,k}-v_{\alpha,j}}{v_{\alpha,k}-v_{\alpha,i}}
=\frac{\lambda_{\alpha,k}-\lambda_{\alpha,j}}{\lambda_{\alpha,k}-\lambda_{\alpha,i}}\,\frac{\lambda_{\alpha,i}}{\lambda_{\alpha,j}}
=\frac{k-j}{k-i}\,\frac{\langle v_{\alpha,i},v_{\alpha,j}\rangle}
{||v_{\alpha,i}||^2} \]
as the hyperplanes $\{H_{\alpha,k}\}$ form an 
equidistant family for each fixed $\alpha$. Therefore
Corollary~\ref{norm1} yields:

\begin{prop}\label{gamma-rk1}
For each $\alpha\in\mathcal A$ and $i$, $j\in\Z$ with $i\neq j$ and 
for $X_{\alpha,i}\in E_{\alpha,i}$, 
$Y_{\alpha,j}\in E_{\alpha,j}$ we have
\[ \sum_{k\in\mathbb Z\setminus\{i\}} 
\frac{k-j}{k-i} 
||(\Gamma_{X_{\alpha,i}}Y_{\alpha,j})_{E_{\alpha,k}}||^2
+ ||(\Gamma_{X_{\alpha,i}}Y_{\alpha,j})_{E_{\mathbf 0}}||^2
=\frac12||v_{\alpha,i}||^2\,||X_{\alpha,i}||^2\,||Y_{\alpha,j}||^2. \]
\end{prop}

Unless the $E_{\alpha,k}$-components all vanish, 
a bound on $||\Gamma_{X_{\alpha,i}}Y_{\alpha,j}||$
is not immediately clear from this formula because the
sign of the factor $\frac{k-j}{k-i}$ changes with $k$, and its size  
could be arbitrarily small for certain values of $k$, $j$.
On the other hand, if $\Gamma_{E_{\alpha,i}}E_{\alpha,j}\subset E_{\mathbf 0}$
for all $j\in\mathbb Z\setminus\{i\}$, then 
Corollary~\ref{gamma-eiej-eiek-perp} says that 
$\Gamma_{E_{\alpha,i}}E_{\alpha,j}\perp\Gamma_{E_{\alpha,i}}E_{\alpha,k}$ 
unless $k=2i-j$. Therefore Lemma~\ref{finite-non-orthog}
can be applied and continuity of $\Gamma_{X_{\alpha,i}}$ already follows
from Proposition~\ref{gamma-rk1}. In general, the applicability
of both of these results relies on the vanishing of 
sufficiently many components of $\Gamma_{E_{\alpha,i}}E_{\alpha,j}$.
As it turns out that,
we shall see in section~\ref{image} that
$(\Gamma_{E_{\alpha,i}}E_{\alpha,j})_{E_{\alpha,k}}$ vanishes
for all $j$, $k\in\Z$, $j\neq i$, in most cases, and for sufficiently
many in the remaining ones.

\section{Root systems of isoparametric
submanifolds}\label{section:ars}
\setcounter{equation}{0}

As mentioned in section~\ref{prelim},
for a finite or infinite dimensional isoparametric submanifold,
Terng proved that the group generated by reflections in the focal hyperplanes 
in a fixed affine normal space is a finite or affine Weyl 
group, respectively~\cite{Te,Te4}.
Under our assumptions, such isoparametric submanifolds are always 
irreducible and 
homogeneous of rank at least two. It turns out then to be possible to   
refine the data given by the  
Weyl group into a root system associated to the 
focal hyperplanes. The refinement amounts to specializing the Coxeter graph of the Weyl 
group
to a Dynkin diagram by adjoining arrows to the double and triple links,
and adding concentric circles around certain vertices in the nonreduced case.
This root system is unique up to scaling,
and is a root system in the ordinary sense if $M$ is finite
dimensional and an affine root system otherwise. To distinguish 
one from the other, we also call the former a finite root system. 
The aim of this section is to describe this construction. 
We start by considering root systems attached to a set of hyperplanes, 
discuss affine root systems, including 
an outline of their classification and then associate a root system to the set of focal 
hyperplanes of an isoparametric submanifold. 

\smallskip

\textbf{Weyl groups and their Coxeter graphs.}
Let~$E$ be an affine Euclidean space and denote by~$T$ its 
group of translations (a finite dimensional real vector space).
Let $\mathcal H$ be a given set of affine hyperplanes in $E$ 
which is invariant under the group $W$ generated by all the orthogonal
reflections in the elements of~$\mathcal H$. It is assumed that 
the normal vectors to the $H\in\mathcal H$ span $T$ and that 
$W$ is a finite or an affine Weyl group.
In the first case $W$ has a fixed point that 
necessarily is contained in all $H\in\mathcal H$; taking this point  
as the origin in $E$, we can identify $E$ with $T$ and view the 
hyperplanes as linear subspaces. The second case may be characterized
as $\mathcal H$ consisting of finitely many families of equidistant
hyperplanes. It is well known~\cite[ch.~VI, \S~2, no.~5, Prop.~8]{bourbaki}
that these actually are described by a unique reduced root system 
$\Delta$ in $T$ as the set of hyperplanes
\begin{equation}\label{L-alpha-k}
 L_{\alpha,k}=\{x\in T\;|\;\alpha(x)=k\},\quad\alpha\in\Delta,\quad
k\in\mathbb Z 
\end{equation}
after choosing an appropriate point in $E$ (any \emph{special point},
namely, one through which there passes one hyperplane from each family
of parallel hyperplanes) to identify $E$ with $T$. We also recall that 
$W$ acts simply transitively on the set of Weyl chambers, which are 
the connected components of the complement of the union of the hyperplanes
in $\mathcal H$. Weyl chambers are also called alcoves in case $W$ is an 
affine Weyl group. The Coxeter graph of $W$ (or~$\mathcal H$) is obtained
from a Weyl chamber $\mathcal C$ by taking as vertices the walls
of $\mathcal C$ (hyperplanes bounding $\mathcal C$) and linking 
two vertices by $0$, $1$, $2$, $3$ or 
infinitely many edges according to whether the 
corresponding walls make an angle $\pi/2$, $\pi/3$, $\pi/4$, ~$\pi/6$
or are parallel 
(other cases cannot occur). $W$ is called irreducible if its 
Coxeter graph is connected. In this case $\mathcal C$ is a 
simplicial cone (resp.~simplex) if $W$ is finite (resp.~affine) 
and hence the Coxeter graph has~$n$ (resp.~$n+1$) vertices,
where $n=\dim E$ is called the rank of $W$. The isomorphism type
of the Coxeter graph is independent of the chosen Weyl chamber and determines
$\mathcal H$ and $W$ up to isomorphism. Here an isomorphism 
between two sets of hyperplanes $\mathcal H\subset E$, 
$\mathcal H'\subset E'$ as above with irreducible Weyl groups is a map $f:E\to E'$ that 
is the 
composition of an isometry with a homothety and takes $\mathcal H$ 
onto $\mathcal H'$. It turns out that the isomorphism classes
of irreducible finite (resp.~affine) Weyl groups correspond bijectively
to the irreducible reduced root systems $\Delta$ in $T$ and are 
correspondingly 
denoted by $A_n$ ($n\geq1$), $B_n$ ($n\geq2$), $C_n$ ($n\geq3$), 
$D_n$ ($n\geq4$), $E_n$ ($n=6$, $7$, $8$), $F_4$ and $G_2$ 
in case $W$ is finite, and 
$\tilde A_n$ ($n\geq1$), $\tilde B_n$ ($n\geq3$), 
$\tilde C_n$ ($n\geq2$), $\tilde D_n$ ($n\geq4$), $\tilde E_n$ 
($n=6$, $7$, $8$), $\tilde F_4$, $\tilde G_2$ in case $W$ is affine.        

It will be important to understand the orbits of the action
of~$W$ on $\mathcal H$ in case $W$ is affine. 
For a fixed Weyl chamber~$\mathcal C$, it
follows from the transitiveness of~$W$ on the set of Weyl chambers 
that each $H\in\mathcal H$ is conjugate under $W$ to some
wall of $\mathcal C$ and thus to (at least) one vertex of the 
Coxeter graph defined by~$\mathcal C$. The $W$-orbits in $\mathcal H$
are thus described by the next result, which is a simple consequence of
~\cite[ch.~IV, \S~1, no.~3, Prop.~3]{bourbaki}.

\begin{prop}
Two vertices of the Coxeter graph lie in the same $W$-orbit 
if and only if they belong to a connected subgraph containing 
only simple links. In particular, any two hyperplanes in 
$\mathcal H$ are conjugate if $W$ is of type $\tilde A_n$ ($n\geq2$),
$\tilde D_n$ ($n\geq4$), $\tilde E_n$ ($n=6$, $7$, $8$), and
$W$ acts with two (resp.~three) orbits in $\mathcal H$
if $W$ is of type $\tilde A_1$, $\tilde B_n$ ($n\geq3$),
$\tilde F_4$, $\tilde G_2$ (resp.~$\tilde C_n$ ($n\geq2$)).
\end{prop}
\smallskip
\textbf{Root systems associated to a set of hyperplanes.}
We now come to the main definition in this section. Let~$E$,
$T$, $\mathcal H$ and~$W$ be as above. 

\begin{defn}\label{ars}
A \emph{root system} associated to $\mathcal H$ is a subset~$R$
of $T\times\mathcal H$ such that 
\begin{enumerate}
\item[(i)] $v\neq0$ and $v\perp H$ for all $(v,H)\in R$.
\item[(ii)] $2{\inn{v}{v'}}/{||v||^2}\in\Z$ for all $(v,H)$,
$(v',H')\in R$.
\item[(iii)] The projection $T\times\mathcal H\to\mathcal H$ 
maps $R$ onto $\mathcal H$.
\item[(iv)] $R$ is invariant under $W$, that is, $(w_*v,wH)\in R$ 
for all $(v,H)\in R$ and $w\in W$. 
\end{enumerate}
\end{defn}

The \emph{rank} of $R$ is defined to be the rank of $W$. Each~$(v,H)\in R$
may be identified with the nonconstant affine mapping $E\to\mathbb R$ 
whose gradient is~$v$ and whose zero set is~$H$. In case $W$ is finite
and $E$ is identified with $T$ by taking the point 
$x\in\cap_{H\in\mathcal H}H$ as origin, these affine mappings are linear 
functionals and one gets a root system in the ordinary sense. 
In case $W$ is affine, one gets an affine root system
in the sense of Macdonald~\cite{Mac}. In this case
we thus call $R$ an \emph{affine root system} associated to~$\mathcal H$.
An equivalent definition, under a different name, has been given 
by Bruhat and Tits~\cite{B-T}.

It follows immediately from~(ii) that if $(v,H)$, $(v',H')\in R$ and 
$v'=\lambda v$ for some $\lambda\in\mathbb R$, then  
$\lambda\in\{\pm\frac12,\pm1,\pm2\}$. In particular, $R$ associates
to each $H\in\mathcal H$ either a pair $\{\pm v\}$ or a quadruple 
$\{\pm v,\pm 2v\}$ of nonzero normal vectors. 
Moreover it is clear that together with $R$ also 
$R_{\mathrm {red}}:=\{(v,H)\in S\;|\;(\frac12v,H)\not\in S\}$,
$R_{\mathrm {red}'}:=\{(v,H)\in S\;|\;(2v,H)\not\in S\}$, and
$\check R=\{(\check v,H)\;|\;(v,H)\in S\}$ are root 
systems associated to $\mathcal H$, where $\check v:=2v/||v||^2$. 
$R$ is called \emph{reduced} if it coincides with $R_{\mathrm{red}}$,
i.~e.~ if $R$ associates a pair $\{\pm v\}$ of normal vectors
to each $H\in\mathcal H$, and \emph{nonreduced} otherwise. 
$R$ is called \emph{irreducible}
if $W$ is irreducible.  

\smallskip
\textbf{Examples of affine root systems.} Assume $W$ is an 
irreducible affine Weyl group associated to a family of hyperplanes
$\mathcal H$ and let $\Delta$ denote the 
unique reduced root system satisfying~(\ref{L-alpha-k}). 
It follows from~\cite[ch.~VI, \S~2, Prop.~2]{bourbaki}
that 
\[ R=\{(h_\alpha,L_{\alpha,k})\;|\;\alpha\in\Delta,\ k\in\Z\} \]
is a reduced affine root system, where $h_\alpha\in T$ is defined
by $\langle h_\alpha,x\rangle=\alpha(x)$ for all $x\in T$. 
Since the distance between $L_{\alpha,k}$ and $L_{\alpha,k+1}$ 
is $1/||h_{\alpha}||$, $R$ can be equivalently 
described as 
\[ R=\{(v,H)\;|\;H\in\mathcal H,\ v\perp H,\ ||v||=1/d_H\}, \]
where $d_H$ denotes the minimal distance from $H$ to a 
parallel $H'\in\mathcal H\setminus\{H\}$. Thus to each $\mathcal H$ 
is associated a canonical reduced affine root system. 

\smallskip
\textbf{Restrictions of root systems.}
Let $E$, $T$, $\mathcal H$ and $W$ as above and suppose 
$R$ is a root system associated to~$\mathcal H$.
Let $E'$ be an affine subspace of $E$, $T'$ its group of translations
and $\mathcal H'=\{H\cap E'\;|\;H^\perp\subset T'\}$. 
Assume that the set of $v\in T'$ such that $(v,H)\in R$ for some
$H\in\mathcal H$ spans $T'$. Then we have the following
result whose proof is simple.

\begin{lem}
The set $R'$ of pairs $(v',H')\in T'\times\mathcal H'$ such that 
there exists $(v,H)\in R$ with $v=v'$ and $H\cap E'=H'$ is a 
root system in $E'$ associated to $\mathcal H'$.
\end{lem}

$R'$ is called the root system obtained from~$R$ by restriction 
to~$E'$.

\smallskip

\textbf{The classification of irreducible affine root systems.}
Let $R$ be an irreducible affine root system associated to~$\mathcal H$. 
For simplicity, we assume that the rank of $R$ is at least two
(the rank one case can also be done easily, but is not relevant to us). 
Suppose first that $R$ is \emph{reduced}. Then it is completely 
determined by its \emph{length function} $\ell:\mathcal H\to\R$, 
which is given by $\ell(H):=||v||$ for $(v,H)\in R$. Namely, 
$R=\{(v,H)\in T\times\mathcal H\;|\;v\perp H,\ ||v||=\ell(H)\}$.
Since $\ell$ is invariant under $W$ and $W$ acts transitively on the 
set of alcoves, $\ell$ is determined by its values on the walls of 
one fixed alcove, and thus on the vertices of the Coxeter graph. 
Moreover $\ell$ has to take the same value on any two vertices
that are linked by a single edge. More generally, it follows
from~(ii) in Definiton~\ref{ars} that $\ell(H)=c \,\ell(H')$ with 
$c\in \{2^{\pm\frac12}\}$ or $\{3^{\pm\frac12}\}$ if the vertices 
associated to $H$ and $H'$ are linked by $2$ or $3$ edges, 
respectively. In each case, we encode the actual choices
of $c\neq1$  in the Coxeter graph by adding to the  corresponding link
an arrow pointing to the vertex of shorter length. The so obtained
diagram is called the Dynkin diagram of $R$. It determines the root
system uniquely up to scaling once the correspondence between
vertices of the diagram and walls of the fixed alcove is given, and up 
to similarity without this extra piece of information. Here 
two irreducible root systems~$R$, $R'$ associated to families
$\mathcal H$, $\mathcal H'$ in $E$, $E'$, respectively, 
are called similar if there exists $\lambda>0$ and an affine mapping
$\varphi:E\to E'$ which is the composition of an isometry with 
a homothety such that $R'=\{(\lambda\varphi_*v,\varphi H)\;|\;
(v,H)\in R\}$. From the above restrictions, one obtains precisely the list 
in Table~1 for the possible Dynkin diagrams of irreducible reduced affine root 
systems of rank at least 2 . 

\setlength{\extrarowheight}{.15cm}
\[ \begin{array}{|c|c|}
\hline
Type & Diagram \\
\hline
\tilde A_n,\ n\geq2 & \Ant{}{}{}{}\\ & \\ 
\hline 
\vspace{-.35in} \mbox{} & \mbox{}\\
\tilde B_n,\ n\geq3 & \Bnta{}{}{}{}{} \\ 
\tilde B_n^{\mathsf v},\ n\geq3 & \Bntaii{}{}{}{}{} \\ & \\
\hline 
\tilde C_n,\ n\geq2 & \Cntaiii{}{}{}{}{}{}\\ 
\tilde C_n^{\mathsf v},\ n\geq2 & \Cntaii{}{}{}{}{}{}\\ 
\tilde C_n^{\prime},\ n\geq2 & \Cnta{}{}{}{}{}{}\\ 
\hline
\tilde D_n,\ n\geq4 & \Dnt{}{}{}{}{}{} \\
\hline
\tilde E_6 & \Evit{}{}{}{}{}{}{}\\
\hline
\tilde E_7 & \Eviit{}{}{}{}{}{}{}{}\\
\hline
\tilde E_8 & \Eviiit{}{}{}{}{}{}{}{}{}\\
\hline
\tilde F_4 & \Fivta{}{}{}{}{} \\
\tilde F_4^{\mathsf v} & \Fivtaii{}{}{}{}{} \\
\hline
\tilde G_2 & \Giitaii{}{}{}{}\\
\tilde G_2^{\mathsf v} & \Giita{}{}{}{}\\
\hline
\end{array} \]
\smallskip
\begin{center}
\textbf{Table~1:} Dynkin diagrams of reduced irreducible
affine root systems of rank at least~$2$.
\end{center}

\bigskip

That all diagrams in Table~1 occur indeed as diagrams of affine root systems
follows from the examples discussed above. The canonical affine root system 
associated to~$\mathcal H$ (more precisely, its similarity type) 
is denoted by the same symbol as the corresponding affine Weyl group. 
The geometric meaning of the arrows in this case can be explained as follows. 
Each vertex is associated to a family of equally spaced parallel 
hyperplanes, and an arrow always points to a family that is wider spaced. 
The difference between the diagrams of~$R$ and~$\check R$ 
is that the directions of the arrows are all reversed as the product
of their length functions is $2$. In this 
way, all diagrams listed above are obtained, except $\tilde C_n'$ 
which can be obtained by a modification of the construction
of the canonical root system~\cite{Mac}. 

Suppose now~$R$ is \emph{nonreduced}. In this case $R$ is completely determined
by $R_{\mathrm{red}}$ and the information for which $H\in\mathcal H$ 
there exists $v\in T$ with $(v,H)$, $(2v,H)\in\mathcal H$. This property
is invariant under the Weyl group and thus can be encoded in the Dynkin diagram
of $R_{\mathrm{red}}$ by adding a second, larger concentric circle 
around the corresponding vertex, following the notation of~\cite{loos2}. 
The diagram so obtained is called the Dynkin diagram of $R$ and 
determines $R$ as before up to scaling (similarity). Again~(ii) in
Definition~\ref{ars} restricts which vertices can admit a second
concentric circle. Namely, if $(v,H)$, $(2v,H)$, $(v',H')\in R$ 
then $m=\frac{\langle v,v'\rangle}{||v||^2}$,
$n=2\frac{\langle v,v'\rangle}{||v'||^2}\in\Z$ implying 
$mn\leq2$ and hence $||v'||^2=2||v||^2$ unless $v$, $v'$ are parallel or 
orthogonal. As we are assuming the rank of $R$ to be at least two,
a concentric circle can be added only to a vertex $v$ in the Dynkin
diagram of $R_{\mathrm{red}}$ which is only doubly linked to other 
vertices and for which the arrows point to $v$ . 
The possible diagrams are listed in Table~2, where
the type refers to $(R_\mathrm{red},R_\mathrm{red'})$.

\setlength{\extrarowheight}{.5cm}
\[ \begin{array}{|c|c|}
\hline
 Type & Diagram 
\\

\hline
\vspace{-1.3cm} \mbox{}&\mbox{} \\

(\tilde B_n,\tilde B_n^{\mathsf v}) & \Bntar{}{}{}{}{} \\

(\tilde C_n^{\mathsf v},\tilde C_n^\prime) & \Cntar{}{}{}{}{}{} \\

(\tilde C_n^\prime,\tilde C_n) &  \Cntarii{}{}{}{}{}{} \\

(\tilde C_n^{\mathsf v},\tilde C_n) & \Cntariii{}{}{}{}{}{}\\

(\tilde C_2,\tilde C_2^{\mathsf v}) & \Ciitariii{}{}{}\\
\hline
\end{array} \]
\smallskip
\begin{center}
\textbf{Table~2:} Dynkin diagrams of nonreduced irreducible
affine root systems of rank at least~$2$.
\end{center}

\bigskip

That all diagrams in Table~2 actually occur can be seen by the 
following construction. Enlarge $R_{\mathrm{red}}$ with the 
elements $(2v,H)$ for all $(v,H)\in R_{\mathrm{red}}$ such that 
$H$ is conjugate under $W$ to a vertex with a double circle. Then 
we only need to check condition~(ii) in Definition~\ref{ars}
to show that this enlarged set is a root system with the 
required properties. If there is only one vertex with an additional 
concentric circle this follows by the remark that for the two diagrams one obtains 
by deleting the additional circle and keeping or reversing the direction 
of the arrows to that vertex there exists always a root system. If there are 
two additional circles one applies the same argument to either of them
using the first step.

\smallskip
\textbf{Root systems of isoparametric submanifolds.}
Consider first the case of a finite dimensional homogeneous
compact isoparametric submanifold $M$ of rank at least two
in Euclidean space, which we may assume to be 
contained in a sphere around the origin.
By Dadok's theorem~\cite{D},
$M$ can be identified with a principal orbit of the 
isotropy representation of an irreducible
symmetric space $G/K$ of compact type, 
where $(G,K)$ is an effective symmetric pair
and $K$ is connected. Let $\Lg=\Lk+\Lp$
be the corresponding Cartan decomposition and fix $x\in M$. 
The affine and linear normal spaces of $M$ at $x$ coincide, and 
$\nu_xM$ is a maximal Abelian subalgebra $\La\subset\Lp$.
Recall that the root system $\Delta$ of $G/K$ with respect to~$\La$ 
is given by $\Delta=\{\alpha\in\La^*\setminus\{0\}\;|\;\Lp_\alpha\neq0\}$,
where $\Lp_\alpha=\{X\in\Lp\;|\;ad_H^2X=-\alpha(H)^2X\}$\ \mbox
for all $H\in\La\}$. One sees that the focal hyperplanes 
of $M$ in $\nu_xM=\La$ coincide with the kernels of the roots
in~$\Delta$. We identify each $\alpha\in\Delta$ with
the pair~$(h_\alpha,\ker\alpha)$, where $h_\alpha\in\La$
is defined by $\langle h_\alpha,v\rangle=\alpha(v)$ for all $v\in\La$,
in order to associate a root system to the set $\mathcal H$ of focal
hyperplanes of $M$ in $\nu_xM$.  We call this the root system
of~$M$.
It follows from the next result that it 
is independent of the identification of $M$ with a principal
orbit of the isotropy representation of a symmetric space. 
Hence it is well defined, up to scaling.

\begin{prop}
Let $G/K$, $G'/K'$ be irreducible symmetric spaces of compact type
of rank at least two, where $(G,K)$, $(G',K')$ are effective symmetric 
pairs and $K$, $K'$ are connected. Let $\Lg=\Lk+\Lp$, $\Lg'=\Lk'+\Lp'$
be the corresponding decompositions of the Lie algebras of $G$, $G'$ into the 
$\pm1$-eigenspaces of the involutions, respectively. Assume 
$\varphi:\Lp\to\Lp'$ is an isometry that maps a principal $K$-orbit 
onto a principal $K'$-orbit. Then, after multiplication by a suitable 
constant, $\varphi$ maps the root system of $G/K$ with respect to $\La$
onto the root system of $G'/K'$ with respect to~$\varphi(\La)$, where
$\La$ is any normal space to the principal $K$-orbit.
\end{prop}

\Pf By irreducibility, after multiplying $\varphi$ by a suitable constant, 
we may assume that $\langle\varphi x,\varphi y\rangle'=\langle x,y\rangle$
for all $x$, $y\in\Lp$, where $\langle,\rangle$, $\langle,\rangle'$
denote the negatives of the Killing forms of~$\Lg$, $\Lg'$,
respectively. Now it suffices to show that $\varphi$ extends to a
Lie algebra isomorphism from $\Lg=\Lk+\Lp$ to $\Lg'=\Lk'+\Lp'$ preserving
the decompositions. 

By effectiveness, we can identify $K$, $K'$ with subgroups of 
$O(\Lp)$, $O(\Lp')$, respectively.  Then, using the rank assumption,
$K$ (and similarly $K'$) is the maximal connected subgroup of 
$O(\Lp)$ with its orbits
(this follows from~\cite{Sim}, cf.~\cite[Prop.~4.3.9]{berndt-console-olmos},
or~\cite{EH2}). Since an isoparametric foliation is determined by 
any regular leaf, $\varphi$ has to map the $K$-orbit foliation 
to the $K'$-orbit foliation and thus conjugates $K$ to $K'$. 
Now $\varphi$ can be extended to a linear bijective map
$\Lg\to\Lg'$ preserving the decompositions by putting
$\varphi(A)=\varphi A\varphi^{-1}$ for all $A\in\Lk$. 
This is clearly a Lie algebra isomorphism from~$\Lk$ to~$\Lk'$, but 
also from $\Lg$ to $\Lg'$ as
\[ [\varphi(A),\varphi x] = \varphi A \varphi^{-1}(\varphi x)=\varphi[A,x], \]
and
\begin{eqnarray*}
\langle \varphi(A),[\varphi x,\varphi y]\rangle' & = &
\langle [\varphi(A),\varphi x],\varphi y\rangle' \\
&=&\langle \varphi[A,x],\varphi y\rangle' \\
&=&\langle [A,x],y\rangle \\
&=&\langle A,[x,y]\rangle \\
&=&\langle \varphi(A),\varphi[x,y]\rangle' 
\end{eqnarray*}
for all $A\in\Lk$ and $x$, $y\in\Lp$. 
Here we have used that $\varphi$ preserves the inner products as
$\Lk\perp\Lp$, $\Lk'\perp\Lp'$ and 
\begin{eqnarray*}
\langle \varphi(A),\varphi(B)\rangle' & = &
\langle \varphi(A),\varphi(B)\rangle_{\mathfrak k'}-\mathrm{tr}_{\mathfrak p'}
\varphi(A)\varphi(B)\\
& = & \langle A,B\rangle_{\mathfrak k}-\mathrm{tr}_{\mathfrak p}
AB\\
&=& \langle A,B\rangle
\end{eqnarray*}
for all $A$, $B\in\Lk$ where $\langle,\rangle_{\mathfrak k}$,
$\langle,\rangle_{\mathfrak k'}$ denote the negatives of the 
Killing forms of $\Lk$, $\Lk'$ and $\mathrm{tr}_{\mathfrak p}$,
$\mathrm{tr}_{\mathfrak p'}$ denote the traces of operators
on $\Lp$, $\Lp'$, respectively. \EPf
 
\medskip

By construction it is clear that the root system is invariant under isometries.
More precisely we have

\begin{lem}\label{inv-isom-rs}
Assume $M$ is finite dimensional.
Let $\varphi$ be an isometry from $V$ to another
Euclidean space $V'$ and let $M'=\varphi M$. Then, for any $x\in M$,
$\varphi$ maps the root system of $M$ in $x+\nu_xM$ to  
the root system of $M'$ in $\varphi(x)+\nu_{\varphi(x)}M'$, up 
to a scaling factor.
\end{lem}

Since a finite dimensional isoparametric submanifold $M$ is congruent
to an orbit of the isotropy representation of a symmetric space and the root
system of $M$ coincides with that of the symmetric space, we get
from the standard theory of symmetric spaces

\begin{prop}\label{restr-rs}
Assume $M$ is a finite dimensional homogeneous compact isoparametric
submanifold of rank at least two in an Euclidean space $V$ 
and let $L$ be an irreducible slice of $M$ of rank at least two 
through $x\in M$ with affine span $W$. Then the root 
system of $L$ associated to the focal hyperplanes in the affine normal
space of $L$ at $x$ in $W$ is obtained by restriction from the root
system of $M$ associated to the focal hyperplanes in $x + \nu_xM$,
up to scaling.    
\end{prop}

After this preparation, we come the main result of this section.

\begin{thm}
For each infinite dimensional connected complete full irreducible
isoparametric submanifold $M$ of rank at least two in a separable
Hilbert space~$V$ and $x\in M$,
there exists a naturally defined affine root system associated to the
family of focal hyperplanes in $x+\nu_xM$ which is unique up to scaling. 
\end{thm}

\Pf Let $E=x+\nu_xM$ and let $\mathcal H$ be the family of focal 
hyperplanes in $E$. Then the group~$W$ generated by the
reflections in the elements
of $\mathcal H$ is an affine Weyl group. 
Fix an alcove $\mathcal A$ and denote its walls by 
$H_{\mathbf i_1},\ldots H_{\mathbf i_{n+1}}$, which also parametrize
the vertices of the Coxeter graph of $W$. 
In order to have an affine root system 
associated to $\mathcal H$, up to scaling, we need to specify
in the Coxeter graph the arrows attached to the double and triple links,
and the possible additional circles around vertices that are only doubly 
linked to other vertices (with arrows pointing to the given one). 
We proceed as follows. For any two vertices $H_{\mathbf i_a}$, 
$H_{\mathbf i_b}$ that are doubly or triply linked, 
consider the finite dimensional
rank two slice $L_P(x)$ where $P$ is the affine span 
of the curvature normals $v_{\mathbf i_a}$, $v_{\mathbf i_b}$ and transfer
the information about arrows and additional circles from
the Dynkin diagram of $L_P(x)$ to the subdiagram of $M$
with vertices $H_{\mathbf i_a}$, $H_{\mathbf i_b}$. 
In this way we get well defined arrows. That also 
additional circles are well defined --- the only problem 
arises for the middle vertex in the $\tilde C_2$ graph
which lies in two such subgraphs --- follows from the fact that 
the information about additional circles can also be read 
off from the multiplicities (cf.~section~\ref{prelim}).

The construction of the affine root system of $M$ is also independent
of the chosen alcove. If $\mathcal A'$ is another alcove then 
there exists an element $w\in W$ with $\mathcal A'=wA$. 
The walls of $\mathcal A'$ are $ wH_{\mathbf i_1},\ldots,
 wH_{\mathbf i_{n+1}}$. Let $P$ be the affine span 
of $v_{\mathbf i_a}$, $v_{\mathbf i_b}$ for $a\neq b$. There 
exists an isometry $\varphi$ of $V$ that preserves $M$ and maps $x$
to $wx$. Then $\varphi(L_P(x))=L_P(wx)$ and $\varphi H_{\mathbf i_a}=
w H_{\mathbf i_a}$, $\varphi H_{\mathbf i_b}=
w H_{\mathbf i_b}$ as $\varphi$ coincides with $w$ on the affine normal
space $x+\nu_xM$ due to $\varphi v_{\mathbf i}(x)=v_{\mathbf i}(wx)
=wv_{\mathbf i}(x)$ for all $\mathbf i\in\mathbf I^*$. Therefore $\varphi$,
which maps the root system of $L_P(x)$ isometrically onto that 
of $L_P(wx)$, maps the Dynkin diagram of $L_P(x)$ with 
vertices $H_{\mathbf i_a}$, $H_{\mathbf i_b}$ isomorphically
to the Dynkin diagram of $L_P(wx)$ with vertices 
$w H_{\mathbf i_a}$, $w H_{\mathbf i_b}$. \EPf

\medskip

It is now clear that Lemma~\ref{inv-isom-rs} and Proposition~\ref{restr-rs}
carry over to to the infinite dimensional setting.

\begin{cor}
\begin{enumerate}
\item[(i)] The (finite or affine) root system of an 
irreducible slice of rank at least two of $M$ 
is obtained from that of $M$ by restriction. 
\item[(ii)] Let $\varphi$ be an isometry from $V$ to another
Hilbert space $V'$ and let $M'=\varphi M$. Then, for any $x\in M$,
$\varphi$ maps the root system of $M$ in $x+\nu_xM$ to  
the root system of $M'$ in $\varphi(x)+\nu_{\varphi(x)}M'$, up 
to a scaling factor.
\end{enumerate}
\end{cor}

\section{$\Gamma$ along rank one slices}\label{image}
\setcounter{equation}{0}

In order to use Proposition~\ref{gamma-rk1} effectively,
we have to understand which components
$(\Gamma_{E_{\mathbf i}}E_{\mathbf j})_{E_k}$ 
of $\Gamma_{E_{\mathbf i}}E_{\mathbf j}$ might be non-zero if
$v_{\mathbf i}$ and $v_{\mathbf j}$ are linearly dependent. 

The first result is a reduction to the rank two case.
 
\begin{lem}\label{rk2}
\begin{enumerate}
\item[(i)] If the affine Weyl group $W$ of $M$ is isomorphic 
to $\tilde A_n$ ($n\geq2$), $\tilde D_n$ ($n\geq4$), 
$\tilde E_n$ ($n=6$, $7$ or $8$) or $\tilde F_4$ then any 
rank one slice is contained in a 
slice $L$ of type $\tilde A_2$. 
\item[(ii)] If $W$ is isomorphic 
to $\tilde B_n$ or $\tilde C_n$ ($n\geq2$) then
any rank one slice is contained in a 
slice of type $\tilde A_2$ or
in a slice whose Dynkin diagram has a symbol which is obtained from
that of the Dynkin diagram of $M$ by replacing $n$ by $2$ and (if $W$ is 
isomorphic to $\tilde B_n$) $B$ by $C$ (e.g. if $M$ is of type
($\tilde B_n,\tilde B_n^{\mathsf v}$) then any rank one slice is contained in a slice of type $\tilde A_2$ or $(\tilde C_2,\tilde C_2^{\mathsf v})$).

\end{enumerate}
\end{lem}

\Pf We may assume the rank one slice to be infinite 
dimensional (since otherwise it is simply a curvature sphere
which is contained in an infinite dimensional rank one slice) and thus of the form $L_P$ with $P=\R v_{\mathbf i}$ for 
some ${\mathbf i}\in \mathbf I^*$. 
The focal hyperplane 
$H_{\mathbf i}$ bounds an alcove and thus corresponds to 
a vertex in the Dynkin diagram. If there exists a vertex  
which is joined by a single link
to the vertex corresponding to $H_{\mathbf i}$, then there exists a finite 
dimensional slice $L_Q$ of type $A_2$ containing $S_{\mathbf i}(x)$ 
such that $Q$ is the affine span of $\{v_{\mathbf i},v_{\mathbf j}\}$ for some 
${\mathbf j}\in \mathbf I^*$. If now $\tilde Q$ is the linear span of  $\{v_{\mathbf i},v_{\mathbf j}\}$
then we see that $L_{\tilde Q}$ is of type $\tilde A_2$ and contains
$L_P$. This proves (i) and shows in case (ii) that any rank one slice $L_P,
P=\R v_{\mathbf i}$, is contained in a slice of type $\tilde A_2$ if the 
hyperplane $H_{\mathbf i}$ is not conjugate under the Weyl group to an  
extremal vertex of the Dynkin diagram that is connected by a double link to  
another vertex. 

To study the remaining cases it is convenient to assume, 
after possibly translating $M$ in $V$, that the origin of $V$ lies 
in $\nu_xM$ and in fact 
in one focal hyperplane from each family of parallel
focal hyperplanes. In particular, the affine normal space is
identified with the normal space and the origin is a special
point in the sense of Bourbaki~\cite[ch.~V, \S~3, nr.~10]{bourbaki}.

Consider first the case $W$ isomorphic to $\tilde B_n$. 
Up to rescaling the metric in $\nu_xM$,
we can choose an orthonormal basis $\theta_1,\ldots,\theta_n$
of the dual space $(\nu_xM)^*$ such that
$H_{\alpha,i}=\{\xi\in\nu_xM\;|\;\alpha(\xi)=i\}$  
are the focal hyperplanes, where 
$\alpha\in\Delta^+ :=\{\theta_a,\ \theta_a\pm\theta_b\;|\; 1\leq a<b\leq n\}$ 
and $i\in \mathbb Z$.
In fact, $\Delta := \Delta^+\cup (-\Delta^+)$ is a a root system of type 
$B_n$ with $\Delta^+$ the positive roots 
(cf.~\cite[ch.~V, planche II]{bourbaki}). Since $\theta_1-\theta_2,\dots,\theta_{n-1}-\theta_n$, $\theta_n$ is a basis of $\Delta$ and 
$\theta_1 +\theta_2$ a highest root, the hyperplanes 
\[ \begin{array}{cc} 
 \theta_1-\theta_2=0,\ldots,\theta_{n-1}-\theta_n=0,\\ \\
\theta_n=0,\quad \theta_1+\theta_2=1 
\end{array}\]
are the walls of an alcove whose
associated Coxeter graph has the form \vspace{-.2in}\mbox{}
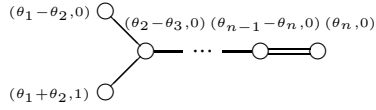
\[ \Bnt{\hspace{-.35in}{\scriptscriptstyle(\theta_1-\theta_2,0)}}{\hspace{-.35in}\scriptscriptstyle(\theta_1+\theta_2,1)}
{\hspace{-.1in}\scriptscriptstyle(\theta_2-\theta_3,0)}{\hspace{-.25in}\scriptscriptstyle(\theta_{n-1}-\theta_n,0)}{\hspace{.05in}\scriptscriptstyle(\theta_n,0)}\bigskip  \]

By the above discussion we may assume that $v_{\mathbf i}=v_{(\theta_n,0)}$
and thus can take $P\subset\nu_xM$ to be the
linear span of $\{v_{(\theta_{n-1}-\theta_n,0)},
v_{(\theta_n,0)}\}$. Then the Coxeter graph of $L_P(x)$ is 
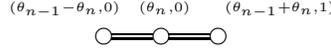
\[ \Ciit{\scriptscriptstyle\hspace{-.48in}(\theta_{n-1}-\theta_n,0)
}{\scriptscriptstyle\hspace{-.1in}(\theta_n,0)}
{\hspace{.05in}\scriptscriptstyle(\theta_{n-1}+\theta_n,1)} \]
Denote the reflection in $H_{\alpha,i}$ by $s_{\alpha,i}$.
Since $n\geq3$, there is an element in $W$ that maps the pair 
$(H_{\theta_{n-1}+\theta_n,1},H_{\theta_n,0})$ to 
$(H_{\theta_{n-1}-\theta_n,0},H_{\theta_n,0})$, namely
the composition $s_{\theta_{n-2}-\theta_{n-1},1}s_{\theta_{n-2}-\theta_{n-1},0}s_{\theta_n,0}$.
Hence both arrows of the Dynkin diagram of the slice point inward or outward in accordance with the direction of the arrow in the Dynkin diagram of $M$.
The desired result follows. 

Next, consider the case $W$ isomorphic to $\tilde C_n$. Up to rescaling the metric in $\nu_xM$,
we can choose an orthonormal basis $\theta_1,\ldots,\theta_n$
of $(\nu_xM)^*$ such that the hyperplanes 
\[ \begin{array}{cc} 
 \theta_1-\theta_2=0,\ldots,\theta_{n-1}-\theta_n=0,\\ \\
2\theta_n=0,\quad 2\theta_1=1 
\end{array}\]
are the walls of an alcove (cf.~\cite[ch.~V, planche III]{bourbaki}).
Now the associated Coxeter graph has the form 
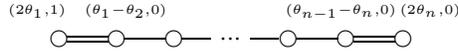
\[ \Cnt{\hspace{-.25in}{\scriptscriptstyle(2\theta_1,1)}}{\hspace{-.15in}\scriptscriptstyle(\theta_1-\theta_2,0)}{}%
{}{\hspace{-.3in}\scriptscriptstyle(\theta_{n-1}-\theta_n,0)}
{\hspace{0in}\scriptscriptstyle(2\theta_n,0)} \]  
We may assume that 
$v_{\mathbf i}=v_{(2\theta_n,0)}$ and can take $P\subset\nu_xM$ to be the
linear span of $\{v_{(\theta_{n-1}-\theta_n,0)},
v_{(2\theta_n,0)}\}$. 
Then the Coxeter graph of $L_P(x)$ is
 
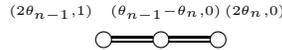
\[ \Ciit{\scriptscriptstyle\hspace{-.48in}(2\theta_{n-1},1)
}{\scriptscriptstyle\hspace{-.25in}(\theta_{n-1}-\theta_n,0)}
{\hspace{.05in}\scriptscriptstyle(2\theta_n,0)} \]
Since the finite Weyl group of type $C_n$ 
contains the full permutation group on the 
$\theta_{\mathbf i}$ (and their flips of signs), 
the pair $(H_{2\theta_{n-1},1},H_{\theta_{n-1}-\theta_n,0})$
is $W$-conjugate to 
$(H_{2\theta_1,1},H_{\theta_1-\theta_2,0})$.
This implies the stated result. \EPf

\medskip

To deal with the case of $\tilde G_2$, we need the following 
result.

\begin{lem}\label{g2}
If the affine Weyl group of $M$ is isomorphic 
to $\tilde G_2$ then $\Gamma_{E_{\mathbf i}}E_{\mathbf j}=0$ for all $\mathbf i$, ${\mathbf j}\in \mathbf I^*$ 
with $v_{\mathbf i}\perp v_{\mathbf j}$. 
\end{lem}

\Pf Let $P$ be the affine span of $\{v_{\mathbf i},v_{\mathbf j}\}$ for 
some $\mathbf i$, ${\mathbf j}\in \mathbf I^*$ with $v_{\mathbf i}\perp v_{\mathbf j}$ (fixing some point 
$x\in M$ as usually). Then $L_P$ is finite dimensional (since 
$0\not\in P$) and focalizes at the intersection point 
of $H_{\mathbf i}$ and $H_{\mathbf j}$ in the affine normal space. 
If there is no further focal line passing through this point 
then $L_P$ is of type $A_1\times A_1$ and the statement is clear.
Otherwise there must be exactly six focal lines passing through 
this point and $L_P$ is of type $G_2$. Therefore it suffices
to check the statement for finite dimensional homogeneous 
isoparametric submanifolds of type $G_2$, in which case 
it follows from simple properties of the associated root 
system. In fact, such a submanifold is congruent to a principal
orbit of the isotropy representation of a symmetric space 
$G/K$ and we can apply the discussion in subsection~\ref{fin-dim}.
We employ the notation from there and recall that we are 
dealing with a reduced root system. 
For $X_\lambda\in E_\lambda$, $Y_\mu\in E_\mu$ we have
\[ \Gamma_{X_\lambda}Y_\mu
=[\check X_\lambda,Y_\mu]\in[\Lk_\lambda,\Lp_\mu]\subset
\Lp_{\lambda+\mu}+\Lp_{\lambda-\mu}, \]
where $\check X_\lambda$ is the unique element
in $\Lk_\lambda$ satisfying $[\check X_\lambda,x]=X_\lambda$. 
In the specific case of $G_2$, it is a standard fact that 
orthogonal roots are always strongly orthogonal, i.e.~$\lambda\perp\mu$
implies that $\lambda\pm\mu$ are not roots, which can also
be immediately seen from the explicit description
of the root system: $\Lambda=\{\pm\lambda_1,\pm\lambda_2,\pm(\lambda_1+\lambda_2),
\pm(2\lambda_1+\lambda_2),\pm(3\lambda_1+\lambda_2),\pm(3\lambda_1+2\lambda_2)\}$.
The desired result follows. \EPf

\medskip

One of the main ingredients to describe the image of $\Gamma$ 
is surprisingly the following elementary lemma in Euclidean plane
geometry.

\begin{lem}\label{euclid}
Let $\ell_1$, $\ell_2$, $\ell_3$ be three different 
parallel lines in the plane and let $x_i$ be a point 
in $\ell_i$ for each $i$. Let $\ell_{ij}$ for 
$1\leq i<j\leq 3$ be the line through $x_i$ and $x_j$. 
\begin{enumerate}
\item[(i)] If the angles between any two of the six lines above
are multiples of $\pi/6$ but never $\pi/2$ then $x_1$, $x_2$, $x_3$ 
are colinear. 
\item[(ii)] If the angles between any two of the six lines above
are multiples of $\pi/4$ then either $x_1$, $x_2$, $x_3$ 
are colinear or one of the lines $\ell_1$, $\ell_2$, $\ell_3$, say $\ell_2$,
lies exactly in half way distance in between the other two. In the later
case $\ell_{13}$ is orthogonal to $\ell_2$.
\end{enumerate}
\end{lem}

\Pf Without loss of generality
we may assume that $\ell_2$ lies between $\ell_1$ and $\ell_3$.
Suppose that $x_1$, $x_2$, $x_3$ are not colinear, which is equivalent
to saying that $\ell_{12}$, $\ell_{13}$, $\ell_{23}$ are pairwise
different. Plainly from the formula for the angle sum,
it is readily seen that the angles of the triangle 
with vertices $x_1$, $x_2$, $x_3$ can only be $\pi/6$, $\pi/6$,
$2\pi/3$, or all $\pi/3$ in case~(i) and are necessarily 
$\pi/4$, $\pi/4$, $\pi/2$ in case~(ii). Therefore 
neither $\ell_{12}$ nor $\ell_{23}$ is orthogonal to $\ell_2$ 
and $\ell_2$ bisects the angle between $x_1-x_2$ and $x_3-x_2$. This 
implies that the angles of the triangle at $x_1$ and $x_3$ are equal
and that $\ell_{13}$ is orthogonal to $\ell_2$. This is a contradiction
in case~(i) and proves the result in case~(ii). \EPf

\medskip

Next we prove Theorem~C, as stated in the introduction.
It will be further sharpened by Theorem~\ref{image2}.

\begin{thm}\label{image1}
If the affine Weyl group $W$ of $M$ is isomorphic 
to $\tilde A_n$ ($n\geq2$), $\tilde D_n$ ($n\geq4$), 
$\tilde E_n$ ($n=6$, $7$ or $8$), $\tilde F_4$ or $\tilde G_2$
then
\[ \Gamma_{E_{\mathbf i}}E_{\mathbf j} \subset E_{\mathbf 0} \]
for all $\mathbf i$, ${\mathbf j}\in \mathbf I^*$ with $v_{\mathbf j}\in\R v_{\mathbf i}$. 
\end{thm}

\Pf By Lemma~\ref{rk2} we may assume $W$ is isomorphic to
$\tilde A_2$ or $\tilde G_2$ and by Proposition~\ref{gamma-ei-ei}
we may assume $\mathbf i\neq\mathbf j$. Assume that there 
exist $X_{\mathbf i}\in E_{\mathbf i}$, $Y_{\mathbf j}\in E_{\mathbf j}$ and 
$\mathbf m\in \mathbf I^*$ with 
$(\Gamma_{X_{\mathbf i}}Y_{\mathbf j})_{E_{\mathbf m}}\neq0$. Since $L_P$ for $P=\R v_{\mathbf i}(x)$ is totally 
geodesic, we must then have 
$v_{\mathbf m}\in \R v_{\mathbf i}$, and by Propositions~\ref{gamma-ei-ei} and~\ref{gamma-perp}, 
$\mathbf m\neq \mathbf i$, $\mathbf j$. According to Theorem~\ref{p}
there exist $\mathbf k$, $\mathbf l\in \mathbf I^*$ with $v_{\mathbf k}$, $v_{\mathbf l}$ 
linearly independent and $Z_{\mathbf k}\in E_{\mathbf k}$, $W_{\mathbf l}\in E_{\mathbf l}$ 
such that $\langle(\Gamma_{X_{\mathbf i}}Y_{\mathbf j})_{E_{\mathbf m}},\Gamma_{Z_{\mathbf k}}W_{\mathbf l}\rangle\neq0$.
Since the line through $v_{\mathbf k}$ and $v_{\mathbf l}$ meets $P$ in 
at most one point, we have 
$\langle\Gamma_{X_{\mathbf i}}Y_{\mathbf j},\Gamma_{Z_{\mathbf k}}W_{\mathbf l}\rangle=
\langle(\Gamma_{X_{\mathbf i}}Y_{\mathbf j})_{E_{\mathbf m}},\Gamma_{Z_{\mathbf k}}W_{\mathbf l}\rangle\neq0$
and no three among $v_{\mathbf i}$, $v_{\mathbf j}$, $v_{\mathbf k}$, $v_{\mathbf l}$ are colinear. 
Thus Corollary~\ref{stand-cpt} yields 
\[ 0\neq\langle(\Gamma_{X_{\mathbf i}}Y_{\mathbf j})_{E_{\mathbf m}},\Gamma_{Z_{\mathbf k}}W_{\mathbf l}\rangle
=\langle \Gamma_{Z_{\mathbf k}}Y_{\mathbf j},\Gamma_{X_{\mathbf i}}W_{\mathbf l}\rangle
+c\cdot\langle \Gamma_{X_{\mathbf i}}Z_{\mathbf k},\Gamma_{Y_{\mathbf j}}W_{\mathbf l}\rangle \]
for some $c\in \R$. At least one of the terms on the right hand side
is nonzero, so at least one of the following two cases must be true:
(i) the lines through 
$v_{\mathbf k}$, $v_{\mathbf j}$ and $v_{\mathbf i}$, $v_{\mathbf l}$ meet at some 
curvature normal, say $v_{\mathbf n}$;
(ii) the lines through 
$v_{\mathbf i}$, $v_{\mathbf k}$ and $v_{\mathbf j}$, $v_{\mathbf l}$ meet at some 
curvature normal. The analysis is completely 
similar in both cases, so we assume ~(i) is true. 
Now we have the following picture for the focal lines:
$H_{\mathbf i}$, $H_{\mathbf j}$, $H_{\mathbf m}$ are parallel
and pairwise different, and $H_{\mathbf k}$, $H_{\mathbf l}$, 
$H_{\mathbf n}$ are three other lines such that 
$H_{\mathbf k}$, $H_{\mathbf l}$ meet at $H_{\mathbf m}$,
$H_{\mathbf k}$, $H_{\mathbf n}$ meet at $H_{\mathbf j}$,
and $H_{\mathbf l}$, $H_{\mathbf n}$ meet at $H_{\mathbf i}$.
The intersection points are not colinear as $\mathbf k\neq\mathbf l$. 
However this is in contradiction with Lemma~\ref{euclid}
since the angle between any two of these lines is a 
multiple of $\pi/6$ and no two among them are orthogonal.
The non-orthogonality is of course automatic in the 
$\tilde A_2$ case and follows from Lemmata~\ref{gamma-colinear}(ii) 
and~\ref{g2} in the $\tilde G_2$ case. \EPf

\medskip

In view of the discussion at the end of section~\ref{sec:red},
Theorem~\ref{image1}
already yields continuity of $\Gamma$ in all cases 
but $\tilde B_n$ and $\tilde C_n$. 

\begin{cor}\label{adefg}
If the affine Weyl group of $M$ 
is of type $\tilde A$, $\tilde D$, $\tilde E$, $\tilde F$ 
or $\tilde G$ then $\Gamma_{X_{\mathbf i}}$ is continuous for all
$X_{\mathbf i}\in E_{\mathbf i}$ and~${\mathbf i}\in \mathbf I^*$.

\end{cor}

The next theorem extends Theorem~\ref{image1} to the remaining cases
of $\tilde B$ and $\tilde C$. 
However, to get continuity of $\Gamma_{X_{\mathbf i}}$ 
in those cases, it will be necessary
to refine its information. The required refinements are given by 
Propositions~\ref{gamma-image-1} (case $E_{\mathbf i}$ is irreducible) 
and~\ref{gamma-image-2} (case $E_{\mathbf i}$ is reducible).
 
As already in the last section, we identify the index set 
$\mathbf I^*$ with $\mathcal A\times\Z$.

\begin{thm}\label{gamma-image-general}
Let $\alpha\in\mathcal A$ and $i$, $j\in\Z$ with $i\neq j$. 
Then
\[ \Gamma_{E_{\alpha,i}}E_{\alpha,j}
\subset E_{\mathbf 0} \oplus E_{\alpha,2i-j} \oplus E_{\alpha,2j-i}
\oplus E_{\alpha,\frac{i+j}2}, \]
where the last term is to be omitted if $i+j$ is odd. 
\end{thm}

\Pf We may assume $W\cong\tilde C_2$ due to Theorem~\ref{image1}
and Lemma~\ref{rk2}(ii). Since the rank one slice $L_P$
for $P=\R v_{\alpha,i}(x)$ is totally geodesic,
$\Gamma_{E_{\alpha,i}}E_{\alpha,j}\subset E_{\mathbf 0}\oplus
\bigoplus_{k\in\mathbb Z}E_{\alpha,k}$.
Assume $(\Gamma_{X_{\alpha,i}}Y_{\alpha,j})_{E_{\alpha,m}}\neq0$ for some
$X_{\alpha,i}\in E_{\alpha,i}$, $Y_{\alpha,j}\in E_{\alpha,j}$ and some
$m\in \Z$, where $m\neq i$, $j$ necessarily. Then we find by 
Theorem~\ref{p}
$(\beta,k)$, $(\gamma,\ell)\in \mathcal A\times\Z$ with $\beta\neq\gamma$ 
and $Z_{\beta,k}\in E_{\beta,k}$, $W_{\gamma,\ell}\in E_{\gamma,\ell}$ 
such that $\langle(\Gamma_{X_{\alpha,i}}Y_{\alpha,j})_{E_{\alpha,m}},
\Gamma_{Z_{\beta,k}}W_{\gamma,\ell}\rangle\neq0$ and deduce 
as in the proof of Theorem~\ref{image1} that 
$(\Gamma_{E_{\alpha,i}}E_{\gamma,\ell})_{E_{\delta,n}}$ and
$(\Gamma_{E_{\alpha,j}}E_{\beta,k})_{E_{\delta,n}}$ are not zero
for some $(\delta,n)\in \mathcal A\times I$. Therefore we have again 
three parallel lines $H_{\alpha,i}$,
$H_{\alpha,j}$, $H_{\alpha,m}$ in the affine
normal space, and three further lines 
$H_{\beta,k}$,
$H_{\gamma,\ell}$, $H_{\delta,n}$
which are pairwise different and such that two of which  
intersect on $H_{\alpha,r}$ for each $r\in\{i,j,m\}$. Thus 
Lemma~\ref{euclid} implies that one of the parallel lines has
to lie exactly in the middle between the other two lines or,
equivalently, that one of the indices $i$, $j$, $m$ is the 
arithmetic mean of the other two, that is $m=\frac{i+j}2$,
$2i-j$ or $2j-i$.\hfill\mbox{ } \EPf

\medskip

Combining Theorem~\ref{gamma-image-general} 
with Corollary~\ref{gamma-eiej-eiek-perp} yields
the following result, which makes Lemma~\ref{finite-non-orthog}
applicable in general. 

\begin{cor}\label{gamma-eiej-eiek-perp-2}
Let $\alpha\in\mathcal A$ and $i$, $j$, $k\in\mathbb Z$ 
with $i\neq j$. Then
\[ \Gamma_{E_{\alpha,i}}E_{\alpha,j}\perp\Gamma_{E_{\alpha,i}}E_{\alpha,k} \]
if $k$ is not one of
\[ 4j-3i,\ 2j-i,\ j,\ \frac{i+j}2,\ \frac{3i+j}4,\ \frac{3i-j}2,\ 2i-j,\ 
3i-2j. \]
\end{cor}

\Pf We may assume that $k\neq i$. 
Let $\mathbf i=(\alpha,i)$, $\mathbf j=(\alpha,j)$, $\mathbf k=(\alpha,k)$. 
Now the condition
$\frac{v_{\mathbf j}}{v_{\mathbf k}}\neq\frac{v_{\mathbf i}-v_{\mathbf j}}%
{v_{\mathbf k}-v_{\mathbf i}}$ in Corollary~\ref{gamma-eiej-eiek-perp}
is equivalent to $k\neq 2i-j$. Moreover
$(\Gamma_{E_{\mathbf i}}E_{\mathbf j})_{E_{\mathbf m}}\perp \Gamma_{E_{\mathbf i}}E_{\mathbf k}$
for all 
$\mathbf m\in \mathbf I^*$ follows from Theorem~\ref{gamma-image-general} 
together with 
$\{2i-j,2j-i,\frac{i+j}2\}\cap\{2i-k,2k-i,\frac{i+k}2\}=\varnothing$. 
Finally $(\Gamma_{E_{\mathbf i}}E_{\mathbf j})_{E_{\mathbf k}}=0$ follows 
from Theorem~\ref{gamma-image-general}
and $k\neq 2i-j$, $2j-i$, $\frac{i+j}2$. \EPf 
 
\medskip

Let $P\subset x+\nu_xM$ be an affine line containing 
precisely four curvature normals, say $v_{\mathbf i}$, $v_{\mathbf j}$, $v_{\mathbf k}$, $v_{\mathbf l}$.
Then the corresponding slice is necessarily finite
dimensional and either of type $B_2$ or $(BC)_2$. 
We may assume $v_{\mathbf i}\perp v_{\mathbf j}$ and $v_{\mathbf k}\perp v_{\mathbf l}$ 
after an eventual permutation of the indices. The slice is
of type $B_2$ if $E_{\mathbf i}$, $E_{\mathbf j}$, $E_{\mathbf k}$, $E_{\mathbf l}$ are all irreducible,
and of type $(BC)_2$ if one pair among $E_{\mathbf i}$, $E_{\mathbf j}$ and $E_{\mathbf k}$, $E_{\mathbf l}$
is irreducible and the other is reducible.

\begin{lem}\label{b2bc2}
Let $P\subset x+\nu_xM$ be an affine line containing 
precisely four curvature normals, say $v_{\mathbf i}$, $v_{\mathbf j}$, $v_{\mathbf k}$, $v_{\mathbf l}$
with $v_{\mathbf i}\perp v_{\mathbf j}$ and $v_{\mathbf k}\perp v_{\mathbf l}$.
\begin{enumerate}
\item[(i)] $\Gamma_{E_{\mathbf i}}E_{\mathbf j}=0$ if the corresponding slice
is of type $B_2$ and $v_{\mathbf i}$, $v_{\mathbf j}$ correspond to the long roots in this slice.
\item[(ii)] $(\Gamma_{E_{\mathbf i}}E_{\mathbf j})_{E_{\mathbf k}'}=0$ if $E_{\mathbf k}$ is reducible
(and $E_{\mathbf i}$, $E_{\mathbf j}$ are irreducible). 
\item[(iii)] $\Gamma_{E_{\mathbf i}''}E_{\mathbf j}=\Gamma_{E_{\mathbf i}}E_{\mathbf j}''=0$ if $E_{\mathbf i}$ and 
$E_{\mathbf j}$ are reducible. 
\end{enumerate}
\end{lem}

\Pf The corresponding slice is homogeneous and thus 
congruent to a principal orbit of the isotropy representation 
of a symmetric space. The roots of the symmetric space are,
up to sign, in a natural bijection with the focal lines
and therefore with $v_{\mathbf i}$, $v_{\mathbf j}$, $v_{\mathbf k}$, $v_{\mathbf l}$. We may assume 
that the root system is of the form $\{\pm\theta_1,\pm\theta_2,
\pm(\theta_1\pm\theta_2)\}$ in the $B_2$ case and 
 $\{\pm\theta_1,\pm\theta_2,\pm2\theta_1,\pm2\theta_2,
\pm(\theta_1\pm\theta_2)\}$ in the $(BC)_2$ case.
Now (i)-(iii) follow from the discussion in subsection~\ref{fin-dim}
by using the bracket relations
\begin{enumerate}
\item[(i')] $[\Lk_{\theta_1+\theta_2},\Lp_{\theta_1-\theta_2}]
=[\Lk_{\theta_1-\theta_2},\Lp_{\theta_1+\theta_2}]=0$ (case $B_2$),
\item[(ii')] $[\Lk_{\theta_1+\theta_2},\Lp_{\theta_1-\theta_2}]
=[\Lk_{\theta_1-\theta_2},\Lp_{\theta_1+\theta_2}]
\subset\Lp_{2\theta_1}+\Lp_{2\theta_2}$ (case $BC_2$), 
\item[(iii')] $[\Lk_{2\theta_1},\Lp_{\theta_2}+\Lp_{2\theta_2}]
=[\Lk_{2\theta_2},\Lp_{\theta_1}+\Lp_{2\theta_1}]=0$ (case $BC_2$),
\end{enumerate}
respectively. \EPf

\medskip

\begin{rem}\label{refl}
Let $W$ be an affine Weyl group isomorphic to $\tilde C_2$ 
acting on an Euclidean plane, which is generated by reflections 
on a family of lines $\mathcal H$. We call a point
an \emph{intersection point} if it lies on at least
two different lines in $\mathcal H$. 
It is clear that along any reflection line
the intersection points are equally spaced, and there are 
exactly two possibilities for their spacing, one being wider
than the other. Let $H$ and 
$H'$ be two nonparallel lines in $\mathcal H$. 
\begin{enumerate}
\item[(i)] If the spacing of intersection points along $H$
is wide (or $H\perp H'$) then there passes
through each intersection point on $H$ a line in $\mathcal H$ 
which is parallel to $H'$. 
\item[(ii)] In general, there passes at least through each second intersection
point on $H$ a line in $\mathcal H$ 
which is parallel to $H'$.  
\end{enumerate}
\end{rem}

\begin{prop}\label{gamma-image-1}
Let $i$, $j\in\Z$ with $i-j$ even and $m=\frac{i+j}2$. 
Let $\alpha\in\mathcal A$ with $E_{\alpha,i}$ irreducible.
Then
\begin{equation}\label{gamma-mid-zero}
 (\Gamma_{E_{\alpha,i}}E_{\alpha,j})_{E_{\alpha,m}}=0 
\end{equation}
unless the diagram of $M$ is 
\[ \Cnta{}{}{}{}{}{} \quad\mbox{or}\quad \Cntar{}{}{}{}{}{} \]
and $(\alpha,i)$ corresponds to the right extremal vertex while
$(\alpha,m)$ corresponds to the left one. In particular
equation~(\ref{gamma-mid-zero}) holds if $i-j$ is divisible by~$4$ or
the affine Weyl group of $M$ 
is not of type $\tilde C_n$ ($n\geq2$). 
\end{prop}

\Pf If $(\alpha,i)$ and $(\alpha,m)$ correspond to the 
two extremal vertices of a $\tilde C_n$ diagram then 
$i-m$ is odd necessarily, that is, $i-j$ is not 
divisible by~$4$. Hence it suffices to prove the first 
assertion, i.e.~that $(\Gamma_{E_{\alpha,i}}E_{\alpha,j})_{E_{\alpha,m}}
\neq0$ can only occur in the two special cases described above. 

By Theorem~\ref{image1} and Lemma~\ref{rk2} we may assume 
$M$ to be of type $\tilde C_2$. Note that according to
Lemma~\ref{rk2} the arrows in the diagram of a slice of type
$\tilde C_2$ both point inward or both point outward 
if $M$ is of type $\tilde B_n$ and do not change direction 
if $M$ is of type $\tilde C_n$. 

So let $M$ be of type $\tilde C_2$ and assume that 
$(\Gamma_{E_{\alpha,i}}E_{\alpha,j})_{E_{\alpha,m}}\neq0$. Exactly
as in the proof of Theorem~\ref{gamma-image-general}
we find $\beta$, $\gamma$, $\delta\in\mathcal A$ 
and $k$, $\ell$, $n\in\Z$ with $\beta\neq\gamma$ and such that 
$(\Gamma_{E_{\beta,k}}E_{\gamma,\ell})_{E_{\alpha,m}}$,
$(\Gamma_{E_{\beta,k}}E_{\alpha,i})_{E_{\delta,n}}$ and 
$(\Gamma_{E_{\gamma,\ell}}E_{\alpha,j})_{E_{\delta,n}}$ are all nonzero. 
Lemma~\ref{euclid} yields that the focal lines $H_{\beta,k}$
and $H_{\gamma,\ell}$ are orthogonal to each other and 
make an angle $\pi/4$ with 
$H_{\alpha,m}$, and that $H_{\delta,n}$ is
orthogonal to $H_{\alpha,i}$. Let $x_i\in H_{\alpha,i}$
and $x_m\in H_{\alpha,m}$ be the intersection points of
$H_{\delta,n}$, $H_{\beta,k}$ and 
of $H_{\beta,k}$, $H_{\gamma,\ell}$, respectively. 

We first observe that $E_{\beta,k}$ and $E_{\gamma,\ell}$
are necessarily irreducible. In fact, if one of them were reducible, then 
also the other one would be, since reflection at $H_{\alpha,m}$
maps $H_{\beta,k}$ to $H_{\gamma,\ell}$. 
By Lemma~\ref{b2bc2}(iii) we have $\Gamma_{E_{\beta,k}}E_{\gamma,\ell}=
\Gamma_{E'_{\beta,k}}E'_{\gamma,\ell}$ and $E_{\beta,k}$, $E_{\gamma,\ell}$ 
could be replaced in our argument by their components 
$E'_{\beta,k}$, $E'_{\gamma,\ell}$ as they have only to satisfy
$\langle(\Gamma_{E_{\alpha,i}}E_{\alpha,j})_{E_{\alpha,m}},
\Gamma_{E_{\beta,k}}E_{\gamma,\ell}\rangle\neq0$. 
However, $(\Gamma_{E'_{\beta,k}}E_{\alpha,i})_{E_{\delta,n}}\neq0$
would then be in contradiction to Lemma~\ref{b2bc2}(ii). 

If $i-m$ is even or $(\alpha,i)$ corresponds to the vertex in 
the middle of the diagram of $M$ (of $\tilde C_2$-type) 
then there exists due to Remark~\ref{refl}
an integer $\ell_1$ such that $H_{\gamma,\ell_1}$ passes
through the intersection point of $H_{\alpha,m}$ and 
$H_{\delta,n}$. Reflection at $H_{\gamma,\ell_1}$ 
maps~$x_m$ to~$x_i$, $H_{\alpha,m}$ to $H_{\delta,n}$,
$H_{\gamma,\ell}$ to a parallel focal line 
$H_{\gamma,\ell_2}$ through~$x_i$, and preserves 
$H_{\beta,k}$. This implies 
$\Gamma_{E_{\beta,k}}E_{\gamma,\ell_2}\neq0$ in contradiction to 
$\Gamma_{E_{\alpha,i}}E_{\delta,n}\neq0$ and Lemma~\ref{b2bc2}(i). 
Note that $E_{\alpha,i}$, $E_{\beta,k}$, $E_{\delta,n}$ and
$E_{\gamma,\ell_2}$ all belong to the slice centered
at $x_i$, which is of type $B_2$ as $E_{\alpha,i}$ is irreducible
by assumption and $E_{\beta,k}$ is irreducible by the above. 
Thus $i-m$ is odd and $(\alpha,i)$ corresponds to one
of the extremal vertices of the diagram of $M$. 
Since $i-m$ is odd, 
$(\alpha,m)$ corresponds to the other extremal vertex of the 
diagram. We also see that the root corresponding to 
$(\alpha,i)$ in the slice centered at $x_i$ must be short
by Lemma~\ref{b2bc2} as $\Gamma_{E_{\alpha,i}}E_{\delta,n}\neq0$. 
Hence an arrow points to the vertex corresponding to $(\alpha,i)$. 

Finally, consider $E_{\alpha,m}$. If it is irreducible, and thus 
the slice centered at $x_m$ is also of type~$B_2$, then the root
corresponding to $(\alpha,m)$ must be long again by Lemma~\ref{b2bc2},
as $\Gamma_{E_{\beta,k}}E_{\gamma,\ell}\neq0$. 
Hence in this case the diagram of $M$ is
$\Ciita{\!\!\!\!\!\!(\alpha,m)}{}{\!\!\!\!\!\!(\alpha,i)}$
where $(\alpha,m)$ and $(\alpha,i)$ correspond to the 
vertices indicated. On the other hand, if $E_{\alpha,m}$ is reducible,
the diagram necessarily is 
$\Ciitar{\!\!\!\!\!\!(\alpha,m)}{}{\!\!\!\!\!\!(\alpha,i)}$
as the arrow between two vertices always points to that one which 
corresponds to a reducible eigenspace, if such a vertex occurs. \EPf

\medskip

The next theorem is not necessary for 
the proof of the continuity of $\Gamma$, but it 
contains interesting information
that sharpens Theorem~\ref{image1}.

\begin{thm}\label{image2}
Assume $E_{\alpha,i}$ is irreducible. Then
\begin{enumerate}
\item[(i)] \[\Gamma_{E_{\alpha,i}}E_{\alpha,j}\subset E_{\mathbf 0}\]
if $i-j$ is divisible by $4$ or the Weyl group of~$M$ 
is not of type~$\tilde C_n$ ($n\geq2$). 
\item[(ii)] \[\Gamma_{E_{\alpha,i}}E_{\alpha,j}\subset E_{\mathbf 0}+
E_{\alpha,\frac{i+j}2}\] if $i-j$ is even. 
\end{enumerate}
\end{thm}

\Pf If $i-j$ is even or $W$ is not of type $\tilde C_n$ then 
$(\Gamma_{E_{\alpha,2i-j}}E_{\alpha,j})_{E_{\alpha,i}}
=(\Gamma_{E_{\alpha,i}}E_{\alpha,2j-i})_{E_{\alpha,j}}=0$ 
by Proposition~\ref{gamma-image-1} (note that 
$(2i-j)-j=2(i-j)$ is divisible by $4$ if $i-j$ is even,
and that $E_{2i-j}$ is irreducible) and hence 
$(\Gamma_{E_{\alpha,i}}E_{\alpha,j})_{E_{\alpha,2i-j}}=
(\Gamma_{E_{\alpha,i}}E_{\alpha,j})_{E_{\alpha,2j-i}}=0$
by Lemma~\ref{gamma-colinear}. Thus~(i) and~(ii) follow from 
Theorem~\ref{gamma-image-general} and Proposition~\ref{gamma-image-1}. \EPf

\begin{prop}\label{gamma-image-2}
Let $i$, $j\in\Z$ with $i-j$ even, $i\not=j$, and $m=\frac{i+j}2$.
Let $\alpha\in\mathcal A$ with $E_{\alpha,i}$ reducible.
Then
\begin{enumerate}
\item[(i)] 
$(\Gamma_{E_{\alpha,i}''}E_{\alpha,j})_{E_{\alpha,m}}=0$.
\item[(ii)] 
$E_{\alpha,m}$ is reducible and 
$(\Gamma_{E_{\alpha,i}}E_{\alpha,j})_{E'_{\alpha,m}}=0$
if $i-j$ is divisible by $4$.
\end{enumerate}
\end{prop}

\Pf In both cases we follow essentially the proof of 
Proposition~\ref{gamma-image-1}.

(i) Suppose, to the contrary, that 
$(\Gamma_{E_{\alpha,i}''}E_{\alpha,j})_{E_{\alpha,m}}\neq0$. 
Then we can find  $\beta$, $\gamma$, $\delta\in \mathcal A$ 
and $k$, $\ell$, $n\in\Z$ 
with $(\Gamma_{E_{\beta,k}}E_{\gamma,\ell})_{E_{\alpha,m}}\neq0$,
$(\Gamma_{E_{\alpha,i}''}E_{\beta,k})_{E_{\delta,n}}\neq0$ and
$(\Gamma_{E_{\alpha,j}}E_{\gamma,\ell})_{E_{\delta,n}}\neq0$,
which implies that 
$H_{\beta,k}$ and $H_{\gamma,\ell}$ 
are orthogonal, each of them makes an angle of $\pi/4$ 
with $H_{\alpha,m}$, and
$H_{\delta,n}$ is orthogonal to $H_{\alpha,i}$.
In particular $E_{\delta,n}$ is also reducible as the reflection
on $H_{\beta,k}$ is an element of the affine Weyl group 
that maps $H_{\alpha,i}$ to $H_{\delta,n}$.
However 
$(\Gamma_{E_{\alpha,i}''}E_{\beta,k})_{E_{\delta,n}}\neq0$ 
yields $\Gamma_{E_{\alpha,i}''}E_{\delta,n}\neq0$ 
and this contradicts
Lemma~\ref{b2bc2}(iii).

(ii) If $i-j$ is divisible by $4$, then $i-m$ is divisible by $2$, 
so there is an element in the affine Weyl group
that maps $H_{\alpha,i}$ to $H_{\alpha,m}$ 
and this shows that $E_{\alpha,m}$ is reducible.  

Suppose now that 
$(\Gamma_{E_{\alpha,i}}E_{\alpha,j})_{E'_{\alpha,m}}\neq0$.
Then we can find  $\beta$, $\gamma$, $\delta\in\mathcal A$ 
and $k$, $\ell$, $n\in\Z$ 
with $(\Gamma_{E_{\beta,k}}E_{\gamma,\ell})_{E'_{\alpha,m}}\neq0$,
$(\Gamma_{E_{\alpha,i}}E_{\beta,k})_{E_{\delta,n}}\neq0$ and
$(\Gamma_{E_{\alpha,j}}E_{\gamma,\ell})_{E_{\delta,n}}\neq0$.
However the first inequality is a contradiction to 
Lemma~\ref{b2bc2}(ii) as $H_{\beta,k}$ and
$H_{\gamma,\ell}$ make an angle of $\pi/4$ with 
$H_{\alpha,m}$ and thus $E_{\beta,k}$, $E_{\gamma,\ell}$
are irreducible. Note that in a finite dimensional rank $2$
slice of type $(BC)_2$ exactly one pair of orthogonal focal lines
corresponds to irreducible eigenspaces.\hfill\mbox{ } \EPf

\medskip

Using the same idea as in the proof of Theorem~\ref{image2}, one gets from 
Theorem~\ref{gamma-image-general} and Proposition~\ref{gamma-image-2} the 
statements (i)-(iv) of Theorem~E. Thus Theorems~D and~E follow 
from Theorem~\ref{gamma-image-general} and the last two results.

The following result is a simple application of Theorems~D and~E.

\begin{thm}
If $M$ is infinite dimensional so is $E_{\mathbf 0}$. 
\end{thm}

\Pf Let $\alpha\in\mathcal A$. Then 
Corollary~\ref{gamma-eiej-eiek-perp-2} implies that 
for all $j\in\mathbb Z\setminus\{0\}$,
$\Gamma_{E_{\alpha,0}}E_{\alpha,j}\perp 
\Gamma_{E_{\alpha,0}}E_{\alpha,k}$ for all but finitely many $k\in\mathbb Z$. 
Starting with $j_1=4$ we may thus construct a monotone sequence~$(j_k)_k\geq1$
in $4\mathbb Z$ such that 
$\Gamma_{E_{\alpha,0}}E_{\alpha,j_k}\perp 
\Gamma_{E_{\alpha,0}}E_{\alpha,j_\ell}$ for all $k\neq\ell$. 
If $E_{\alpha,0}$ is irreducible (resp.~reducible)
so are the $E_{\alpha,j_k}$, and 
$\Gamma_{E_{\alpha,0}}E_{\alpha,j_k}\subset E_{\mathbf 0}$ 
(resp.~$\Gamma_{E_{\alpha,0}''}E_{\alpha,j_k}''\subset E_{\mathbf 0}$)
by Theorem~D (resp.~Theorem~E). Since those subspaces are pairwise
orthogonal and never zero due to Corollary~\ref{gamma-cartan},
the result follows. \EPf

\section{Continuity of $\Gamma$ and rigidity}
\setcounter{equation}{0}

In this section, we collect results from previous 
sections to prove the continuity of $\Gamma$ in complete 
generality (Theorem~A). The rigidity theorem (Theorem~B) is 
then a consequence.  

\begin{prop}\label{gamma-estimate}
Let $\alpha\in\mathcal A$, let $i$, $j\in\Z$ with $i\neq j$,
and let $X_{\alpha,i}\in E_{\alpha,i}$, $Y_{\alpha,j}\in E_{\alpha,j}$.
Then 
\begin{enumerate}
\item[(i)] 
\[ ||\Gamma_{X_{\alpha,i}}Y_{\alpha,j}||^2\leq ||v_{\alpha,i}||^2||X_{\alpha,i}||^2||Y_{\alpha,j}||^2+ 3||(\Gamma_{X_{\alpha,i}}Y_{\alpha_{\mathbf j}})_{E_{\alpha,
\frac{i+j}2}}||^2, \]
where the last term is to be omitted in case $i-j$ is odd. 
\item[(ii)] If $i-j$ is not divisible by $2^k$ for some integer $k\geq1$ then 
\[ ||\Gamma_{X_{\alpha,i}}Y_{\alpha,j}||\leq2^{k-1}||v_{\alpha,i}||||X_{\alpha,i}||||Y_{\alpha,j}||. \]
\end{enumerate}
\end{prop}

\Pf We combine Proposition~\ref{gamma-rk1} together with
Theorem~\ref{gamma-image-general} to write
\begin{eqnarray*}
\frac12||v_{\alpha,i}||^2||X_{\alpha,i}||^2||Y_{\alpha,j}||^2 
& = & ||(\Gamma_{X_\alpha,i}Y_{\alpha,j})_{E_{\mathbf 0}}||^2 
+ 2||(\Gamma_{X_\alpha,i}Y_{\alpha,j})_{E_{\alpha,2i-j}}||^2 \\
&&\qquad+\frac12||(\Gamma_{X_\alpha,i}Y_{\alpha,j})_{E_{\alpha,2j-i}}||^2
-||(\Gamma_{X_\alpha,i}Y_{\alpha,j})_{E_{\alpha,\frac{i+j}2}}||^2.
\end{eqnarray*}
Multiplying through by $2$ and adding 
$3||(\Gamma_{X_{\alpha,i}}Y_{\alpha,j})_{E_{\alpha,\frac{i+j}2}}||^2$ 
to both sides yields~(i). 

In the proof of~(ii) we use induction on $k$.
The case $k=1$ is contained in (i). Now we assume that 
(ii) holds for some $k\geq1$ and that $i-j$ is not divisible by $2^{k+1}$.
Then $i-\frac{i+j}2$ is not divisible by $2^k$. Therefore,
for $Z=(\Gamma_{X_{\alpha,i}}Y_{\alpha,j})_{E_{\alpha,\frac{i+j}2}}$, we get
\begin{eqnarray*}
||(\Gamma_{X_{\alpha,i}}Y_{\alpha,j})_{E_{\alpha,\frac{i+j}2}}||^2 & = & 
\langle\Gamma_{X_{\alpha,i}}Y_{\alpha,j},Z\rangle \\
& = & |\langle Y_{\alpha,j},\Gamma_{X_{\alpha,i}}Z\rangle | \\
& \leq & ||Y_{\alpha,j}||\,||\Gamma_{X_{\alpha,i}}Z|| \\
& \leq & 2^{k-1}\,||Y_{\alpha,j}||\,||v_{\alpha,i}||\,||X_{\alpha,i}||\,||Z||,
\end{eqnarray*}
by the Cauchy-Schwarz inequality and the induction hypothesis, and thus 
\[ ||(\Gamma_{X_{\alpha,i}}Y_{\alpha,j})_{E_{\alpha,\frac{i+j}2}}||\leq
2^{k-1}\,||v_{\alpha,i}||\,||X_{\alpha,i}||\,||Y_{\alpha,j}||. \]
The inequality in~(ii) now follows from~(i). \EPf

\medskip

We finally come to one of our main results.

\begin{thm}\label{gamma-continuity}
For all ${\mathbf i}\in \mathbf I^*$ and $X\in E_{\mathbf i}$,
$\Gamma_X$ is continuous. More precisely, 
there exists a constant $C$ such that 
$||\Gamma_X Y|| \leq C\,||v_{\mathbf i}||\, ||X||\, ||Y||$
for all $\mathbf i\in\mathbf I^*$, $X\in E_{\mathbf i}$ and $Y\in T_xM$. 
In particular, the one-parameter 
groups $F_X^t$ are smooth curves in the Banach-Lie group
of isometries of~$V$.
\end{thm}

\Pf Fix $(\alpha,i)\in \mathcal A\times\Z$ and $X_{\alpha,i}\in E_{\alpha,i}$.
For the continuity of $\Gamma_{X_{\alpha,i}}$, 
in view of the discussion
preceeding Lemma~\ref{finite-non-orthog}, 
Corollary~\ref{gamma-eiej-eiek-perp-2} 
and Proposition~\ref{gamma-estimate},
it is enough to show that 
\begin{equation}\label{gamma-estimate-final}
 ||(\Gamma_{X_{\alpha,i}}Y_{\alpha,j})_{E_{\alpha,m}}||
\leq 2\,||v_{\alpha,i}||
\,||X_{\alpha,i}||\,||Y_{\alpha,j}|| 
\end{equation}
for all $j\in\Z$ with $j\neq i$ such that $j-i$ is divisible
by $4$ and all $Y_{\alpha,j}\in E_{\alpha,j}$, where $m=\frac{i+j}2$.

If $E_{\alpha,i}$ is irreducible or $E_{\alpha,i}$ is reducible
and $X_{\alpha,i}\in E''_{\alpha,i}$, the left hand side 
of~(\ref{gamma-estimate-final}) is zero (Theorems~\ref{gamma-image-1}
and~\ref{gamma-image-2}(i)). Thus we may assume $E_{\alpha,i}$ reducible
and $X_{\alpha,i}\in E'_{\alpha,i}$. 

Let $Z=(\Gamma_{X_{\alpha,i}}Y_{\alpha,j})_{E_{\alpha,m}}$.
Due to Proposition~\ref{gamma-image-2}(ii), $Z\in E''_{\alpha,m}$. Thus
\begin{eqnarray*}
||(\Gamma_{X_{\alpha,i}}Y_{\alpha,j})_{E_{\alpha,m}}||^2 & = & 
\langle \Gamma_{X_{\alpha,i}}Y_{\alpha,j},Z\rangle \\
&=& |\langle Y_{\alpha,j},(\Gamma_{X_{\alpha,i}}Z)_{E_{\alpha,j}}\rangle| \\
&\leq& ||Y_{\alpha,j}||\,||(\Gamma_{X_{\alpha,i}}Z)_{E_{\alpha,j}}|| \\
& = & ||Y_{\alpha,j}||\,||(\Gamma_ZX_{\alpha,i})_{E_{\alpha,j}}||
\,\left|\frac{v_{\alpha,i}-v_{\alpha,j}}{v_{\alpha,m}-v_{\alpha,j}}\right| \\
&\leq&||Y_{\alpha,j}||\,||v_{\alpha,m}||\,||Z||\,||X_{\alpha,i}||\,
\left|\frac{v_{\alpha,i}-v_{\alpha,j}}{v_{\alpha,m}-v_{\alpha,j}}\right|,
\end{eqnarray*}
where we have used Cauchy-Schwarz, Codazzi (Proposition~\ref{gamma-permute})
and Proposition~\ref{gamma-estimate}(i) 
(note that $(\Gamma_ZX_{\alpha,i})_{E_{\alpha,\frac{m+i}2}}=0$ 
by Proposition~\ref{gamma-image-2}(i)).

Recall that $v_{\alpha,k}=\frac{a}{b+k}v_0$ for a unit vector
$v_0$ and $a$, $b\in\R$. This implies that 
\[ ||v_{\alpha,m}||\,
\left|\frac{v_{\alpha,i}-v_{\alpha,j}}{v_{\alpha,m}-v_{\alpha,j}}\right|
=2||v_{\alpha,i}|| \]
and hence
\[ ||(\Gamma_{X_{\alpha,i}}Y_{\alpha,j})_{E_{\alpha,m}}||\leq 2
||v_{\alpha,i}||\,||X_{\alpha,i}||\,||Y_{\alpha,j}||, \]
as we wished. 

Our discussion so far shows 
that there exists a constant $C$ such that 
\[ ||\Gamma_X Y|| \leq C\,||v_{\mathbf i}||\, ||X||\, ||Y|| \]
for all $\mathbf i\in\mathbf I^*$, $X\in E_{\mathbf i}$ and 
$Y\in E_{\mathbf 0}^\perp$.  Assume now that $Y\in E_{\mathbf 0}$. 
Since $\Gamma_{X_{\mathbf i}}E_{\mathbf 0}\perp E_{\mathbf 0}$, 
we can find a sequence $Z_n\in \sum_{\mathbf j\in \mathbf I^*}E_{\mathbf j}$
such that $Z_n\to\Gamma_{X_{\mathbf i}}Y$. Then
\[ \langle\Gamma_{X_{\mathbf i}}Y,Z_n\rangle = 
-\langle Y,\Gamma_{X_{\mathbf i}}Z_n\rangle \leq 
||Y||\,||\Gamma_{X_{\mathbf i}}Z_n||\leq
C\,||v_{\mathbf i}||\,||X_{\mathbf i}||\,||Y||\,||Z_n||, \]
which yields in the limit as $n\to\infty$ the desired inequality. \EPf

\begin{cor}
$\Gamma$ is continuous as a bilinear
mapping, that is, there exists $C>0$ such that
$||\Gamma_X Y|| \leq C\, ||X||\, ||Y||$ for
all $X$, $Y\in T_xM$ with $X \perp E_{\mathbf 0}$.
\end{cor}

\Pf If $X=\sum_{\mathbf i\in\mathbf I^*}X_{\mathbf i}$ with 
$X_{\mathbf i}\in E_{\mathbf i}$ and 
$Y\in\sum_{\mathbf i\in\mathbf I}E_{\mathbf i}$ then
\begin{eqnarray*}
||\Gamma_XY||&\leq& \sum_{\mathbf i\in\mathbf I^*}||\Gamma_{X_{\mathbf i}}Y|| \\
&\leq& \sum_{\mathbf i\in\mathbf I^*}C\,||v_{\mathbf i}||\,||X_{\mathbf i}||
\,||Y||\\
&\leq& C\,||Y||\left(\sum_{\mathbf i\in\mathbf I^*}||v_{\mathbf i}||^2\sum_{\mathbf i\in\mathbf I^*}||X_{\mathbf i}||^2\right)^{1/2}\\
&\leq&C'\,||X||\,||Y||,
\end{eqnarray*}
where $C'= C\left(\sum_{\mathbf i\in\mathbf I^*}||v_{\mathbf i}||^2\right)^{1/2}<+\infty$. \EPf

\medskip

The continuity of $\Gamma$ is essential in the proof
of the following rigidity theorem.

\begin{thm}\label{rigidity}
For any point $x\in M$, $\alpha_x$ and $(\nabla\alpha)_x$ 
determine $M$ completely. 
\end{thm}

\Pf Let $\tilde M$ be a second connected complete full 
irreducible isoparametric submanifold of $V$ with 
$x\in M\cap\tilde M$, $T_xM=T_p\tilde M$,
$\alpha_x=\tilde\alpha_x$ and $(\nabla\alpha)_x=(\nabla\tilde\alpha)_x$.
Owing to Theorem~\ref{thm:nabla-alpha-gamma},  
$\Gamma_x=\tilde\Gamma_x$. 

It follows from $T_xM=T_x\tilde M$ and $\alpha_x=\tilde\alpha_x$
that $M$ and $\tilde M$ have the
same normal spaces and the same curvature spheres at $x$.
Moreover, for each ${\mathbf i}\in \mathbf I^*$ and $X\in E_{\mathbf i}(x)$, the one-parameter
groups $F_X^t$ and $\tilde F_X^t$ coincide since they have 
the same infinitesimal generators, defined on the whole of $V$
by Theorem~\ref{gamma-continuity} (notice that it is at this point that the 
continuity of $\Gamma_X$ and $\tilde\Gamma_X$ is crucial since otherwise 
the self-adjoint infinitesimal generators of $(F_X^t)_*$ and $(\tilde F_X^t)_*$
might not coincide). 
Therefore for all $y\in  S_{\mathbf i}(x)$, we have by equivariance
that $y\in M\cap\tilde M$, $T_yM=T_y\tilde M$,
$\alpha_y=\tilde\alpha_y$ and $\Gamma_y=\tilde\Gamma_y$. 
Proceeding by induction we now see
that $M$ and $\tilde M$ coincide along $Q_x=\tilde Q_x$.
Since $Q_x$ and $\tilde Q_x$ are dense in $M$ and 
$\tilde M$, respectively, the desired result follows. \EPf

\medskip

If $M$ is finite dimensional
and homogeneous, the continuity of $\Gamma$ is of course obvious
and the proof of Theorem~\ref{rigidity} 
also applies. Even in this case
the result seems to be new.



\providecommand{\bysame}{\leavevmode\hbox to3em{\hrulefill}\thinspace}
\providecommand{\MR}{\relax\ifhmode\unskip\space\fi MR }
\providecommand{\MRhref}[2]{%
  \href{http://www.ams.org/mathscinet-getitem?mr=#1}{#2}
}
\providecommand{\href}[2]{#2}

\end{document}